\documentclass[review, 12pt]{elsarticle}
\usepackage[dvipsnames]{xcolor}
\usepackage{setspace}
\usepackage{multirow}
\usepackage{multicol}
\usepackage{graphicx}
\usepackage{xcolor}
\RequirePackage{amsmath}
\usepackage{bm}
\usepackage{xfrac}
\usepackage{csquotes}
\usepackage{float}
\usepackage{arydshln} 
\usepackage{subfigure}
\usepackage{amsmath,esint,amsfonts}

\usepackage[most]{tcolorbox}
\usepackage[colorinlistoftodos]{todonotes}
\usepackage{lineno}
\usepackage[colorlinks=true,linkcolor=blue,urlcolor=blue,citecolor=blue,anchorcolor=blue]{hyperref}
\usepackage{cleveref} 
\usepackage{nomencl}
\makenomenclature
\usepackage{etoolbox}
\renewcommand\nomgroup[1]{%
	\item[\bfseries
	\ifstrequal{#1}{X}{Symbols}{%
		\ifstrequal{#1}{Y}{Greek Symbols}{%
			{}}}%
	]}
\nomenclature[X, 08]{$\dot{q}$}{Rate of heat generation}%
\nomenclature[X, 09]{$r$}{Euclidean distance between points (radius)}%
\nomenclature[X, 03]{$d$}{Dimension (1D, 2D or 3D)}%
\nomenclature[X, 10]{$s$}{Interpolated variable}%
\nomenclature[X, 11]{$\bm{s}$}{Set of interpolated variables at the discrete points}%
\nomenclature[X, 01]{$\bm{A}$}{Coefficient matrix}%
\nomenclature[X, 04]{$k$}{Thermal conductivity}%
\nomenclature[X, 12]{$T$}{Temperature}%
\nomenclature[X, 02]{$C$}{Order of convergence}%
\nomenclature[X, 06]{$P$}{Polynomial appended to PHS-RBF}%
\nomenclature[X, 13]{$X$,$Y$,$Z$}{Cartesian coordinates}%
\nomenclature[Y, 04]{$\phi$}{Radial basis function (PHS-RBF)}%
\nomenclature[Y, 05]{$\bm{\Phi}$}{Submatrix of PHS-RBF (cloud)}%
\nomenclature[X, 05]{$p$}{Degree of appended polynomial}%
\nomenclature[X, 07]{$\bm{P}$}{Submatrix composed of monomials appended to PHS}%
\nomenclature[Y, 03]{$\lambda_{i}$,$\gamma_{i}$}{Unknown interpolation coefficients}%
\nomenclature[Y, 03]{$\Delta$r}{Average grid spacing}%
 \usepackage{tikz}
 \usetikzlibrary{shapes, arrows}
 \usepackage{ragged2e}

\journal{Journal}

\biboptions{sort&compress}

\begin{document}

\begin{frontmatter}

\title{Simulation of Heat Conduction in Complex Domains of Multi-material Composites using a Meshless Method}

\author{Naman Bartwal\textsuperscript{a}}
\author{Shantanu Shahane\textsuperscript{b}}
\author{Somnath Roy\textsuperscript{a}\fnref{Corresponding Author}}
\author{Surya Pratap Vanka\textsuperscript{c}}
\address{\textsuperscript{a}Department of Mechanical Engineering\\
	Indian Institute of Technology Kharagpur \\
	Kharagpur, West Bengal 721302}
\address{\textsuperscript{b}Google, 1600 Amphitheatre Parkway\\
	Mountain View, California 94043}
\address{\textsuperscript{c}Department of Mechanical Science and Engineering\\
	University of Illinois at Urbana-Champaign \\
	Urbana, Illinois 61801}
\fntext[Corresponding Author]{Corresponding Author Email: \url{somnath.roy@mech.iitkgp.ac.in}}
\vspace{-0.8cm}
\begin{abstract}
Several engineering applications involve complex materials with significant and discontinuous variations in thermophysical properties. These include materials for thermal storage, biological tissues with blood capillaries, etc. For such applications, numerical simulations must exercise care in not smearing the interfaces by interpolating variables across the interfaces of the subdomains. In this paper, we describe a high accuracy meshless method that uses domain decomposition and cloud-based interpolation of scattered data to solve the heat conduction equation in such situations. The polyharmonic spline (PHS) function with appended polynomial of prescribed degree is used for discretization. A flux balance condition is satisfied at the interface points and the clouds of interpolation points are restricted to be within respective domains. Compared with previously proposed meshless algorithms with domain decomposition, the cloud-based interpolations are numerically better conditioned, and achieve high accuracy through the appended polynomial. The accuracy of the algorithm is demonstrated in several two and three dimensional problems using manufactured solutions to the heat conduction equation with sharp discontinuity in thermal conductivity. Subsequently, we demonstrate the applicability of the algorithm to solve heat conduction in complex domains with practical boundary conditions and internal heat generation. Systematic computations with varying conductivity ratios, interpoint spacing and degree of appended polynomial are performed to investigate the accuracy of the algorithm. 
\end{abstract}


\end{frontmatter}

\section{Introduction}
In several heat transfer engineering applications, thermal properties such as conductivity, specific heat and density vary discontinuously and significantly across regions of complex shape. Some examples are composite materials, metal foams with embedded polymers for thermal storage and biological tissues embedded with blood vessels \cite{ali2019critical,aramesh2022metal,bhowmik2013conventional}. In some applications, the thermal properties can vary by as much as a factor of one hundred between adjacent regions. In such situations, while seeking numerical solutions, it is necessary to discretize the governing equations carefully to avoid smearing of the discontinuities as well as prevent non-monotonic variation of the thermal variables. Such regions can be identified and discretized individually as separate domains with local properties within each region. The unstructured finite volume method is a popular numerical technique used for such problems. However, its accuracy is typically limited to second order unless complex reconstruction schemes or higher order basis functions are used \cite{wang2016compact,wang2016compact2,liu2016high}. Further, reconstruction schemes and high order basis functions must not interpolate across the discontinuous property interfaces. In addition, the grids$\slash$volumes generated must not be excessively skewed and the aspect ratio of the elements must be close to unity to avoid loss of accuracy. Achieving high accuracy with finite volume method can be challenging in problems with large and discontinuous variations in properties. 
In this paper, we describe a high accuracy meshless numerical technique which uses polyharmonic radial basis functions (PHS-RBF) with appended polynomial to interpolate scattered data and solve elliptic partial differential equations encountered in multiphysics simulations. In our recent works \cite{ shahane2021high, shahane2021semi, bartwal2022application, shahane2022consistency, unnikrishnan2022shear}, we have demonstrated the technique for heat conduction and fluid flow problems in several complex geometries and have shown high discretization accuracy whose order increases with the degree of the appended polynomial. Further, we have used a cloud based interpolation, versus global interpolation, to be able to use large numbers of scattered points. The method can be used to simulate thermal transport in large and complex domains with simultaneous refinement in point spacing as well as polynomial degree ($h-p$ refinement). However, when this method is used in a straightforward way for problems with large and discontinuous property variations, significant errors can result because of interpolation across the discontinuities. Therefore, in this work we have extended the algorithm to multiple subdomains with interpolations strictly within individual homogeneous domains. At the interfaces, flux balance conditions are imposed for the discrete points on the interfaces. In this paper, we investigate the accuracy resulting from such $C_1$ condition at the interface by considering manufactured (known) solutions. We then apply the method to several heat conduction problems in complex geometries. 

Meshless methods offer several advantages over mesh based methods and have been pursued since the early 1990s. Meshless methods based on radial basis function (RBF) interpolations to solve partial differential equations were originally demonstrated by Kansa \cite{kansa1990multiquadrics,kansa1991multiquadrics} who used the Hardy multiquadrics \cite{hardy1971multiquadric} to interpolate scattered data. Kansa\textquotesingle s original work used global interpolation which is exponentially accurate but suffers from high condition numbers for modest numbers of scattered points. Subsequently, domain decomposition methods  have been proposed to reduce the condition number with coupling between the subdomains through overlapping or interface flux balances \cite{mayo1984fast,mayo1985fast,leveque1994immersed,chen1998finite,li1999numerical,liu2000boundary,gibou2002second,li2003new,hou2005numerical}. However, such methods also are not suitable for large problems with need for millions of scattered points to represent a complex domain. In contrast, interpolation within clouds of scattered points maintains the condition number within solvability limits and produces good accuracy. Further, multiquadrics, inverse multiquadrics and Gaussians all require prescription of a shape parameter \cite{shu2003local,larsson2003numerical,ding2006numerical,sanyasiraju2008local,chandhini2007local,sanyasiraju2009note,vidal2016direct,zamolo2019solution,mairhuber1956haar,buhmann1993spectral,wong2002compactly,fornberg2004stable,fornberg2004some,larsson2005theoretical,flyer2009radial,fornberg2011stable,fasshauer2012stable,fornberg2013stable,flyer2016enhancing} which dictates interpolation accuracy and the condition number of the interpolation coefficient matrix. On the other hand, interpolation with polyharmonic splines (PHS-RBF) does not need a shape parameter, and when appended with polynomial produces interpolation accuracy consistent with the degree of the appended polynomial \cite{barnett2015robust,flyer2016role,bayona2017role,santos2018comparing,bayona2019comparison,bayona2019role,miotti2021fully,shahane2021high,radhakrishnan2021non,unnikrishnan2022shear, bartwal2022application, shahane2022consistency, shahane2021semi, chu2022rbf}. In our recent works \cite{ shahane2021high, shahane2021semi, bartwal2022application, shahane2022consistency, unnikrishnan2022shear}, we had demonstrated high accuracy in heat conduction and fluid flow computations in several problems with manufactured solutions. These problems were limited to uniform properties in the solution domain and therefore, there were no interfaces that needed special attention in discretization of the derivatives. 
When the solution domain consists of regions with significant differences in properties, it is necessary to isolate the distinct regions and discretize the derivatives individually in each subregion. This requires point clouds that are limited to the homogeneous regions with common points on the interfaces. Some recent works with such heterogeneous media are by \citet{reutskiy2016meshless}, \citet{ahmad2020local} and \citet{gholampour2022efficient}. \citet{reutskiy2016meshless} considered heat conduction in irregular two-dimensional domains while \citet{ahmad2020local} demonstrated a subdomain RBF technique for geometries containing irregular interfaces with sharp corners. \citet{gholampour2022efficient} on the other hand, considered solution of elasticity problems with irregular interfaces by employing PHS-RBF and appended polynomial for 2D problems. It is worth mentioning here that the spectral element method \cite{karniadakis2005spectral} also decomposes a complex domain into elements and uses Chebyshev polynomial of high degree in each element. At element boundaries, the variable and its first derivative are made continuous. However, the motivation for such elementwise discretization and interface condition in spectral element methods is different from current motivation. When properties are modestly varying across the domain, the current cloud based method does not need special domain decomposition, unless needed for parallel computing purposes. 

\Cref{sec:methodology_single_domain} of this paper first briefly describes the PHS-RBF technique for single domains with homogeneous properties. In \cref{sec:methodology_multi_domain}, we present details of the extension of this algorithm to problems with regions of discontinuous and significant property variations. The accuracy of the multidomain algorithm is then evaluated by comparing it with the smearing approach as discussed in \cref{Sec:verif}. A detailed convergence analysis is presented in \cref{Sec:accuracy}. Subsequently, in \cref{Sec:results}, we apply the algorithm to more practical problems. \Cref{Sec:summary} provides summary of current observations and plans for future work.  

\section{Methodology}
\label{Sec:methodology}
\subsection{Single Domain PHS-RBF Algorithm}
\label{sec:methodology_single_domain}
\begin{figure}[H]
	\centering
	\includegraphics[width = 0.55\textwidth]{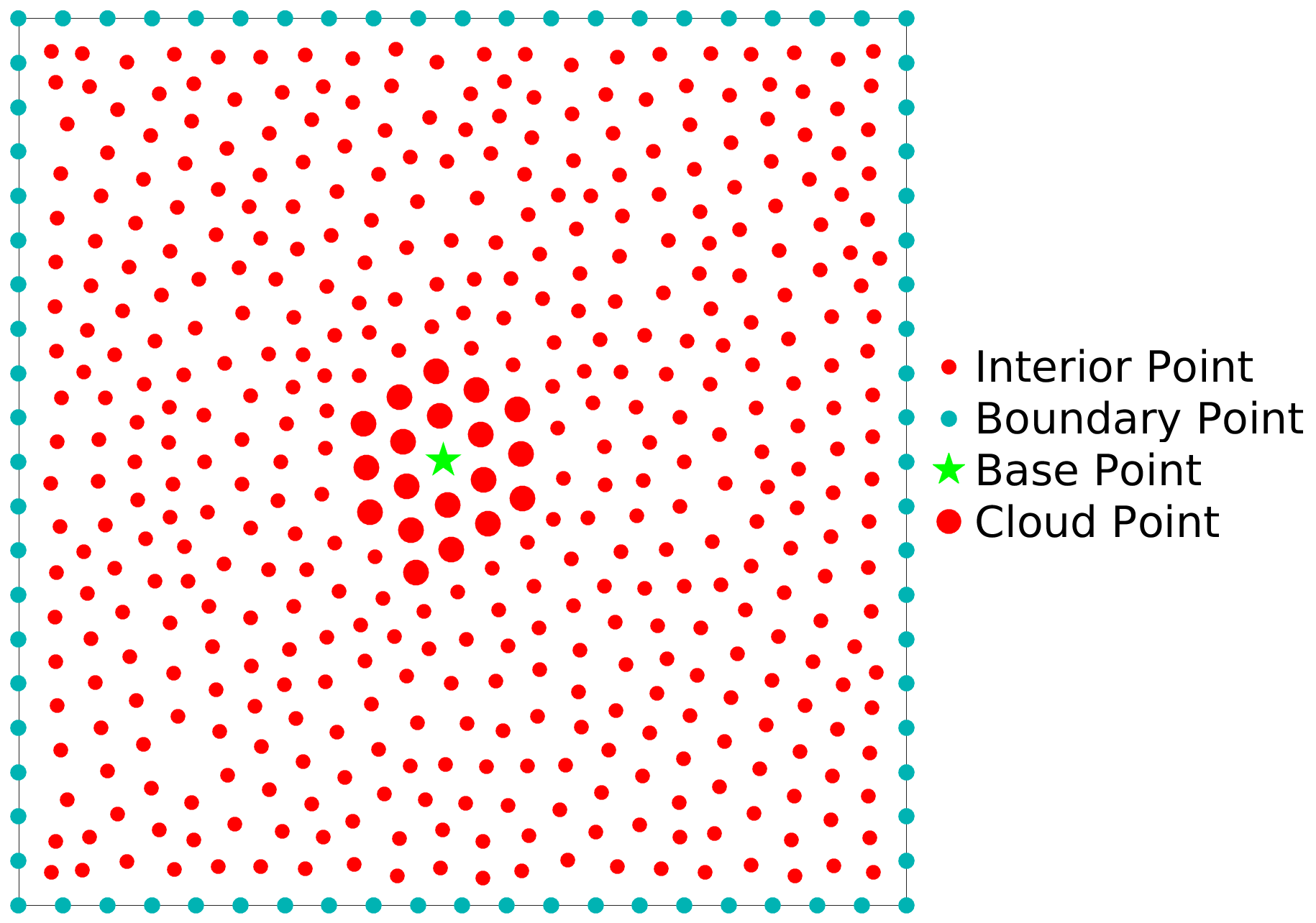}
	\caption{Methodology: Scattered distribution of points within a domain and example of a cloud.}
	\label{Fig:sample_domain}
\end{figure}

In this section, we first briefly describe the single domain cloud-based PHS-RBF algorithm to solve heat conduction and fluid flow problems in complex geometries. As this technique has been documented in our recent works \cite{ shahane2021high, shahane2021semi, bartwal2022application, shahane2022consistency, unnikrishnan2022shear}, as well as of some others \cite{barnett2015robust,flyer2016role,bayona2017role,santos2018comparing,bayona2019comparison,bayona2019role,miotti2021fully}, we describe the details only briefly for completeness and relevance to the multidomain algorithm. The starting point for the algorithm is the layout of scattered points in the complex solution domain. These scattered points can be generated by any of the techniques \cite{geuzaine2009gmsh,shankar2018robust,duh2021fast,van2021fast} including the vertices of any finite element grid generated by standard software. While the latter practice may appear as grid generation, it is to be mentioned that the elements and edges are eventually removed and thus, there are no issues of element skewness or aspect ratio which degrade the discretization accuracy. 
The second step is the definition of clouds for interpolation of neighbor values. For each scattered point, a subset of the nearest neighbor points is defined based on the degree of the appended polynomial. For a polynomial of degree $p$ and dimension of the problem $d$, the number of cloud points is given by $2\binom{p+d}{p}$. When the scattered point does not have enough interpolation points on one side, the cloud is extended in other directions asymmetrically. This happens near the external boundaries, or internal obstacles. For each scattered point, a separate cloud is defined. \Cref{Fig:sample_domain} shows a typical base point denoted by a green star and its cloud.                        
 The next step is computation of interpolation weights. For this interpolation, we use the PHS with appended polynomial as
\begin{equation}
s(\bm{x}) = \sum_{i=1}^{q} \lambda_{i} \phi(||\bm{x} - \bm{x_i}||_2) + \sum_{i=1}^{m} \gamma_i P_i (\bm{x})
\label{Eq:RBF_interp}
\end{equation}
where, $\phi \left(||\bm{x} - \bm{x_i}||_2\right)$ indicates polyharmonic spline radial basis function (PHS-RBF) and $||\bm{x} - \bm{x_i}||_2$ is the  distance between the location at which the the function $s$ is approximated (green coloured star shaped dot in \cref{Fig:sample_domain}) and different data points in its cloud (enlarged red coloured circular dots in \cref{Fig:sample_domain}). $\lambda_{i}$ and $\gamma_{i}$ are weights given to the PHS-RBF kernel and the appended monomials $P_{i}$, respectively. The PHS function $\phi$ is given as
$\phi = (||\bm{x}-\bm{x_{i}}||_2)^{2\alpha + 1}$ , where $\alpha$ is a positive integer.
If the values of $s$ are known at the scattered points, the above interpolation can be evaluated by collocating the data at the scattered points and thus providing a formula to evaluate $s(\bm{x})$ at any intermediate location. Let $\mathcal{L}$ be any linear differential operator of the governing partial differential equation. Applying $\mathcal{L}$ to \cref{Eq:RBF_interp} and collocating at the $q$ cloud points gives $q$ equations. The number of cloud points $q$ is taken to be twice the number of monomials, giving the total number of unknowns to be $(q + m)$. The additional $m$ equations for $\lambda$ are given by the constraints
\begin{equation}
\sum_{i=1}^{q} \lambda_i P_j(\bm{x_i}) =0 \hspace{0.5cm} \text{for } 1 \leq j \leq m
\label{Eq:RBF_constraint}
\end{equation}
Above collocation and constraint \cite{shahane2021high} written together in the matrix-vector form gives:
\begin{equation}
\begin{aligned}
\mathcal{L}[\bm{s}] &=
\left(\begin{bmatrix}
\mathcal{L}[\bm{\Phi}] & \mathcal{L}[\bm{P}]  \\
\end{bmatrix}
\begin{bmatrix}
\bm{A}
\end{bmatrix} ^{-1}\right)
\begin{bmatrix}
\bm{s}  \\
\bm{0} \\
\end{bmatrix}
\end{aligned}
\label{Eq:RBF_interp_mat_vec_L_solve}
\end{equation}
\begin{equation}
\text{where  }  
\begin{bmatrix}
\bm{A}
\end{bmatrix}
=
\begin{bmatrix}
\bm{\Phi} & \bm{P}  \\
\bm{P}^T & \bm{0} \\
\end{bmatrix}
\label{Eq:coeff_matrix}
\end{equation}
is a matrix of geometric values and monomials evaluated at the scattered points. The inverse of matrix $\bm{A}$ in \cref{Eq:coeff_matrix} is written only as a mathematical notation. In practice, a linear solver is used.
\Cref{Eq:RBF_interp_mat_vec_L_solve} can be further simplified as

\begin{equation}
\begin{aligned}
\mathcal{L}[\bm{s}] 
=
\begin{bmatrix}
\bm{B}
\end{bmatrix}
\begin{bmatrix}
\bm{s}  \\
\bm{0} \\
\end{bmatrix}
=
[\bm{B_1}] [\bm{s}] + [\bm{B_2}] [\bm{0}]
= [\bm{B_1}] [\bm{s}]
\end{aligned}
\label{Eq:RBF_interp_mat_vec_L_solve_simplified}
\end{equation}

$[\bm{B_1}]$ is of size $(q \times q)$ and contains coefficients for the linear differential operators such as gradient, Laplacian, etc. $[\bm{B_1}]$ depends only on the geometric distances between the individual points and its cloud points. \Cref{Eq:RBF_interp_mat_vec_L_solve_simplified} is satisfied at each one of the scattered points, generating a row of coefficients in the assembled matrix. This sparse coefficient matrix is of size $n \times n$ with $q$ non-zero coefficients in each row. The sparse equation set is then solved for the unknown discrete values $s_i$, i = 1 to $n$. In this work, we have used a preconditioned BiCGSTAB \cite{eigenweb}, although other solvers such as multilevel methods \cite{radhakrishnan2021non} can also be used. 
The above single domain algorithm was applied to several heat conduction problems in \cite{bartwal2022application} and was shown to achieve high accuracy. By considering different point spacings, we have shown that the discretization accuracy varies (in case of Laplacian) approximately as one order less than the degree of the appended polynomial and of the same order as the polynomial degree in case of gradients \cite{shahane2021high}.  However, if there are regions of significant and distinct variations in thermal properties, such accuracy can be lost. Consider the solution of the steady state heat conduction equation with a source term in a two dimensional space: 
\begin{equation}
\frac{\partial}{\partial x}\bigg(k\frac{\partial T}{\partial x}\bigg) + \frac{\partial}{\partial y}\bigg(k\frac{\partial T}{\partial y}\bigg) = \dot{q}
\label{Eq:smearing_govern}
\end{equation}
\Cref{Eq:smearing_govern} can be discretized either by first discretizing $k\frac{\partial T}{\partial x}$ and $k\frac{\partial T}{\partial y}$ and then taking their derivatives again or by expanding the derivatives as
\begin{equation}
k\bigg(\frac{\partial^{2}T}{\partial x^{2}} + \frac{\partial^{2}T}{\partial y^{2}}\bigg) + \frac{\partial k}{\partial x}\frac{\partial T}{\partial x} + \frac{\partial k}{\partial y}\frac{\partial T}{\partial y} = \dot{q} 
\end{equation}

\begin{figure}[H]
	\centering
	\includegraphics[width = 0.525\textwidth]{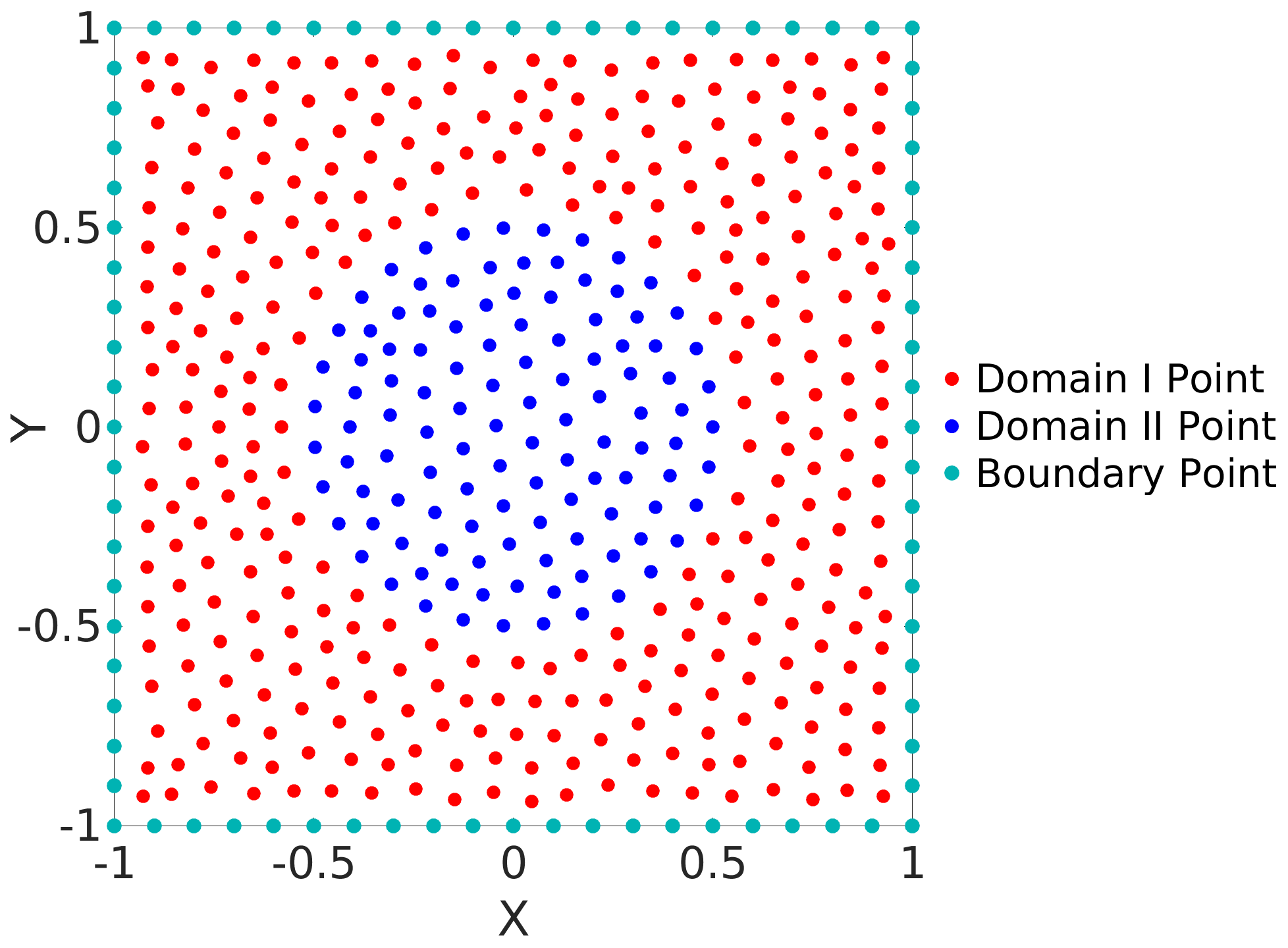}
	\caption{Methodology: Distribution of 531 scattered points within a domain with two regions of dissimilar conductivity.}
	\label{Fig:circle_rect_smearing_points}
\end{figure}

\begin{figure}[H]	
	\centering
	\subfigure[]{\includegraphics[width=0.45\textwidth]{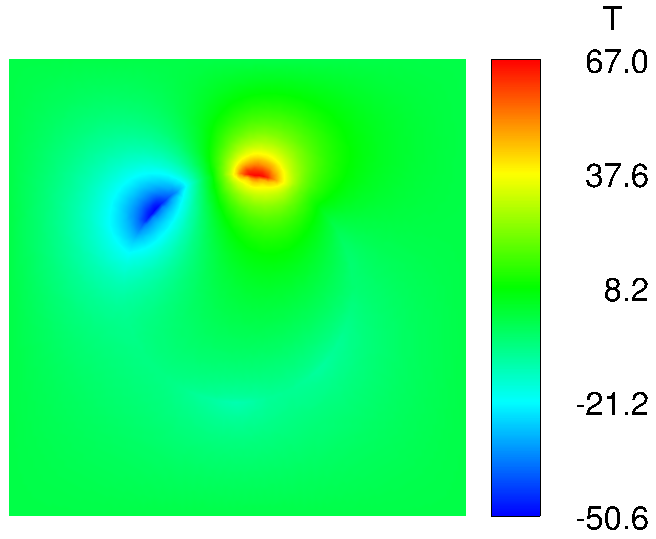}\label{fig:smearing_3}}
	\subfigure[]{\includegraphics[width=0.466\textwidth]{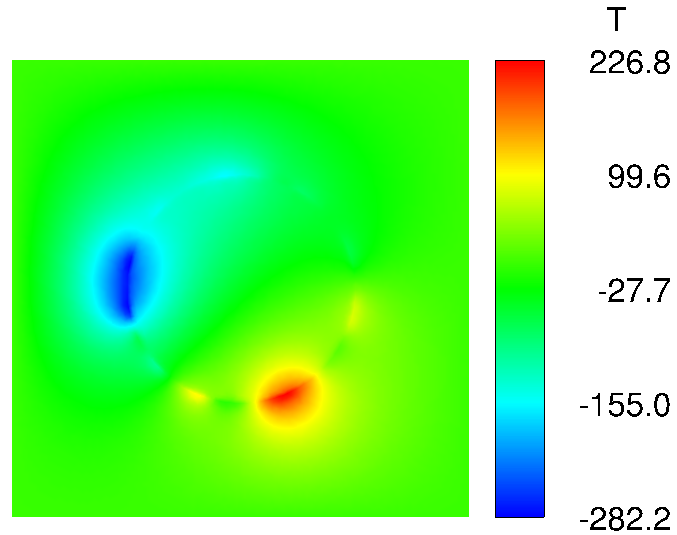}\label{fig:smearing_5}}
	\caption{Methodology: Temperature contours (a) polynomial degree = 3 and (b) polynomial degree = 5 for 12078 scattered points.}
	\label{Fig:circle_rect_smearing_contours}
\end{figure}

\begin{figure}[H]
	\centering
	\includegraphics[width = 0.525\textwidth]{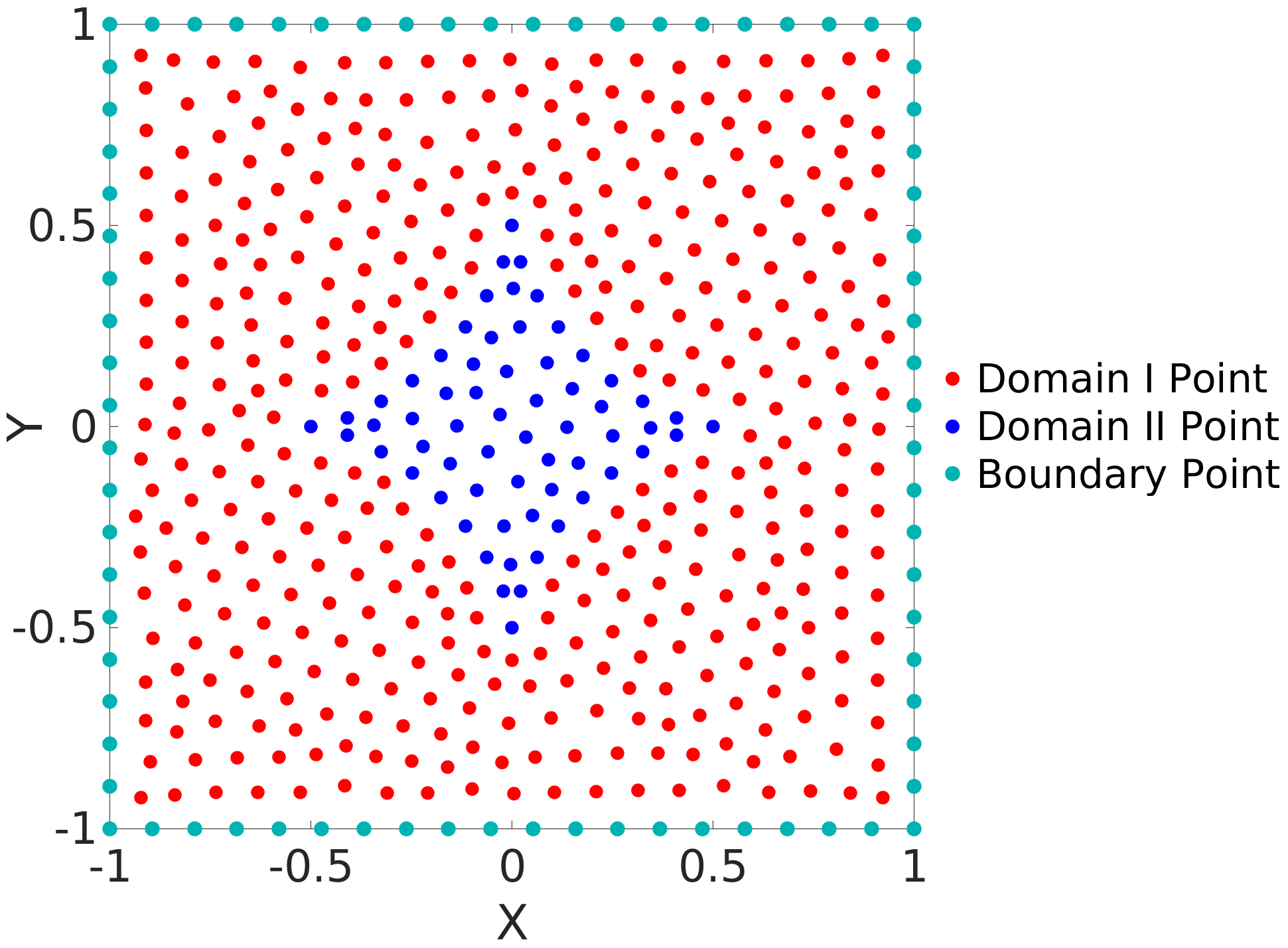}
	\caption{Methodology: Distribution of 540 scattered points within a domain with two regions of dissimilar conductivity.}
	\label{Fig:astroid_rect_smearing_points}
\end{figure}

\begin{figure}[H]	
	\centering
	\subfigure[]{\includegraphics[width=0.47\textwidth]{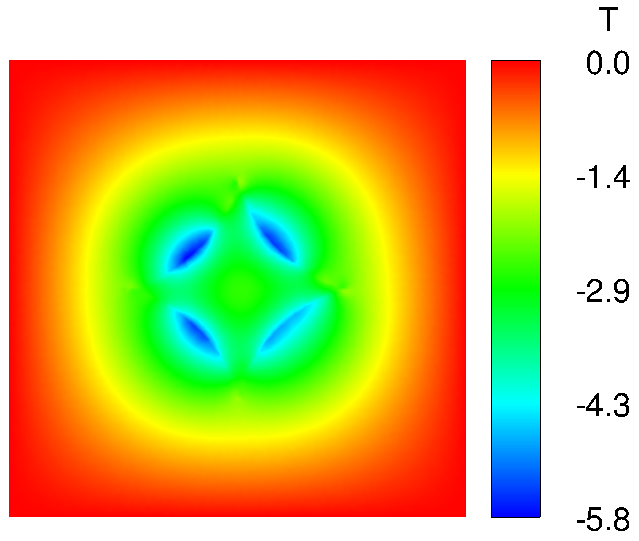}\label{fig:smearing_3}}
	\subfigure[]{\includegraphics[width=0.488\textwidth]{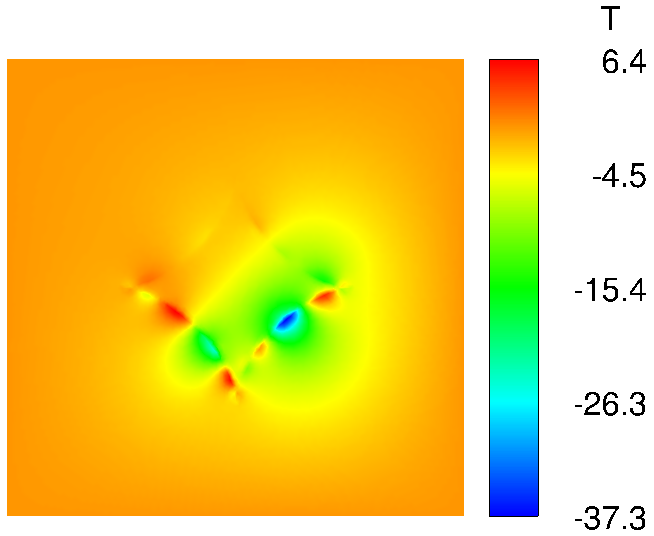}\label{fig:smearing_5}}
	\caption{Methodology: Temperature contours (a) polynomial degree = 3 and (b) polynomial degree = 5 for 12060 scattered points.}
	\label{Fig:astroid_rect_smearing_contours}
\end{figure}

In the above expansion, we can then evaluate the first and second derivatives using the kernel differentiation procedure explained above.  The $\bm{B_1}$ matrix mentioned in \cref{Eq:RBF_interp_mat_vec_L_solve_simplified} now contains the coefficients from the Laplacian and from the first derivatives of the temperature ($T$). The sharp discontinuity in the conductivity is manifested by the derivatives $\frac{\partial k}{\partial x}$ and  $\frac{\partial k}{\partial y}$, which needs to be evaluated using the RBF first derivative matrix. Although, the analytical derivatives of this step function do not exist, the RBF computes it using the stencil for the first derivative, much like computing a derivative across a shock in shock capturing techniques. 
We have evaluated this procedure by solving for a manufactured solution with two regions of distinct conductivity. Three ratios of the thermal conductivities $\frac{k_I}{k_{II}}$ = $5$, $10$ and $100$ with different point spacings and degree of the appended polynomial were considered. We show the results for two different problems: (a) circle inside a square (\cref{Fig:circle_rect_smearing_points}) and (b) astroidal domain inside a square (\cref{Fig:astroid_rect_smearing_points}).  \Cref{Fig:circle_rect_smearing_contours,Fig:astroid_rect_smearing_contours} show the temperature contours for the interior and boundary conditions given by \cref{Eq:smearing_bc}.
\begin{equation}
\begin{aligned}
\dot{q} &= 100 \text{ at every node inside the circular and astroidal domain (blue dots)}\\
T &= 0 \text{ at the rectangular boundary (cyan dots)}\\
k &= 10 \text{ (red and cyan dots) } \text{and } k = 1 \text{ (blue dots)}
\label{Eq:smearing_bc}
\end{aligned}
\end{equation}
For conductivity ratio of 10, we see that the single domain technique which interpolates the temperatures and conductivities through the interface, significantly smears the discontinuity and results in erroneous results as shown by the contours of non-dimensional temperature (\cref{Fig:circle_rect_smearing_contours,Fig:astroid_rect_smearing_contours}). The larger the polynomial, the bigger is the cloud, and generates larger errors contrary to the expected higher accuracy.

\newpage
\subsection{Algorithm for Multiple Domains}
\label{sec:methodology_multi_domain}
As shown above, interpolation across the interfaces results in large errors and destroys the high spatial accuracy observed for domains with homogeneous thermal properties. This is because the clouds are defined isotropically that leads to interpolations through the interfaces.  It is obvious that the heterogeneous regions should be isolated and separately interpolated and treated as separate subdomains. This is similar to the well-known shock patching techniques used in computation of high-speed gas dynamics problems. In this section, we describe the multidomain technique with the cloud based PHS-RBF interpolation.  The verification and accuracy of the technique for several problems with manufactured solutions are then demonstrated followed by applications of the algorithm to two and three dimensional heat conduction in complex geometries.    

\begin{figure}[H]
	\centering
	\subfigure{\includegraphics[width=0.55\textwidth]{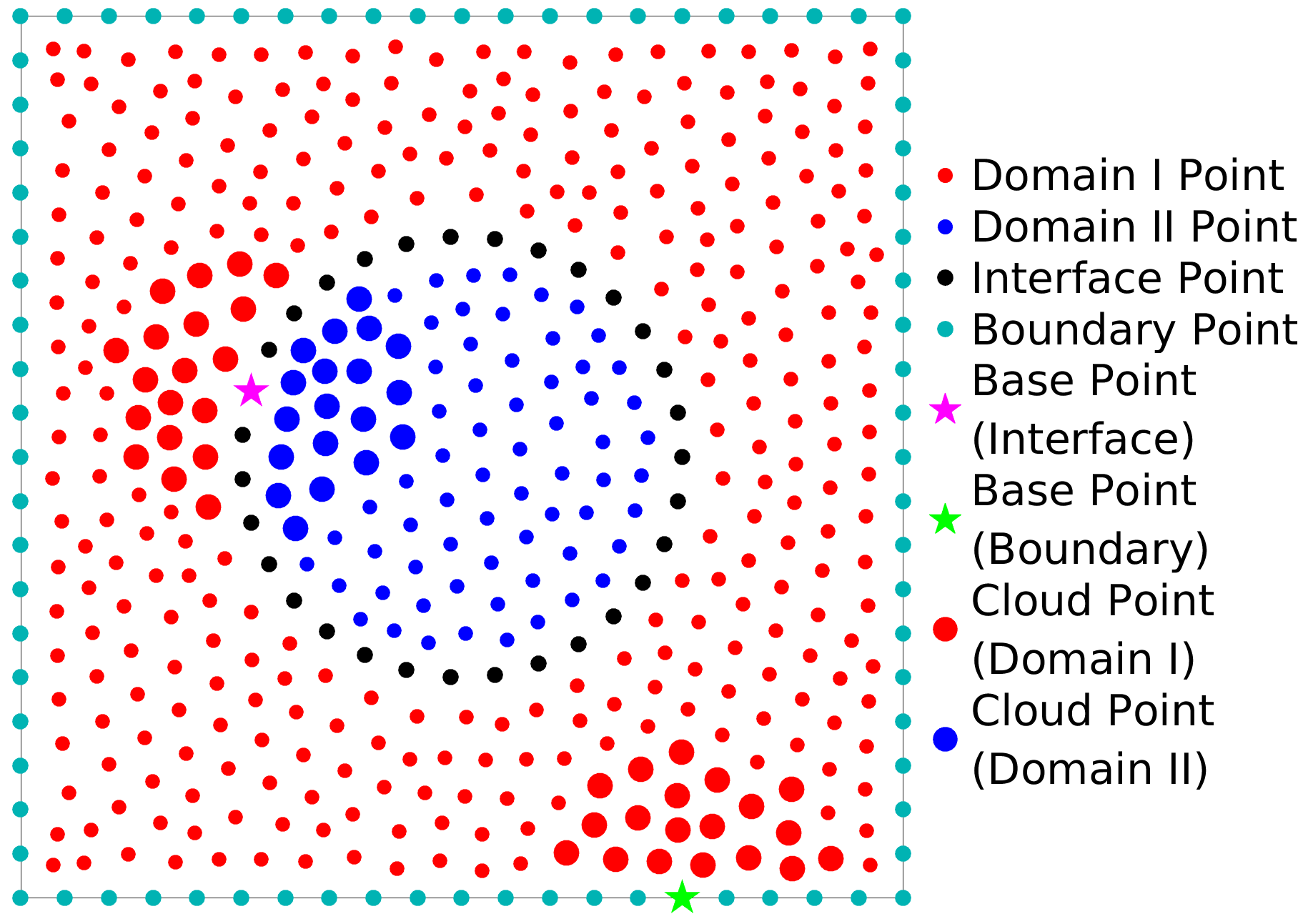}\label{fig:clamp_left}}
	\caption{Methodology: Clouds formed by points located at the interface (magenta coloured star shaped dot) and boundary (green coloured star shaped dot) at polynomial degree = 3} 
	\label{Fig:method_interface_boundary_nodes}
\end{figure}

\begin{figure}[H]
	\centering
	\subfigure{\includegraphics[width=0.55\textwidth]{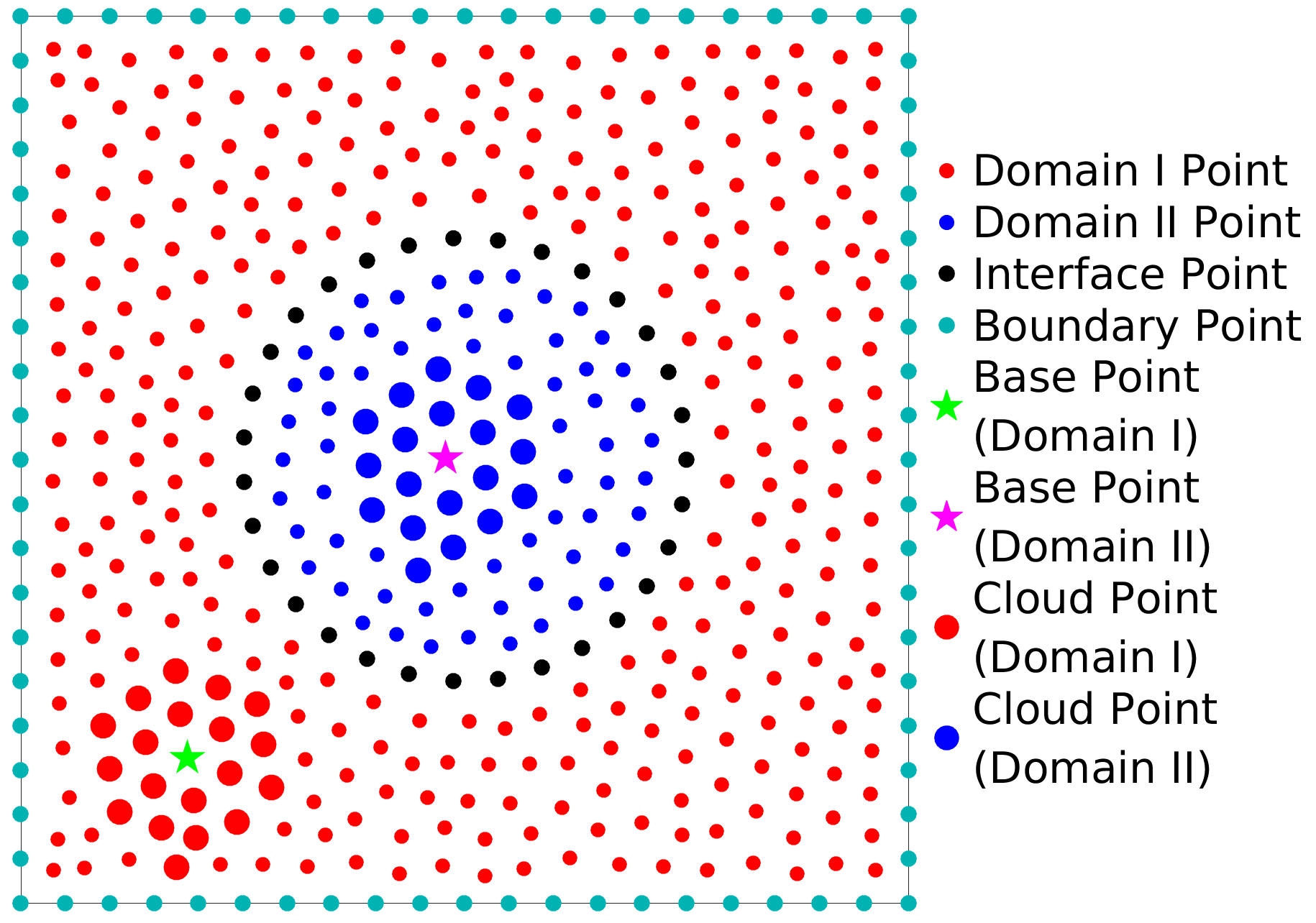}\label{fig:clamp_left}}
	\caption{Methodology: Clouds formed by points located in the inner (magenta coloured star shaped dot) and outer domain (green coloured star shaped dot) at polynomial degree = 3} 
	\label{Fig:method_inner_outer_domain}
\end{figure}

The algorithm again begins with the layout of scattered points in the solution domain. However, we first define the interfaces separating the various regions, and layout scattered points on the interfaces. For 2D problems, we layout points along the line  interfaces. For 3D problems, we first triangulate the surfaces of the interfaces, and remove the edges and elements. Then, using these interface points and the points prescribed on the external boundaries, we generate separate sets of points in all the subdomains. As before, the points can be vertices of elements of a finite element grid. 
The governing heat conduction equation is then individually satisfied in each subdomain, without any interpolation across the domains. On the interface, the partial differential equation is not satisfied. Instead, we impose the $C_0$ and $C_1$ conditions at each of the interface   points. Thus, for all the interface points the flux balance condition is satisfied,
\begin{equation}
k_{1} \frac{\partial T}{\partial \hat{n}}\bigg\vert_1 = k_{2} \frac{\partial T}{\partial \hat{n}}\bigg\vert_2
\label{Eq:method_interface_condition}
\end{equation}
The $C_0$ condition is automatically satisfied by the common temperatures at the interface. 
The next step is the definition of the cloud points. This is where the main difference arises. In the multidomain algorithm, the clouds are defined to contain only points in their own subdomain. Hence near the interfaces, the clouds are defined asymmetrically, thus avoiding interpolation across the interfaces. \Cref{Fig:method_interface_boundary_nodes,Fig:method_inner_outer_domain} show three types of clouds that can result in the simple case of an interface, boundary and interior points for degree of appended polynomial equal to 3.

 Clouds in the interiors of each subdomain are isotropic if there are enough neighbor points available. If sufficient cloud points are not available, points are picked in the other directions and are restricted to be in their own subdomain as in the case of the external boundaries. A point on the interface has two subdomains on the two sides of its normal. Hence, two sets of cloud points, one in each subdomain are defined. These two clouds are used to compute the normal derivatives and the fluxes on each side using the conductivities in the appropriate subdomain.  The flux balance condition at each interface point then gives an additional equation to close the system.  
The final step in the algorithm is the assembly of the coefficient matrix. We assemble one single coefficient matrix for all subdomains combined. This matrix consists of the satisfaction of the partial differential equation at all interior points in each subdomain and the flux balance condition at all the points on the interfaces.  At the external boundaries, based on the type of boundary condition (Dirichlet or Neumann), we fix the value, or include an implicit equation for the derivative. In such a case, clouds are defined for the points on the external boundaries as well. The total number of discrete equations is therefore the sum of all points interior to the subdomains, points on the interfaces and points on external boundaries. The number of discrete equations equals the number of discrete unknowns, making the problem well posed.  As before, we solve the set of equations using a preconditioned BiCGSTAB algorithm, after reorganizing the matrix using the well known reverse Cuthill-Mckee (RCM) algorithm \cite{cuthill1969reducing} which reduces the bandwidth of the sparse matrix.

\begin{figure}[H]
	\centering
	\subfigure[]{\includegraphics[width=0.48\textwidth]{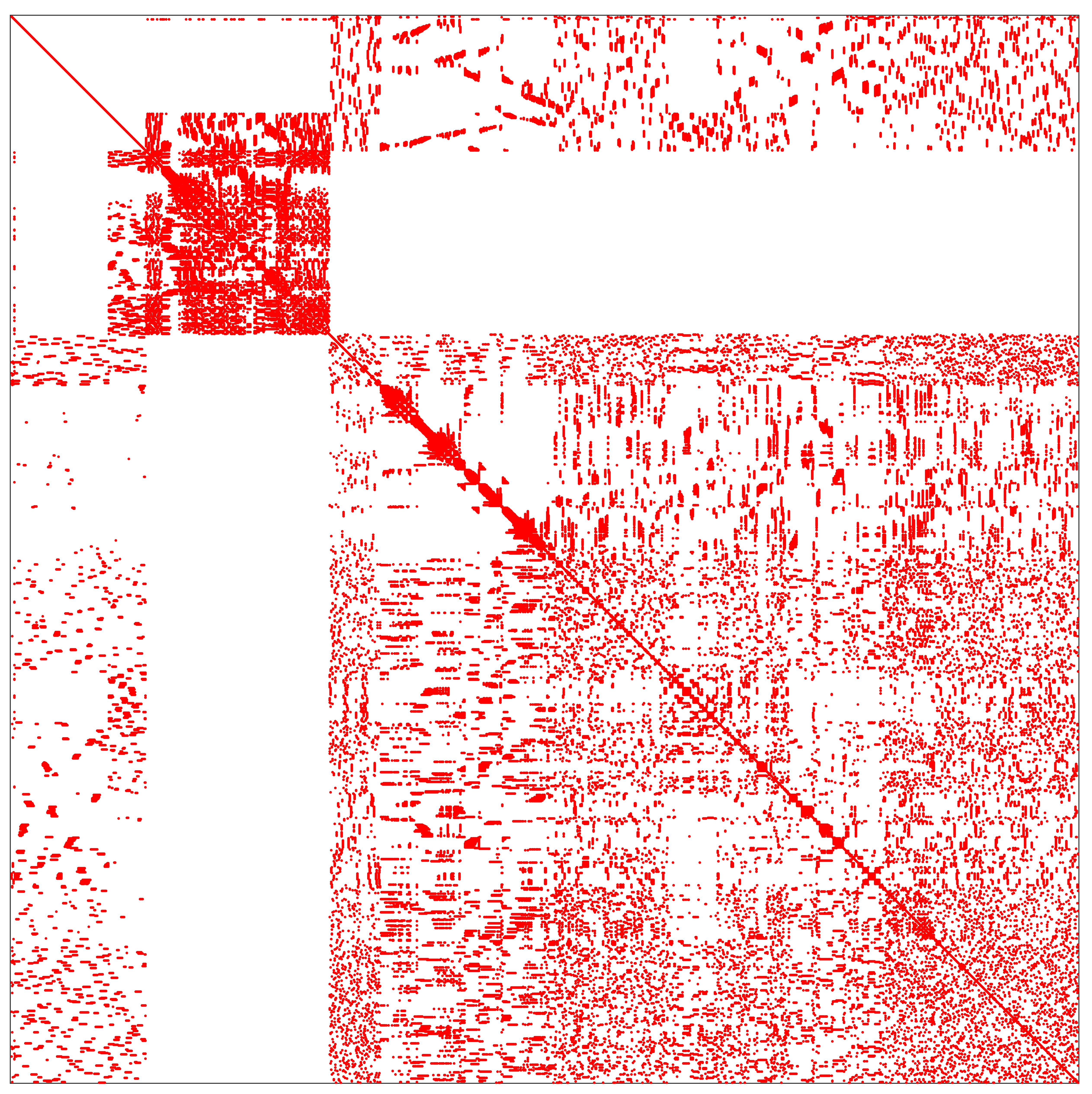}}
	\subfigure[]{\includegraphics[width=0.48\textwidth]{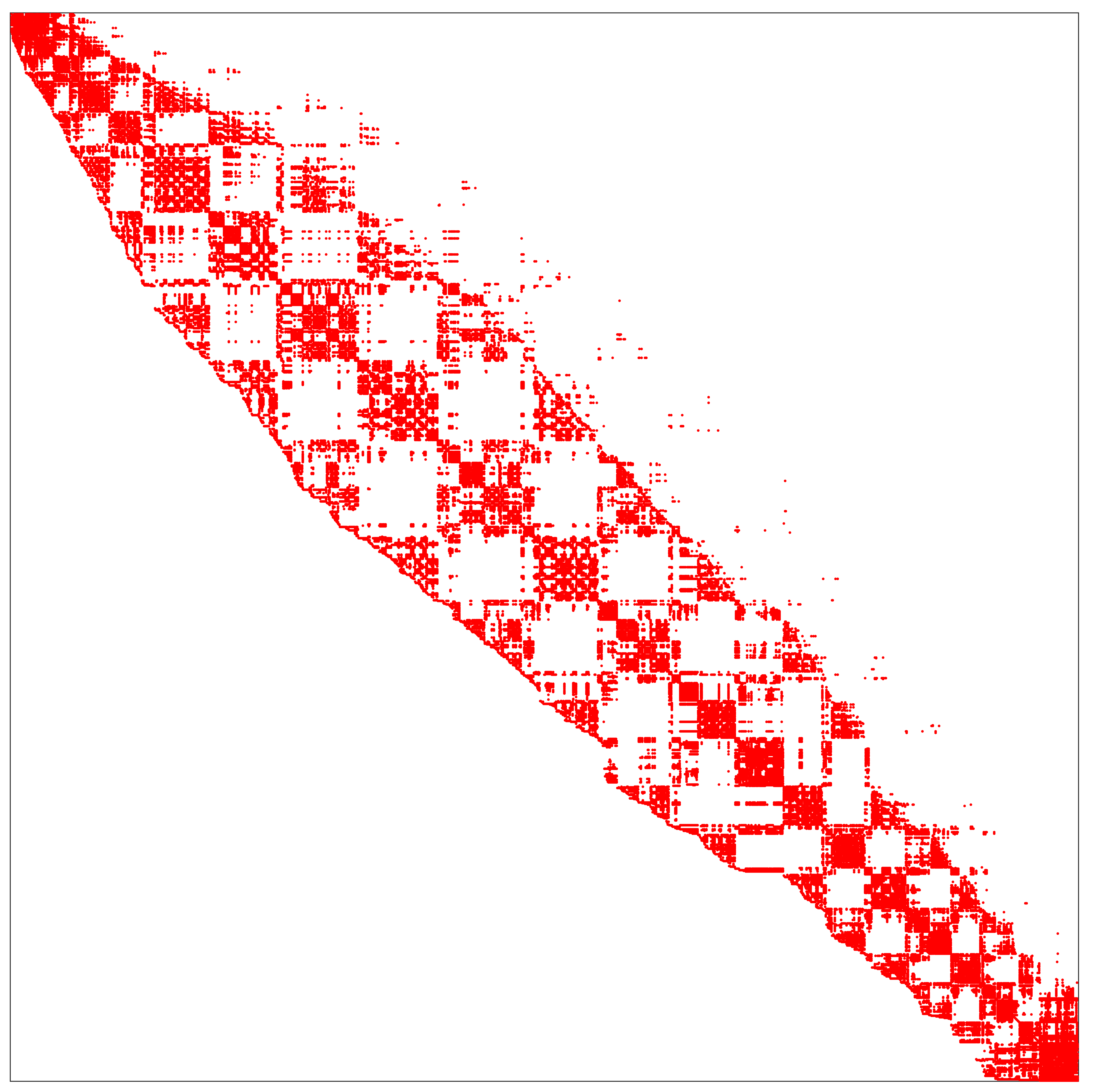}}
	\caption{Methodology: Sparsity pattern of the final coefficient matrix for the multidomain algorithm with 1541 points: (a) before RCM and (b) after RCM}
	\label{Fig:method_rcm_multi_domain}
\end{figure}

 \Cref{Fig:method_rcm_multi_domain} shows the sparsity of the final assembled matrix before and after the application of RCM for 1541 scattered distribution of points within the domain. In summary, the multidomain algorithm differs from our previously published single domain algorithm in the following steps: 

\begin{enumerate}[(a)]
	\item
	The scattered points are generated separately for each subdomain after deciding the points on the interior interfaces.
	\item
	The clouds for each scattered point in the interior of the subdomains are defined so as to include only points within the subdomain (together with the interface points) and not to include points across the interfaces.
	\item
	The governing partial differential equation is satisfied only at the points interior to the subdomains and not on the interfaces.
	\item
	At points on the interfaces, flux balance equations are satisfied using two sets of clouds spanning into the two adjacent domains along the local normal.
\end{enumerate}
 Note that normals at each scattered point on the interface are required for the flux balances. These normals are precomputed at the time of point generation from information of the finite element grid before deleting the elements and edges. \Cref{Fig:flow_chart_methodology} shows the steps adopted for multidomain algorithm in the form of a flowchart. In the following sections we validate the algorithm and evaluate the order of accuracy for different conductivity ratios in two and three dimensional problems.  We then apply the algorithm to compute heat conduction in complex composite material configurations.

\begin{center}
\begin{figure}[H]
\begin{tikzpicture}
\tikzstyle{terminator} = [rectangle, draw, text centered, rounded corners, minimum height=2em]
\node at (0,0) [terminator] (t-id0) {Generation of points within the domain};

\tikzstyle{terminator1} = [rectangle, draw, text centered, rounded corners, minimum height=2em]
\node at (0,-1.7) [terminator1] (t-id1) {Identifying points constituting different interfaces and subdomains};

\tikzstyle{terminator2} = [rectangle, draw, text centered, rounded corners, minimum height=2em, text width=1.0\textwidth]
\node at (0,-4.4) [terminator2] (t-id2) [draw, align=justify]{Defining cloud points for every point within the solution domain. The clouds are defined to contain only points in their own subdomain. Two sets of cloud points are identified for interface points as an interface point has two subdomains on the two sides of its normal};

\tikzstyle{terminator3} = [rectangle, draw, text centered, rounded corners, minimum height=2em, text width=1.0\textwidth]
\node at (0,-7.4) [terminator3] (t-id3) [draw, align=justify]{Calculation of weighted coefficients at every point within the domain for the given linear differential operation};

\tikzstyle{terminator4} = [rectangle, draw, text centered, rounded corners, minimum height=2em]
\node at (0,-9.4) [terminator4] (t-id4) [draw, align=center]{Flux balance condition is satisfied at the interface points};

\tikzstyle{terminator5} = [rectangle, draw, text centered, rounded corners, minimum height=2em,text width=1.0\textwidth]
\node at (0,-11.35) [terminator5] (t-id5) [draw, align=justify]{ Assemble a single sparse matrix for all the subdomains and interfaces combined};

\tikzstyle{terminator6} = [rectangle, draw, text centered, rounded corners, minimum height=2em]
\node at (0,-13.3) [terminator6] (t-id6) [draw, align=center]{Reorganizing the matrix using reverse Cuthill-Mckee (RCM) algorithm};

\tikzstyle{terminator7} = [rectangle, draw, text centered, rounded corners, minimum height=2em,text width=1.0\textwidth]
\node at (0,-15.3) [terminator7] (t-id7) [draw, align=justify]{Solving the final assembled matrix-vector system using preconditioned BiCGSTAB algorithm};

\tikzstyle{connector} = [draw, -latex']
\path [connector] (t-id0) -- (t-id1);
\path [connector] (t-id1) -- (t-id2);
\path [connector] (t-id2) -- (t-id3);
\path [connector] (t-id3) -- (t-id4);
\path [connector] (t-id4) -- (t-id5);
\path [connector] (t-id5) -- (t-id6);
\path [connector] (t-id6) -- (t-id7);
\end{tikzpicture}
\caption{Flowchart showing steps adopted for multidomain algorithm.}
\label{Fig:flow_chart_methodology}
\end{figure}
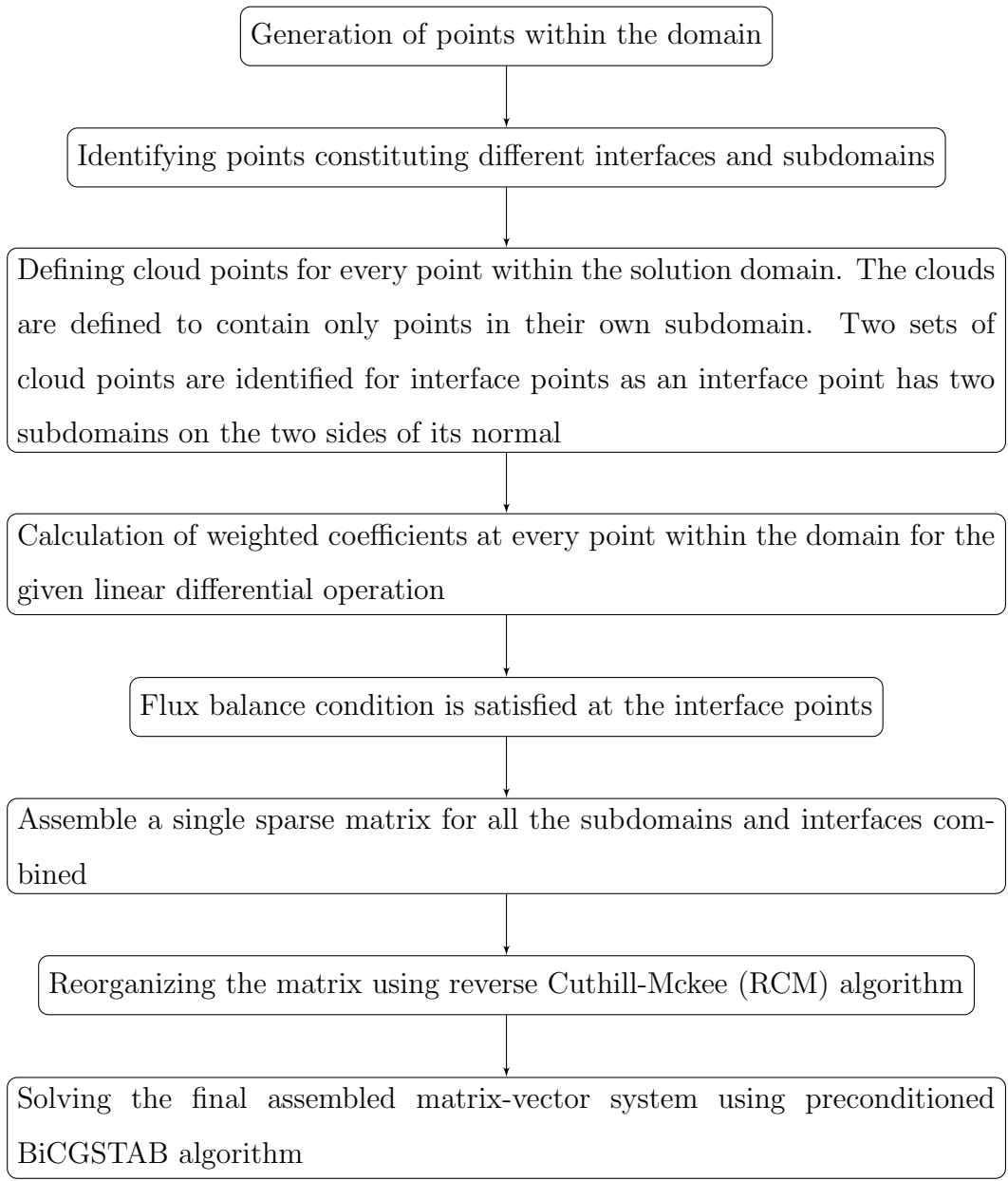
\end{center}

\section{Verification}
\label{Sec:verif}
\begin{figure}[H]
	\centering
	\subfigure[]{\includegraphics[width=0.57\textwidth]{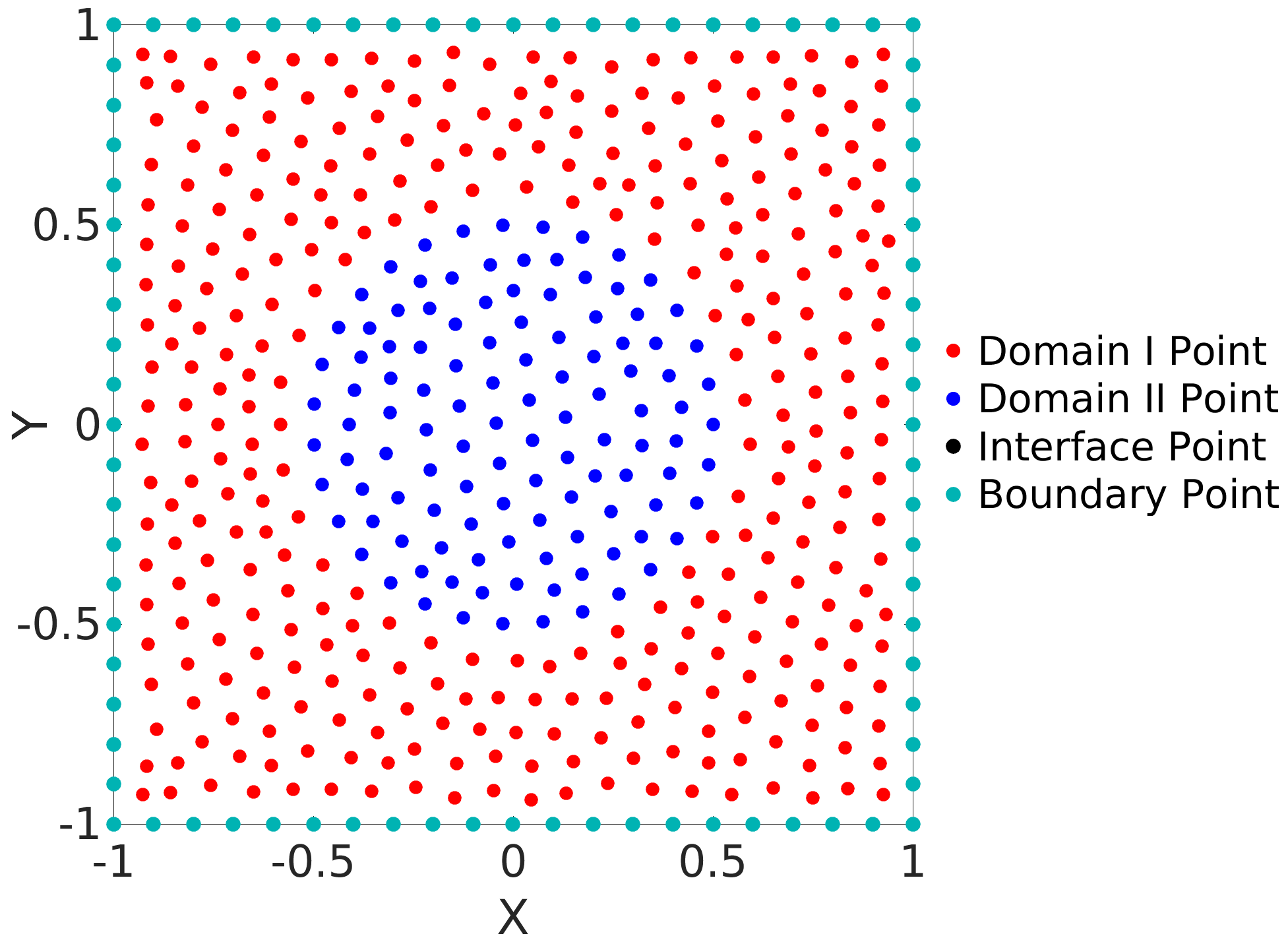}\label{fig:assess_accuracy_smearing}}
	\subfigure[]{\includegraphics[width=0.42\textwidth]{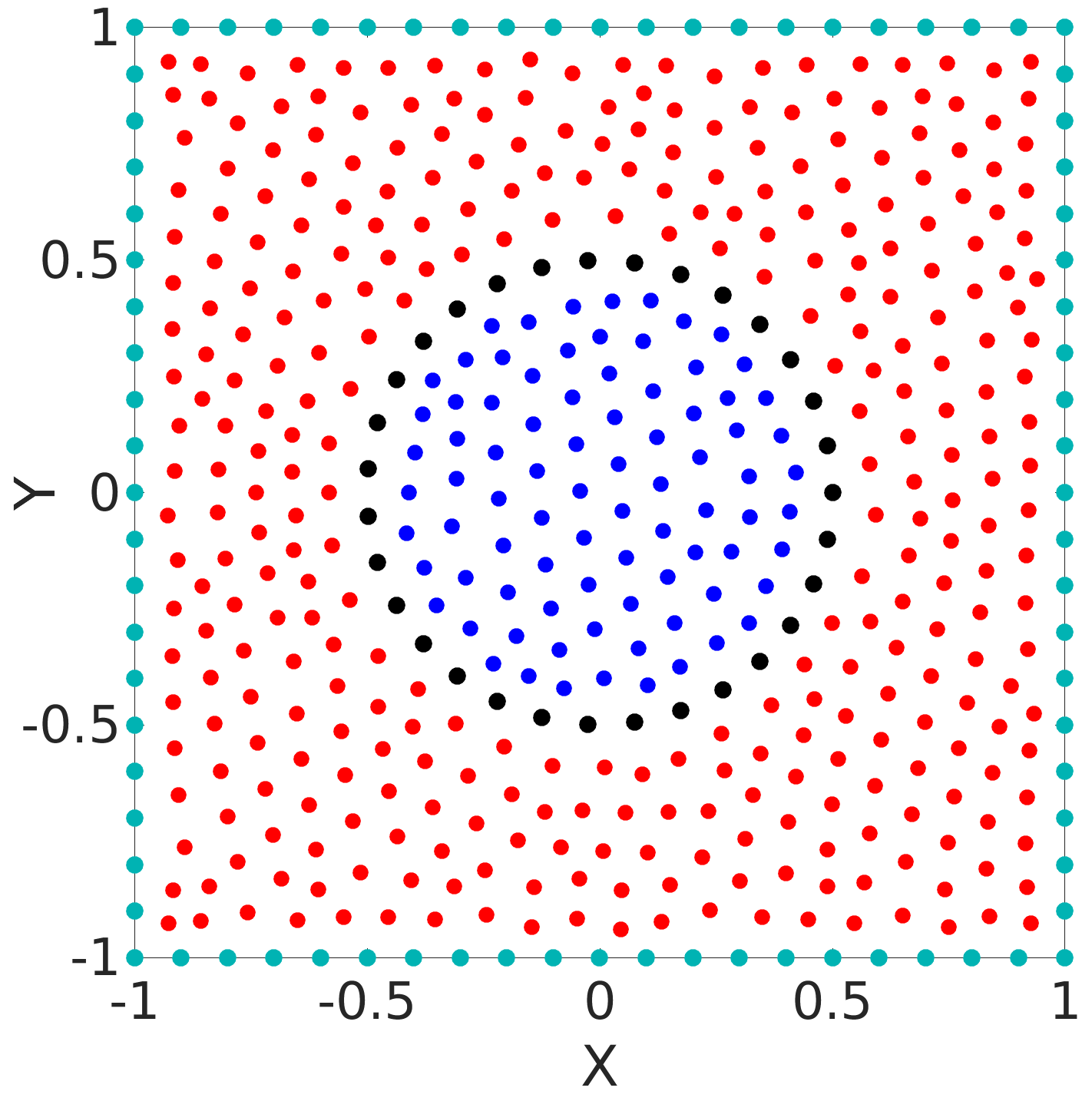}\label{fig:assess_accuracy_patching}}
	\caption{Verification: Point distributions (a) without and (b) with an explicit interface.}
	\label{Fig:smearing_vs_patching_nodes}
\end{figure}

In this section, we verify the implementation of the multidomain algorithm by considering a problem of equal thermal conductivities in two subdomains. An artificial circular domain is separated by distinct scattered points and the results of the single and multidomain algorithms are compared. In the multidomain algorithm, an interface is also defined at which the differential equation is not satisfied but the flux is balanced. \Cref{Fig:smearing_vs_patching_nodes} shows the distribution of 531 points within the domain (Red dots $\rightarrow$ Domain I, Blue dots $\rightarrow$ Domain II and Black dots $\rightarrow$ Interface). Open source software Gmsh \cite{geuzaine2009gmsh} is used to generate the points by the technique of Delaunay triangulation.

\begin{figure}[H]
	\centering
	\includegraphics[width = 0.6\textwidth]{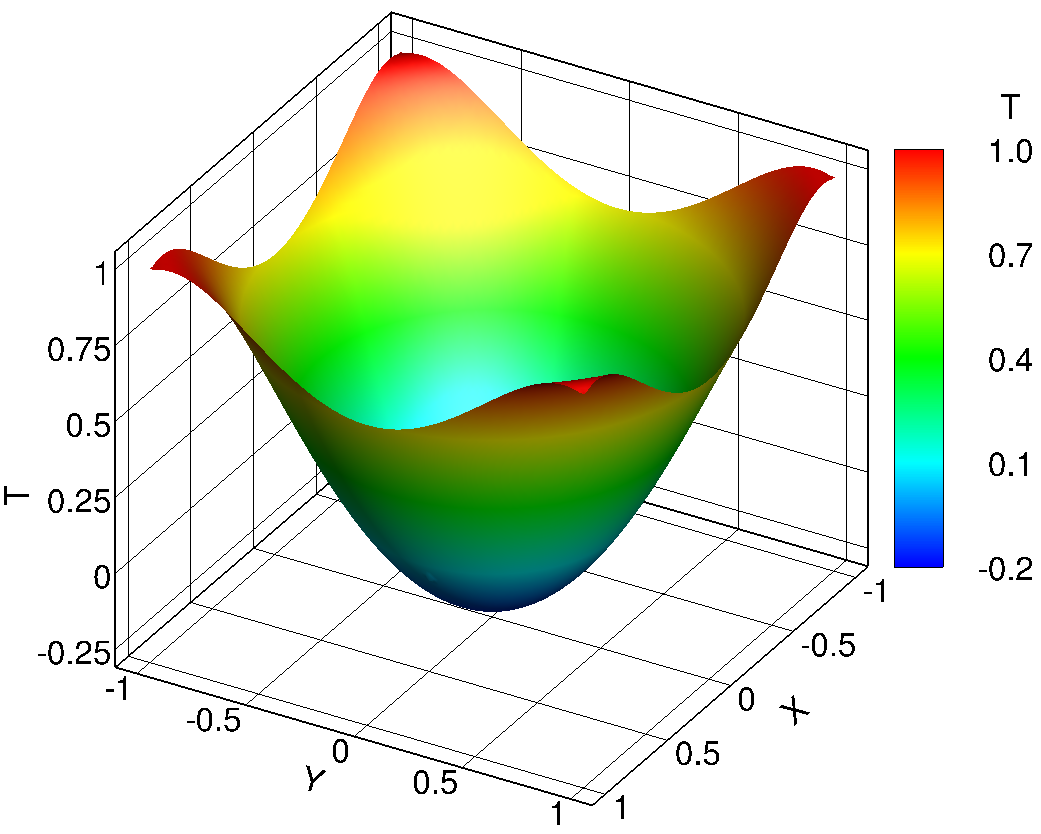}
	\caption{Verification: Surface plot of the manufactured temperature distribution.}
	\label{Fig:circle_rect_manuf_solution_3D_assess_accuracy}
\end{figure}

We consider equal thermal conductivities in domain I and domain II. We take two different scenarios. In the first one as shown in \cref{fig:assess_accuracy_smearing}, there is no interface and the cloud of a particular point is not limited to its respective subdomain but may contain points from other domain. In the second case, as depicted by \cref{fig:assess_accuracy_patching}, the cloud is restricted to individual subdomains, and the flux conditions are satisfied at the interface points. As the thermal conductivities are equal in the subdomains, we must get comparable results using both the approaches. We solve the heat conduction equation for a manufactured solution for temperature given by:
\begin{equation}
\begin{aligned}
T_{I} &= k_{I} \; sin(x^{2} + y^{2} - r^{2})\\ 
T_{II} &= k_{II} \; sin(x^{2} + y^{2} - r^{2})
\label{Eq:assess_accuracy_manuf_solution}
\end{aligned}
\end{equation}
where, radius of circular domain  $r$ = $0.5$. A surface plot of the manufactured temperature distribution is shown in \cref{Fig:circle_rect_manuf_solution_3D_assess_accuracy}.

\begin{figure}[H]
	\centering
	\subfigure[]{\includegraphics[width=0.555\textwidth]{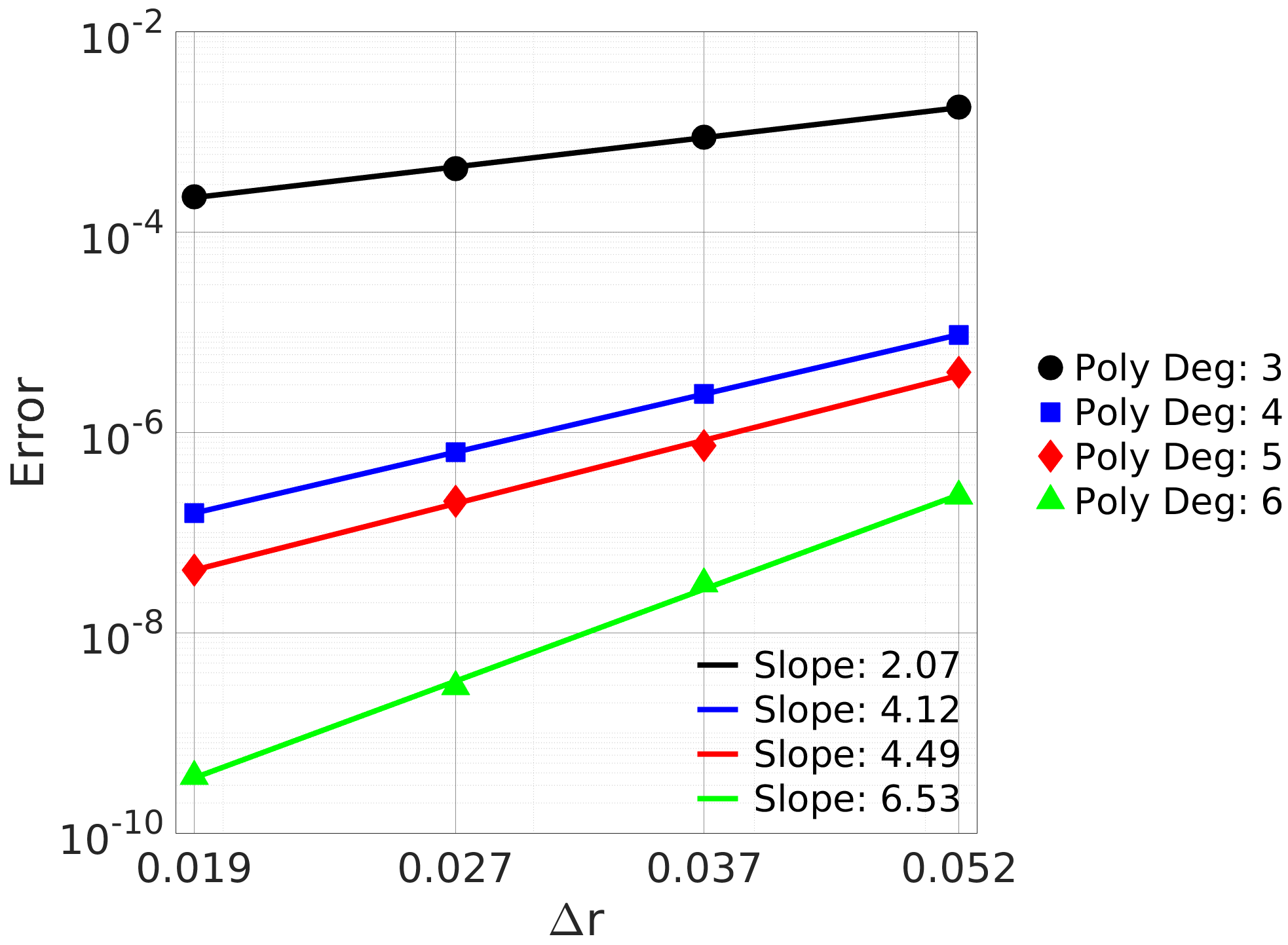}\label{fig:assess_accuracy_smearing_error}}
	\subfigure[]{\includegraphics[width=0.435\textwidth]{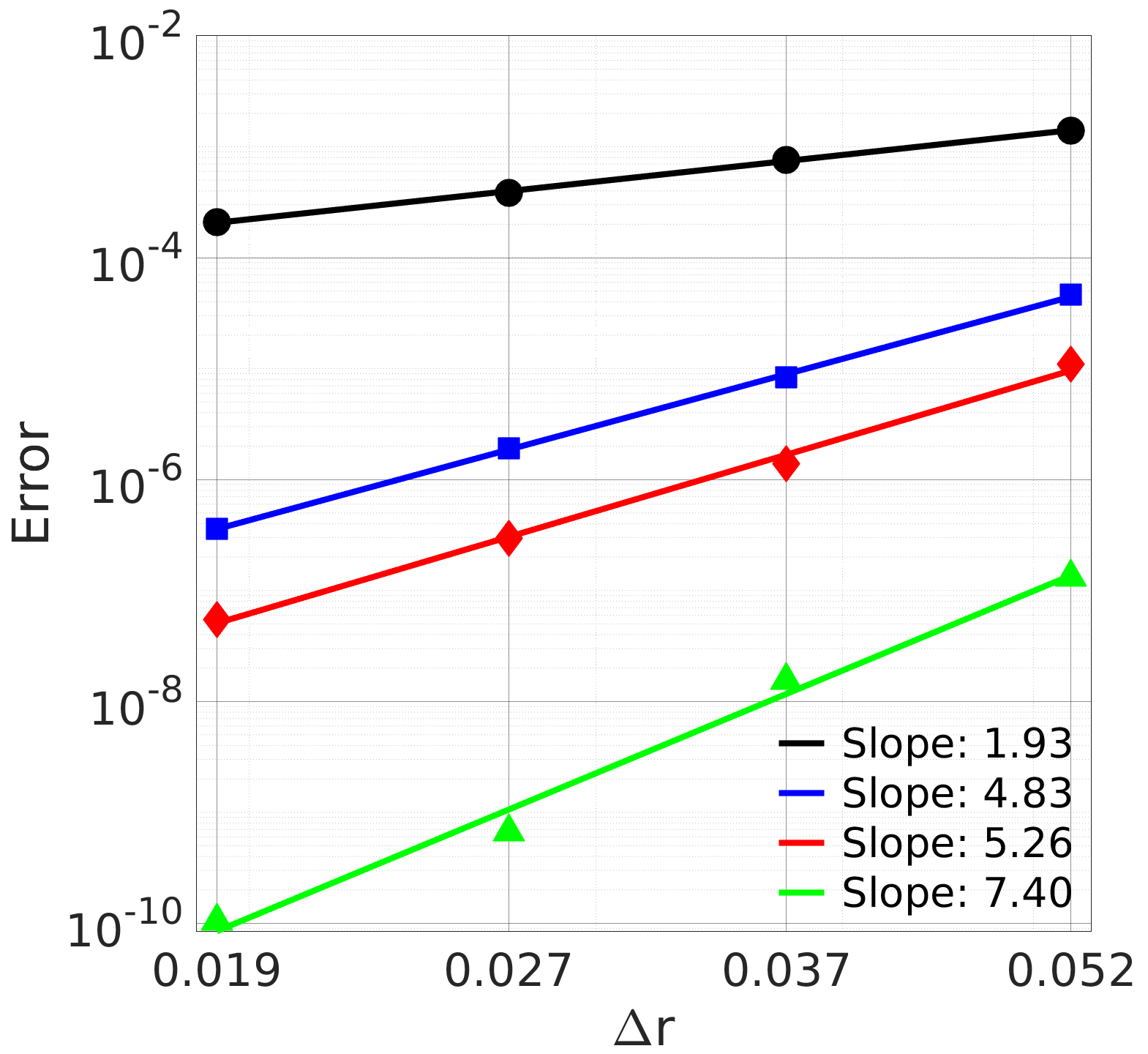}\label{fig:assess_accuracy_patching_error}}
	\caption{Verification: Variation of error with average point spacing for (a) single domain and (b) multi domain calculations with equal conductivities.}
	\label{Fig:smearing_vs_patching_error}
\end{figure}

Four different point distributions with 1541, 3089, 6033 and 12078 points having an average inter-point spacing of 0.052, 0.037, 0.027 and 0.019 are considered and the degree of appended polynomial is varied from 3 to 6. \Cref{Fig:smearing_vs_patching_error} shows the log-log variation of L1-norm of the discretization error with respect to average grid spacing for both the strategies of defining the clouds. The closest node is identified for every point and the average of these distances is defined as $\Delta \text{r}$. 

\par The error for the entire domain given as

\begin{equation}
Error_{Domain} = \frac{ \sum_{n = 1}^{N} |T_{exact} - T_{calc}| }{|T_{max \_ exact}| \times N}
\label{Eq:full_domain_error}
\end{equation}
where, $N$ is the number of points in the combined domain, $T_{exact}$, $T_{calc}$ and $T_{max\_exact}$ are the manufactured solution, computed solution and the maximum exact value, respectively. The order of convergence is approximately determined by a line of best fit. We show that the two approaches give the expected order of convergence for a given degree of appended polynomial. This verifies the domain decomposition algorithm to produce accurate results and verifies the implementation of the flux balance condition. 

\section{Study of Multidomain Accuracy with Manufactured Solutions}
\label{Sec:accuracy}
In this section, we apply the multidomain algorithm to study the robustness and accuracy in two and three dimensional geometries, all with manufactured solutions. The geometries are: (i) circle inside a square (ii) astroidal domain inside a square  and (iii) sphere inside a cube. For each case, the solution domain is divided into domain I, domain II and an interface. Each problem is investigated for three thermal conductivity ratios of 5, 10 and 100. Non-dimensional errors in the computed temperatures are estimated using the manufactured solutions as follows:

\begin{equation}
Error_{I} = \frac{ \sum_{n = 1}^{N_{I}} |T_{exact} - T_{calc}| }{(|T_{exact\_max\_I} - T_{exact\_min\_I}|) \times N_{I}}
\label{Eq:I_error}
\end{equation}
\begin{equation}
Error_{II} = \frac{ \sum_{n = 1}^{N_{II}} |T_{exact} - T_{calc}| }{(|T_{exact\_max\_II} - T_{exact\_min\_II}|) \times N_{II}}
\label{Eq:II_error}  
\end{equation}
\begin{equation}
Error_{Interface} = \frac{ \sum_{n = 1}^{N_{Int}} |T_{exact} - T_{calc}| }{\frac{|T_{exact\_max\_I} - T_{exact\_min\_I}| +  |T_{exact\_max\_II} - T_{exact\_min\_II}|}{2} \times N_{Int}}
\label{Eq:interface_error}
\end{equation}
Here, $N_{I}$, $N_{II}$ and $N_{Int}$ are the number of points in domain I, domain II and the interface, respectively. 

\subsection{Circle inside a square}
\label{Sec:circ_rectangle}
\begin{figure}[H]
	\centering
	\includegraphics[width = 0.6\textwidth]{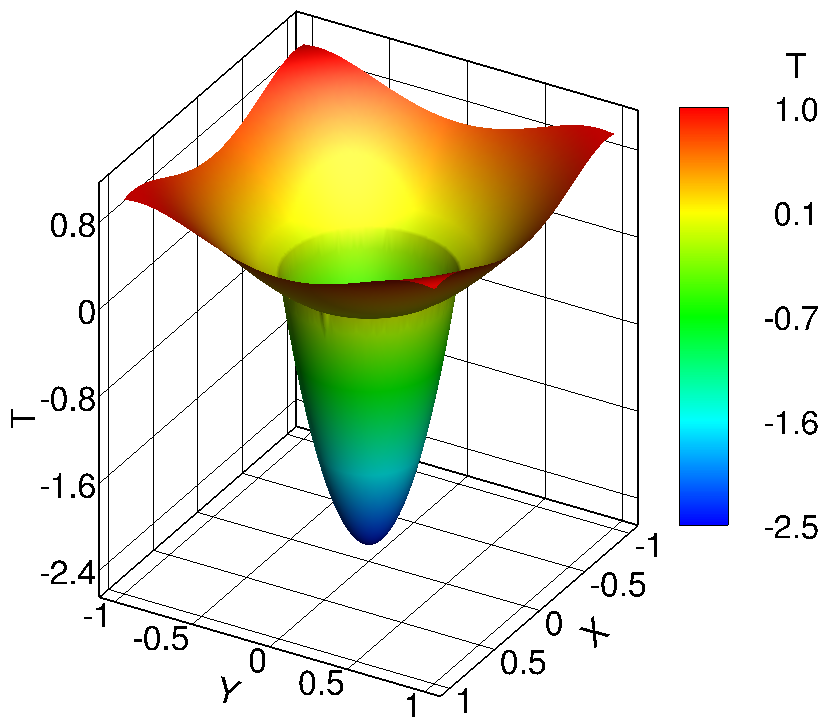}
	\caption{Circle inside a square: Surface plot of the manufactured temperature distribution for circle inside a square and conductivity ratio of 10.}
	\label{Fig:circle_rect_manuf_solution_3D}
\end{figure}

We reconsider the problem of circular domain inside a square as discussed in \cref{Sec:verif}, but now with different thermal conductivities. The manufactured temperature distribution for different subdomains is given by \cref{Eq:circle_rect_manuf_solution} and is shown in \cref{Fig:circle_rect_manuf_solution_3D} for a conductivity ratio $\frac{k_{I}}{k_{II}}$ of 10 . 
\begin{equation}
\begin{aligned}
T_{I} &= k_{II} \; sin(x^{2} + y^{2} - r^{2})\\ 
T_{II} &= k_{I} \; sin(x^{2} + y^{2} - r^{2})
\label{Eq:circle_rect_manuf_solution}
\end{aligned}
\end{equation}

 The manufactured temperature distribution (\cref{Eq:circle_rect_manuf_solution}) is designed such that temperature and heat fluxes are both continuous at the interface points. \Cref{Table:circle_rect_nodes_list} lists the total number of points considered for the error analysis and their distributions within domain I, domain II and the interface.

\begin{table}[H]
	\centering
		\begin{tabular}{|c|c|c|c|c|}
			\hline
			$\Delta \text{r}$ & $N$ & $N_{I}$ & $N_{II}$ & $N_{Int}$ \\ \hline
			0.052 & 1541 & 1221 & 266 & 54 \\ 
			0.037 & 3089 & 2455 & 556 & 78 \\ 
			0.027 & 6033 & 4817 & 1106 & 110 \\ 
			0.019 & 12078 & 9674 & 2247 & 157 \\ \hline
		\end{tabular}%
	\caption{Circle inside a square: Average grid spacing ($\Delta \text{r}$) with the corresponding total number of points distributed within the entire domain ($N$) and subdomains i.e. Domain I ($N_{I}$), Domain II ($N_{II}$) and Interface ($N_{Int}$).}
	\label{Table:circle_rect_nodes_list}
\end{table}

\begin{figure}[H]
  \centering
  \subfigure[]{\includegraphics[width=0.49\textwidth]{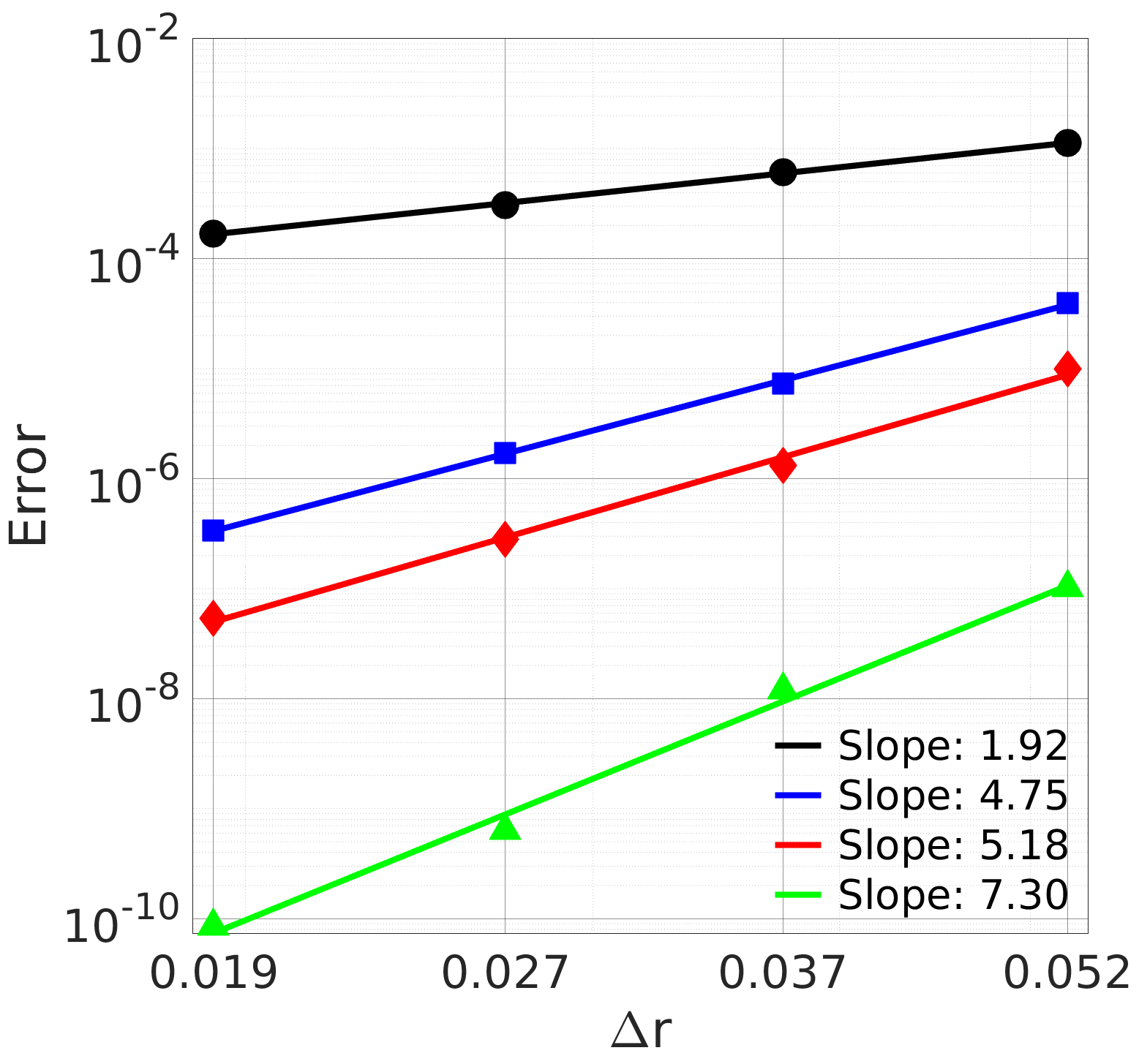}\label{fig:rect_circle_error_k_5_1}}
  \subfigure[]{\includegraphics[width=0.49\textwidth]{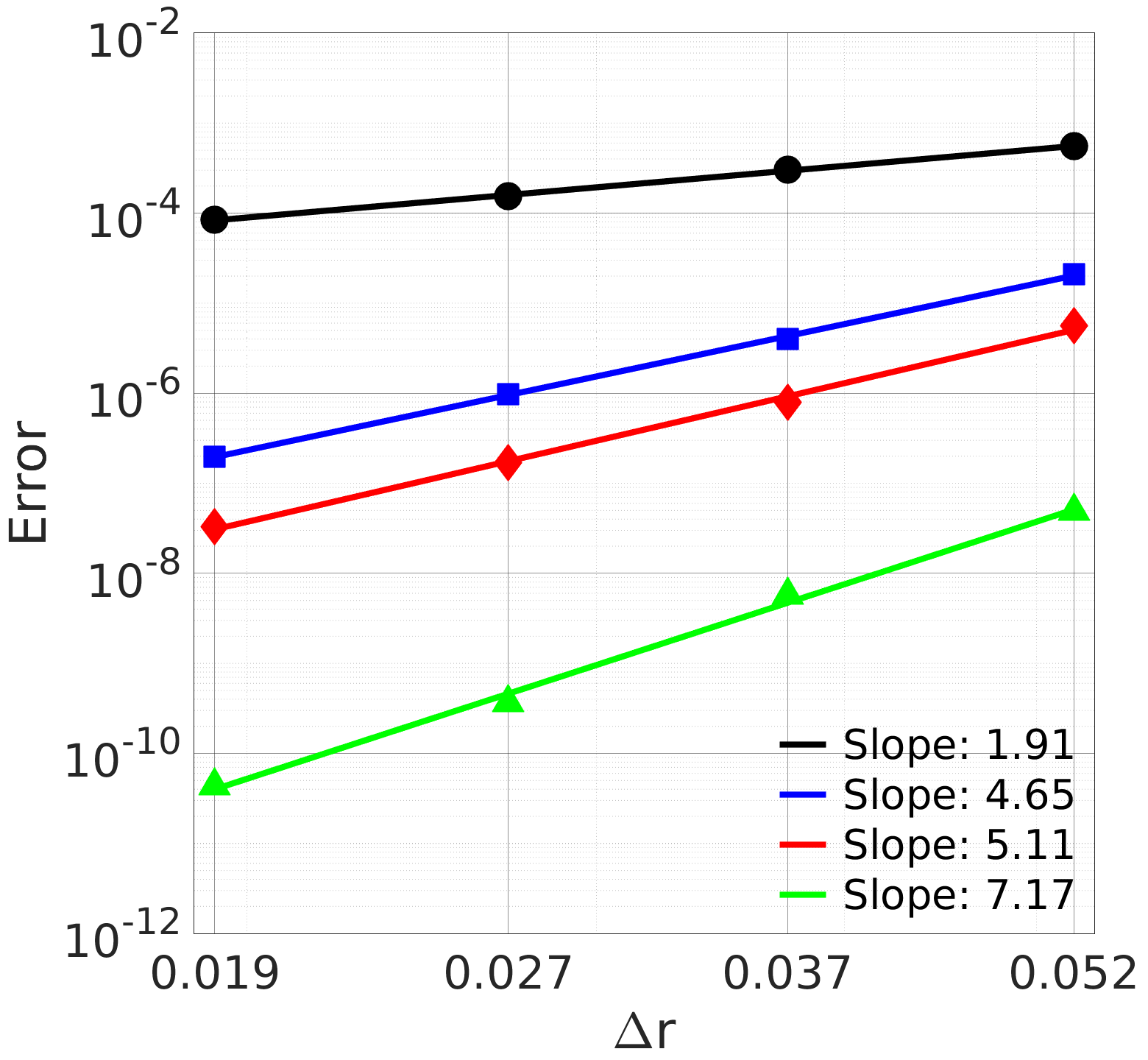}\label{fig:rect_circle_error_k_10_1}}
  \subfigure[]{\includegraphics[width=0.65\textwidth]{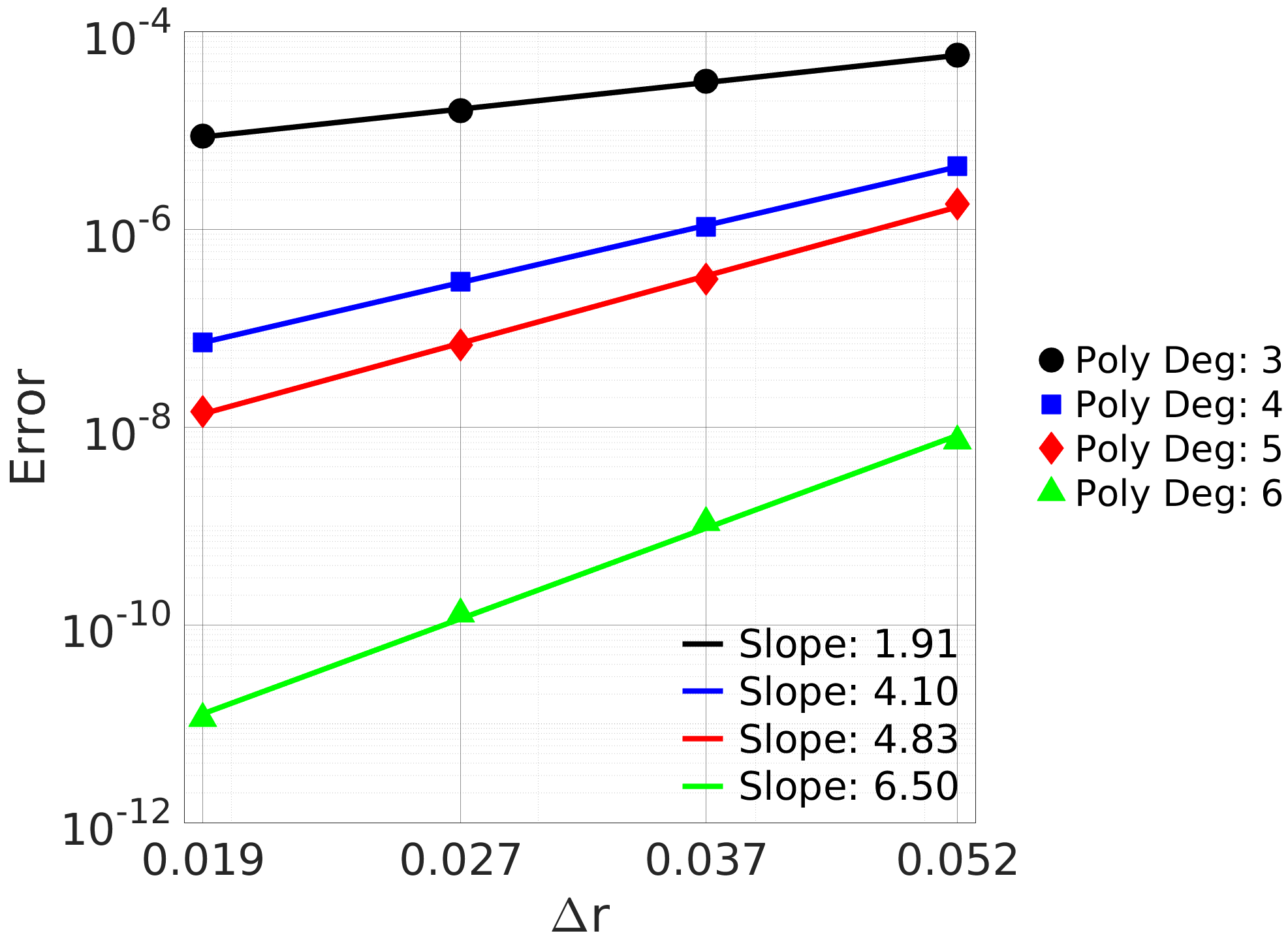}\label{fig:rect_circle_error_k_100_1}}
  \caption{Circle inside a square: Average error vs. $\Delta \text{r}$ in the entire domain for thermal conductivity ratios: (a) 5, (b) 10 and (c) 100.}
  \label{Fig:rect_circle_error_diff_conductivity}
\end{figure}

\begin{figure}[H]
	\centering
	\subfigure[ ]{\includegraphics[width=0.49\textwidth]{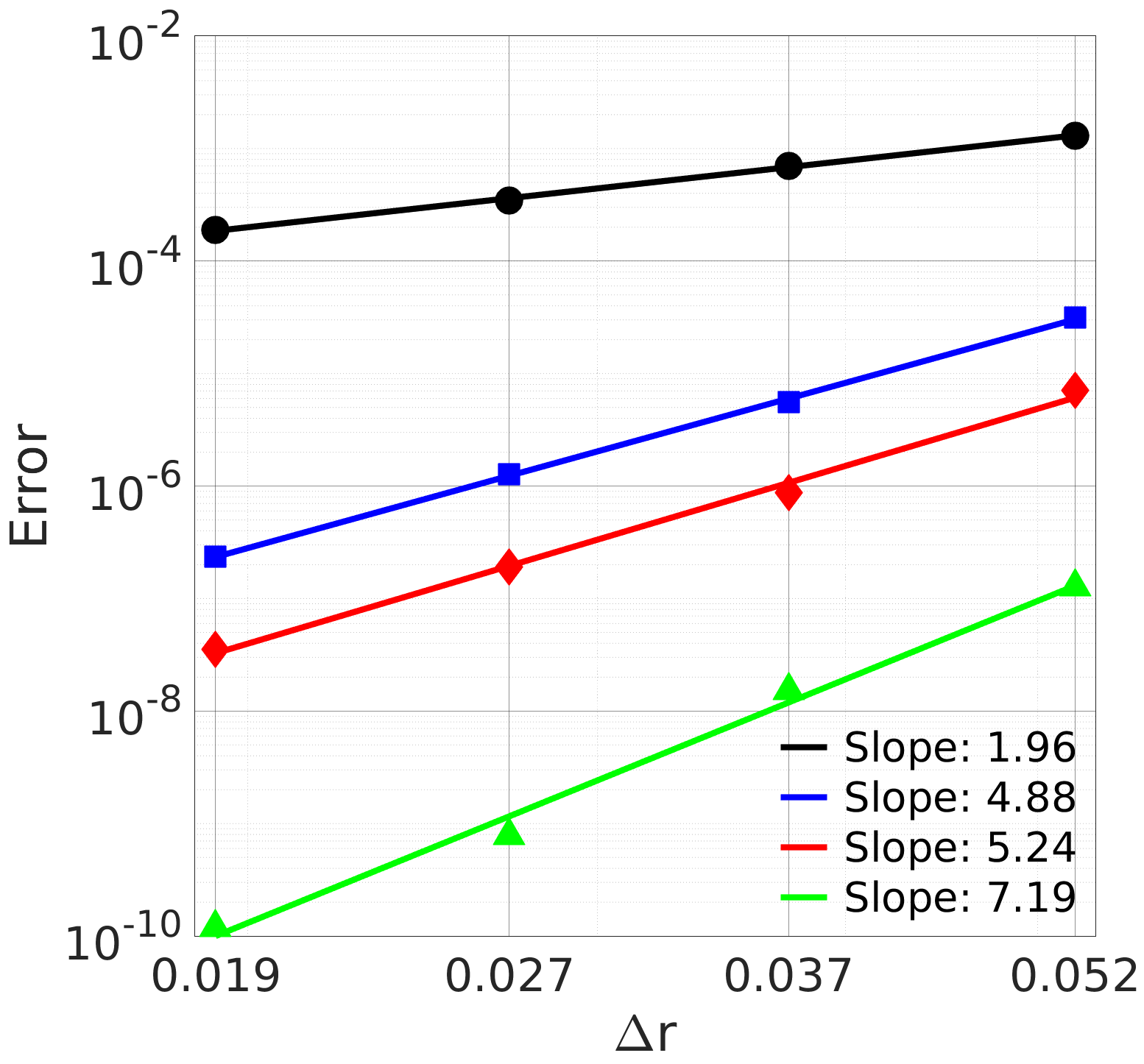}}
	\subfigure[ ]{\includegraphics[width=0.49\textwidth]{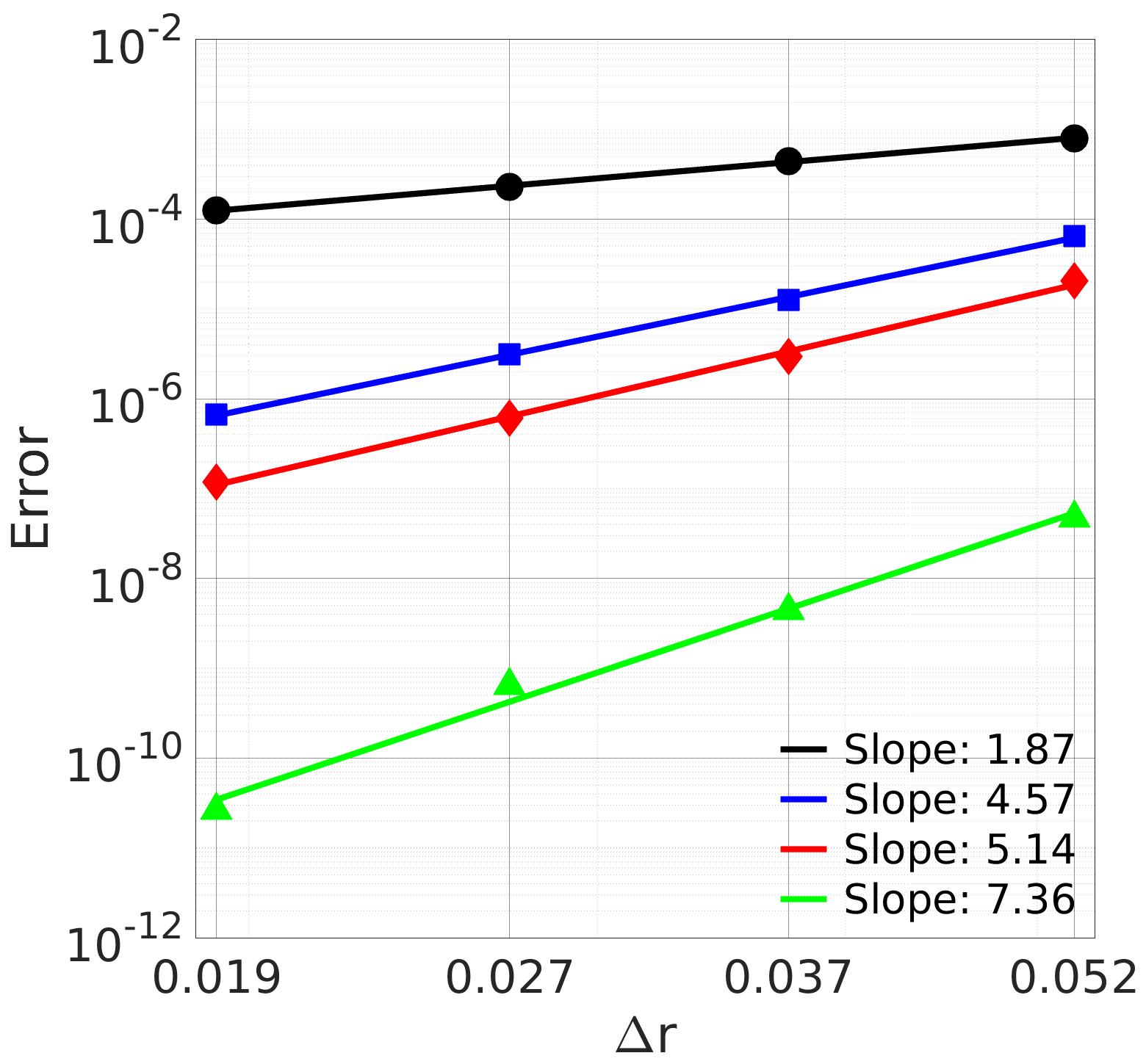}}\\
	\subfigure[ ]{\includegraphics[width=0.65\textwidth]{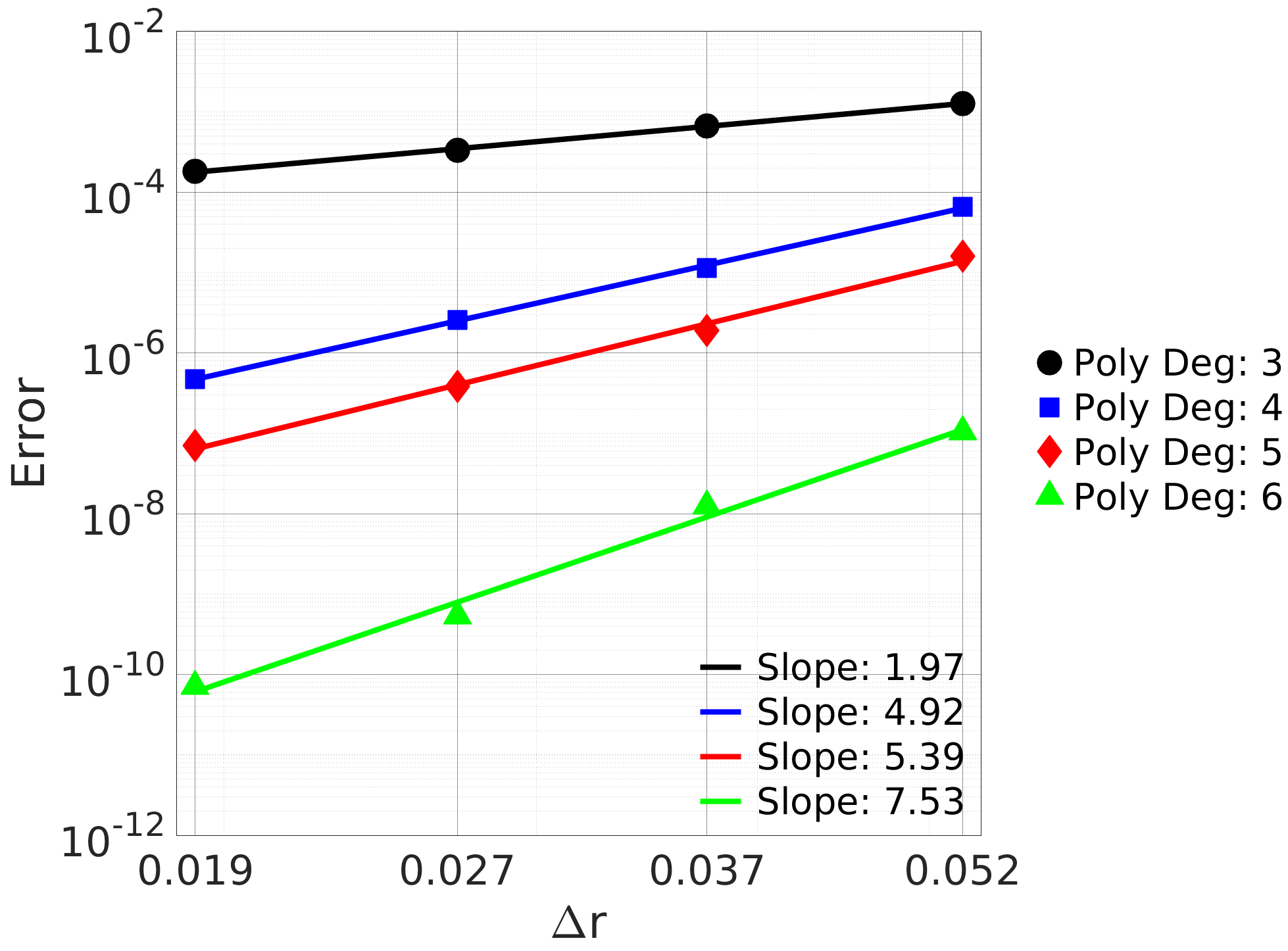}}
	\caption{Circle inside a square: Average error vs. $\Delta \text{r}$ in subdomains for thermal conductivity ratio of 10: (a) domain I, (b) domain II and (c) interface.}
	\label{Fig:rect_circle_error_three_domains}
\end{figure}

\begin{figure}[H]
	\centering
	\subfigure[]{\includegraphics[width=0.45\textwidth]{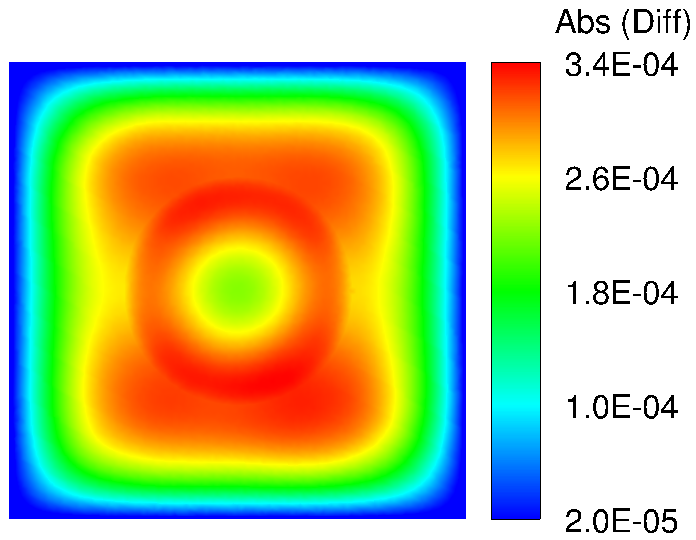}}
	\subfigure[]{\includegraphics[width=0.45\textwidth]{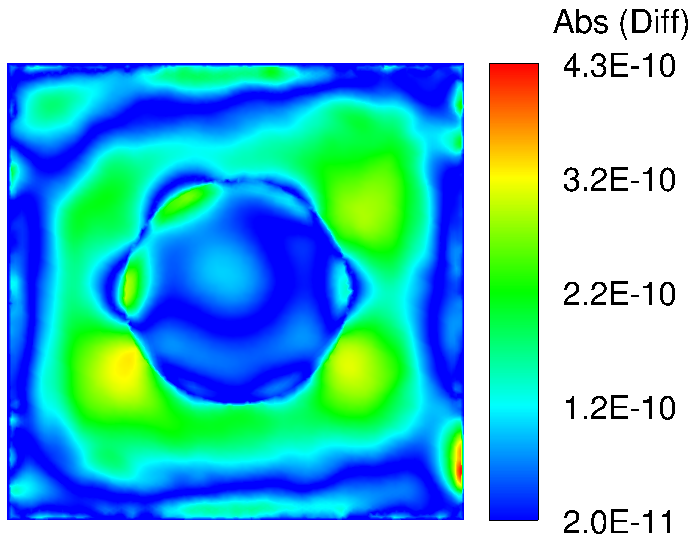}}
	\caption{Circle inside a square: Contours of absolute differences between computed and exact solutions for thermal conductivity ratio of 10 and with 12078 scattered points: (a) polynomial degree = 3 and (b) polynomial degree = 6.}
	\label{Fig:circle_rect_three_domains_error_contours}
\end{figure}

\Cref{Fig:rect_circle_error_diff_conductivity} presents the discretization error as a function of the average grid spacing for the degrees of appended polynomial varied from 3 to 6 for thermal conductivities ratios of 5, 10 and 100. We have presented the average of errors over the three domains including the interface. It can be seen that for each case, the order of convergence has a minimum value of $k-1$ for a polynomial degree of $k$, even for a large conductivity ratio of 100. The rate of convergence increases from 1.91 for polynomial degree of 3 to 6.50 for polynomial degree of 6. \Cref{Fig:rect_circle_error_three_domains} shows the variation of the errors in individual domains for a thermal conductivity ratio of 10. It can be seen that all three domains (including the interface) display the expected rate of convergence for a given degree of the appended polynomial.  It can be seen from \cref{Fig:circle_rect_three_domains_error_contours}  that the errors (Abs (Diff) = $|T_{exact} - T_{calc}|$) reduce very fast with increase in degree of appended polynomial, establishing the high accuracy of the method.

 \newpage
 
\subsection{Astroidal domain inside a square}
\label{Sec:astroid_rect}
\begin{figure}[H]
	\centering
	\includegraphics[width = 0.55\textwidth]{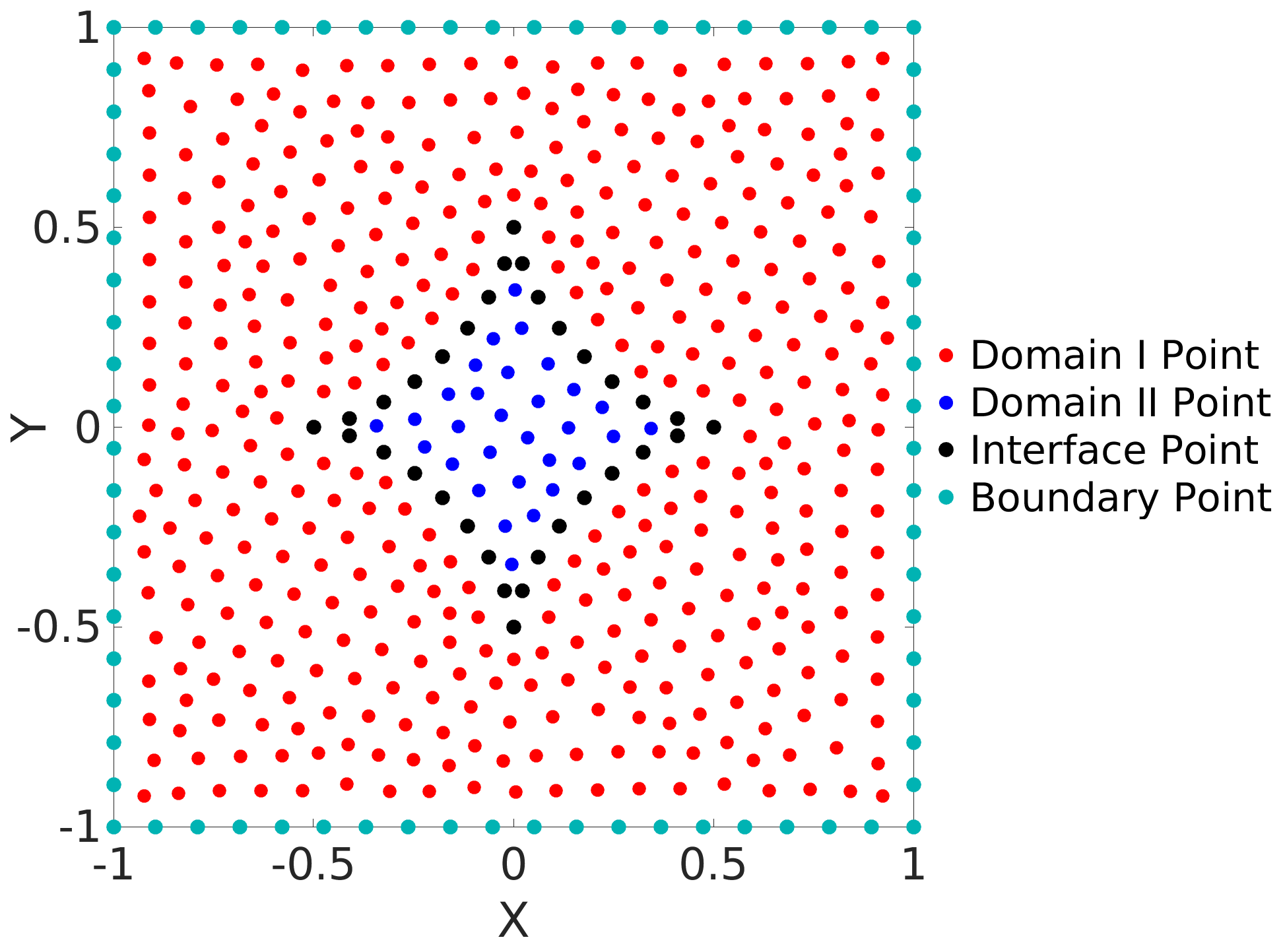}
	\caption{Astroidal domain inside a square: Locations of 540 points within the domain.}
	\label{Fig:astroid_in_rect_nodes}
\end{figure}

The second problem considered for accuracy assessment is a geometry with sharp corners. Modelling interfaces with sharp corners remains a challenge as most of the traditional CFD techniques such as finite difference (FDM), finite volume (FVM), boundary element (BEM) methods, etc may require high grid resolution to give accurate results. To test the applicability of PHS-RBF for such problems, we have considered a 2D astroidal shaped domain inside a square section. \Cref{Fig:astroid_in_rect_nodes} shows the positions of 540 points within the domain.

\begin{figure}[H]
	\centering
	\includegraphics[width = 0.6\textwidth]{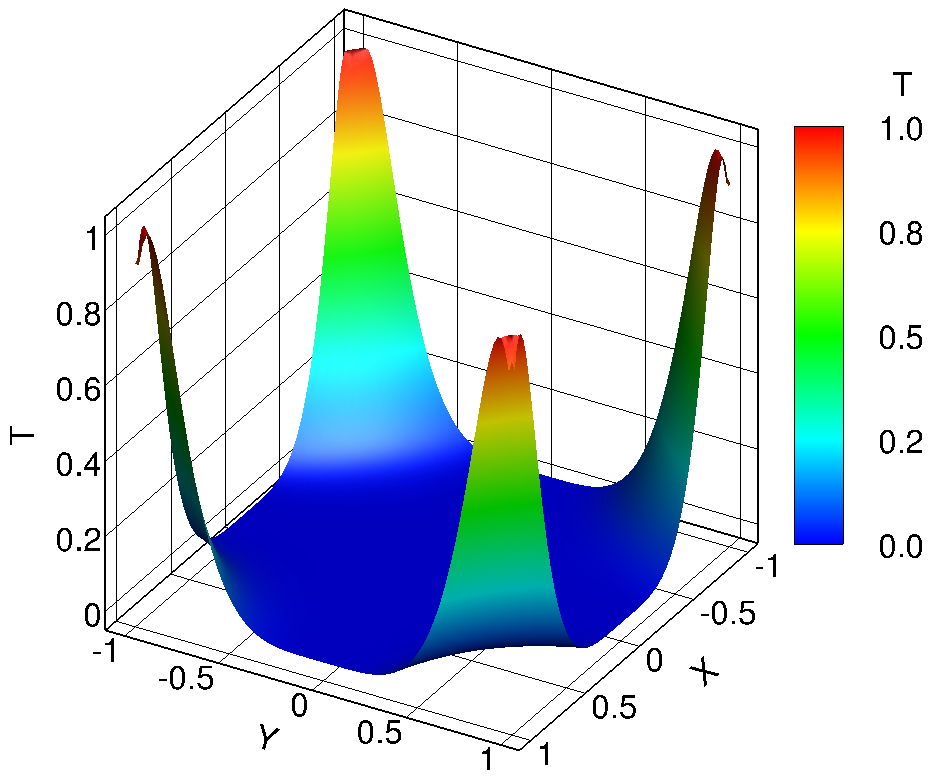}
	\caption{Astroidal domain inside a square: Surface plot of the manufactured temperature distribution for the astroidal geometry inside a square.}
	\label{Fig:astroid_in_rect_manuf_dist_3D}
\end{figure}

The external boundary points are represented by cyan coloured dots and are part of domain I. The parametric equation for the interface is given by
\begin{equation}
\begin{aligned}
x = a \; cos^{3}\theta \\
y = a \; sin^{3}\theta
\end{aligned}
\hspace{1cm} 0\leq\theta\leq 2\pi;\hspace{0.2cm}a = 0.5
\label{Eq:astroid_rect_parametric}
\end{equation} 
 The manufactured temperature distribution given by \cref{Eq:astroid_in_rect_manufactured_temp_distribution} satisfies the flux conditions at the interface. A surface plot of the distribution is shown in \cref{Fig:astroid_in_rect_manuf_dist_3D} for thermal conductivity ratio of 10 i.e. $k_{I}$ =  10 and $k_{II}$ = 1.
\begin{equation}
\begin{aligned}
T_{I} &= k_{II} \; sin\bigg(x^4 y^4 \bigg(\bigg(\frac{x}{a}\bigg)^\frac{2}{3} + \bigg(\frac{y}{a}\bigg)^\frac{2}{3} - 1\bigg)\bigg) \\ 
T_{II} &= k_{I} \; sin\bigg(x^4 y^4 \bigg(\bigg(\frac{x}{a}\bigg)^\frac{2}{3} + \bigg(\frac{y}{a}\bigg)^\frac{2}{3} - 1\bigg)\bigg) 
\end{aligned}
\label{Eq:astroid_in_rect_manufactured_temp_distribution}
\end{equation}

\begin{figure}[H]
	\centering
	\subfigure[]{\includegraphics[width=0.49\textwidth]{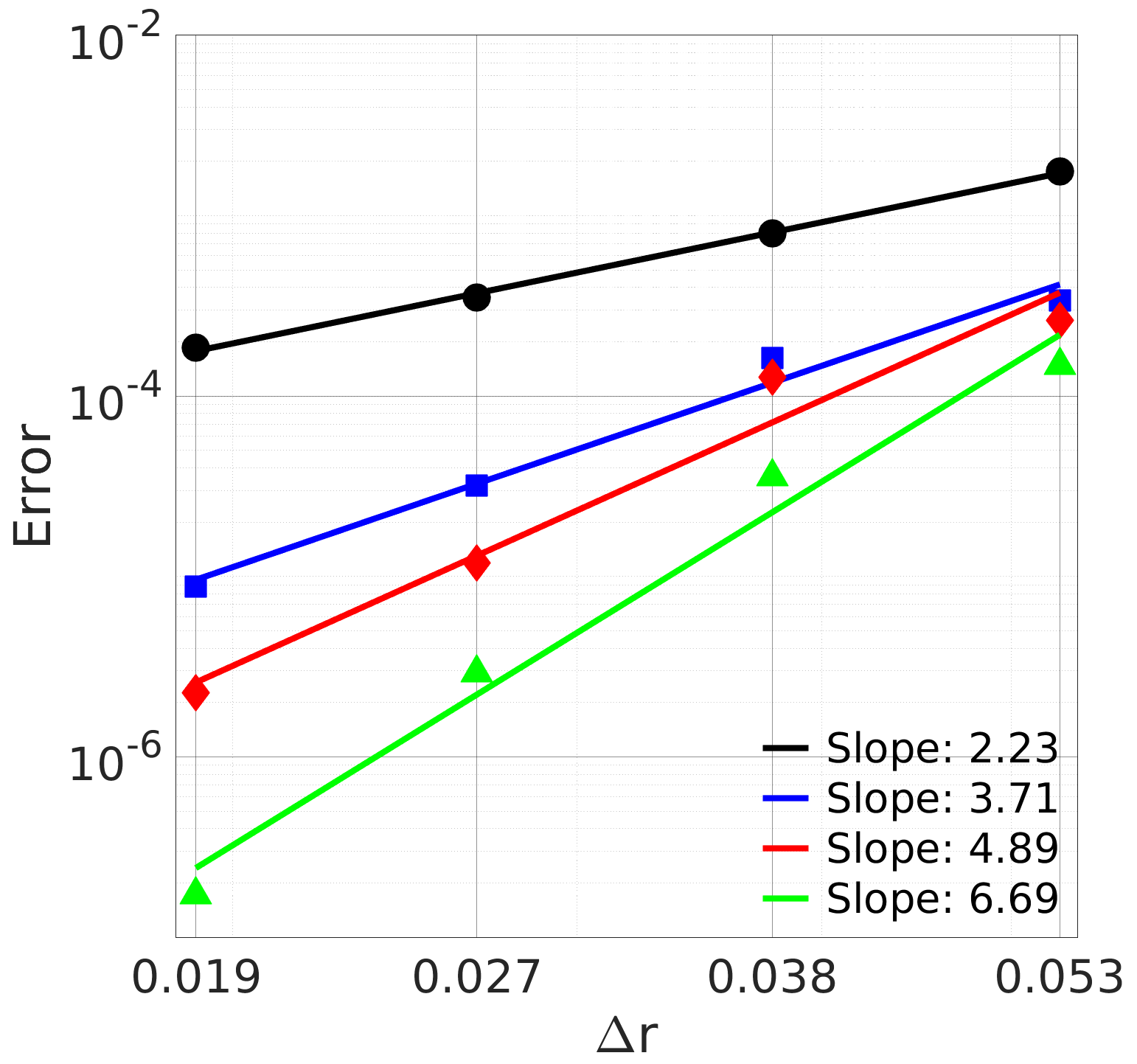}\label{fig:astroid_in_rect_error_full_domain_k_5_1}}
	\subfigure[]{\includegraphics[width=0.49\textwidth]{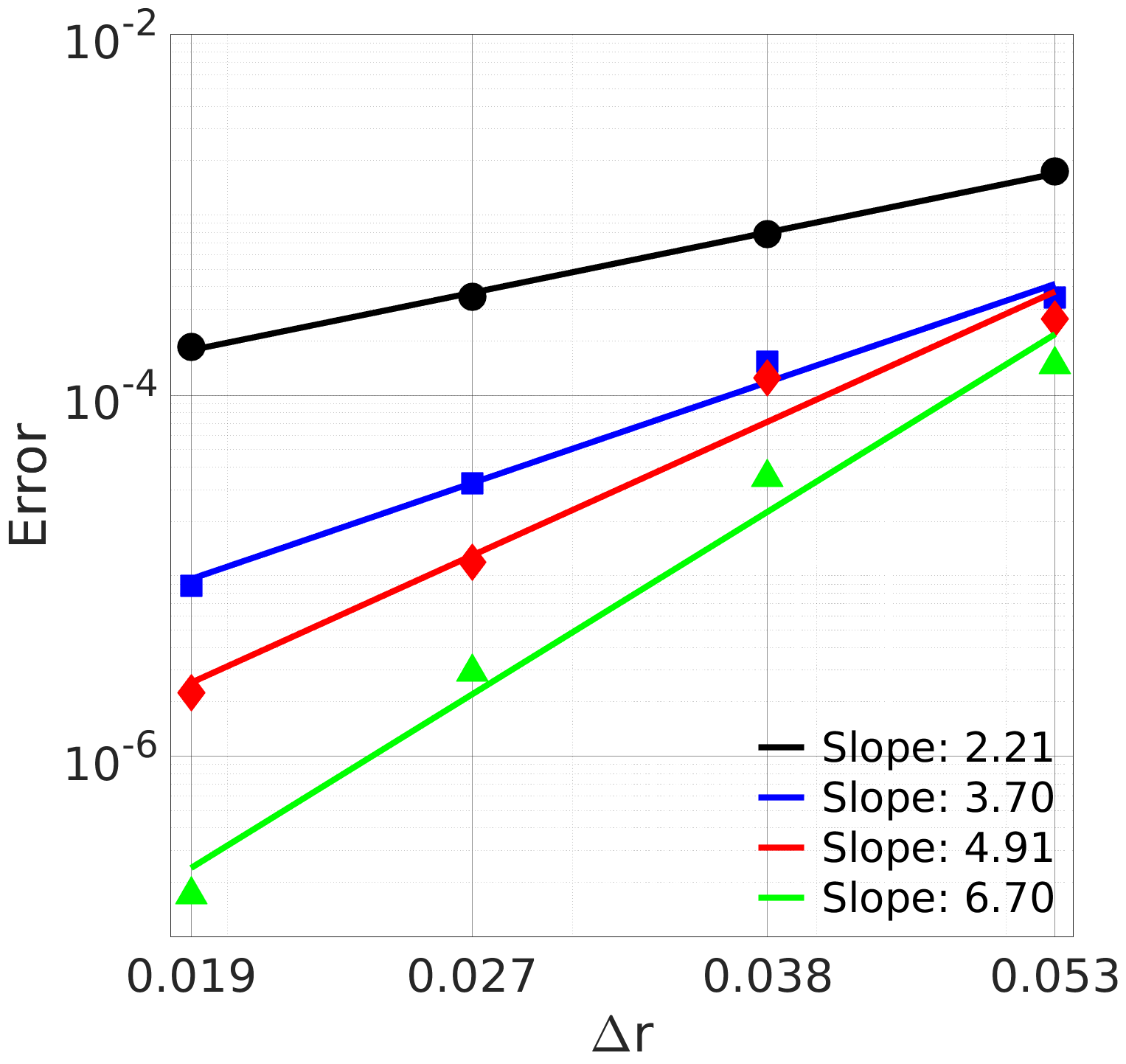}\label{fig:astroid_in_rect_error_full_domain_k_10_1}}
	\subfigure[]{\includegraphics[width=0.65\textwidth]{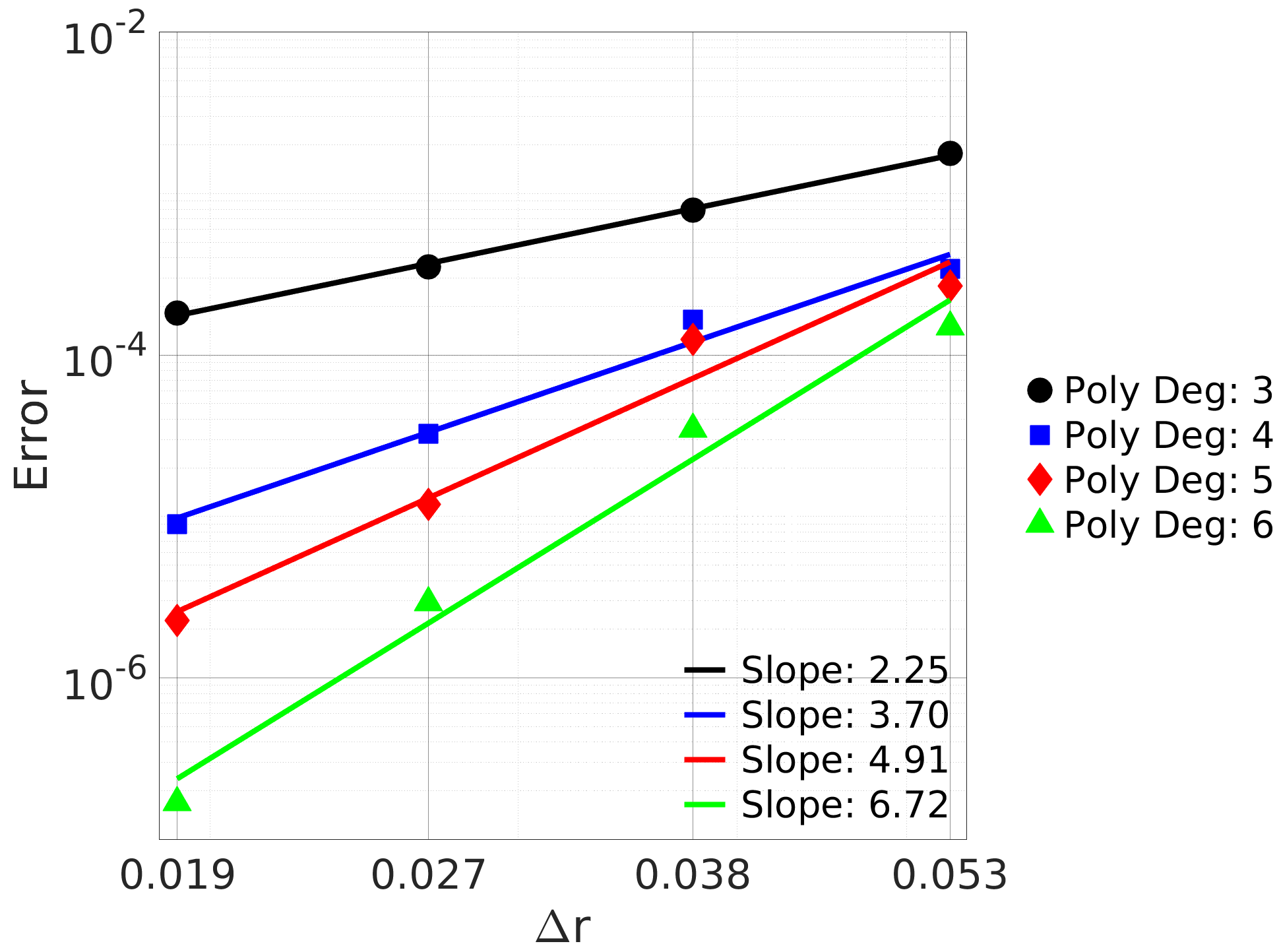}\label{fig:astroid_in_rect_error_full_domain_k_100_1}}
	\caption{Astroidal domain inside a square: Average error vs. $\Delta \text{r}$ in the entire domain for thermal conductivity ratios: (a) 5, (b) 10 and (c) 100.}
	\label{Fig:astroid_in_rect_error_full_domain}
\end{figure}

\Cref{Fig:astroid_in_rect_error_full_domain} plots the variation of the average absolute difference between the computed and the exact temperatures as a function of the average point spacing for four different distribution of data points within each subdomain (\Cref{Table:astroid_rect_nodes_list}). \Cref{fig:astroid_in_rect_error_full_domain_k_5_1,fig:astroid_in_rect_error_full_domain_k_10_1,fig:astroid_in_rect_error_full_domain_k_100_1} are for conductivity ratios of 5, 10 and 100. It can be noted from \cref{Fig:astroid_in_rect_error_full_domain} that error decreases as expected with the degree of the appended polynomial, even for the conductivity ratio of 100. The conductivity ratio does not appear to influence the accuracy. This observation is very encouraging compared to a single domain method which tends to smear the sharp discontinuity and produce erroneous results.

\begin{table}[H]
	\centering
	\begin{tabular}{|c|c|c|c|c|}
		\hline
		$\Delta \text{r}$ & $N$ & $N_{I}$ & $N_{II}$ & $N_{Int}$ \\ \hline
		0.053 & 1523 & 1396 & 79  & 48 \\ 
		0.038 & 3038 & 2769 & 197 & 72 \\ 
		0.027 & 6084 & 5568 & 412 & 104 \\ 
		0.019 & 12060& 11080& 832 & 148 \\ \hline
	\end{tabular}%
	\caption{Astroidal domain inside a square: Average grid spacing ($\Delta \text{r}$) with the corresponding total number of points distributed within the entire domain ($N$) and subdomains i.e. Domain I ($N_{I}$), Domain II ($N_{II}$) and Interface ($N_{Int}$).}
	\label{Table:astroid_rect_nodes_list}
\end{table}

\begin{figure}[H]
	\centering
	\includegraphics[width=0.65\textwidth]{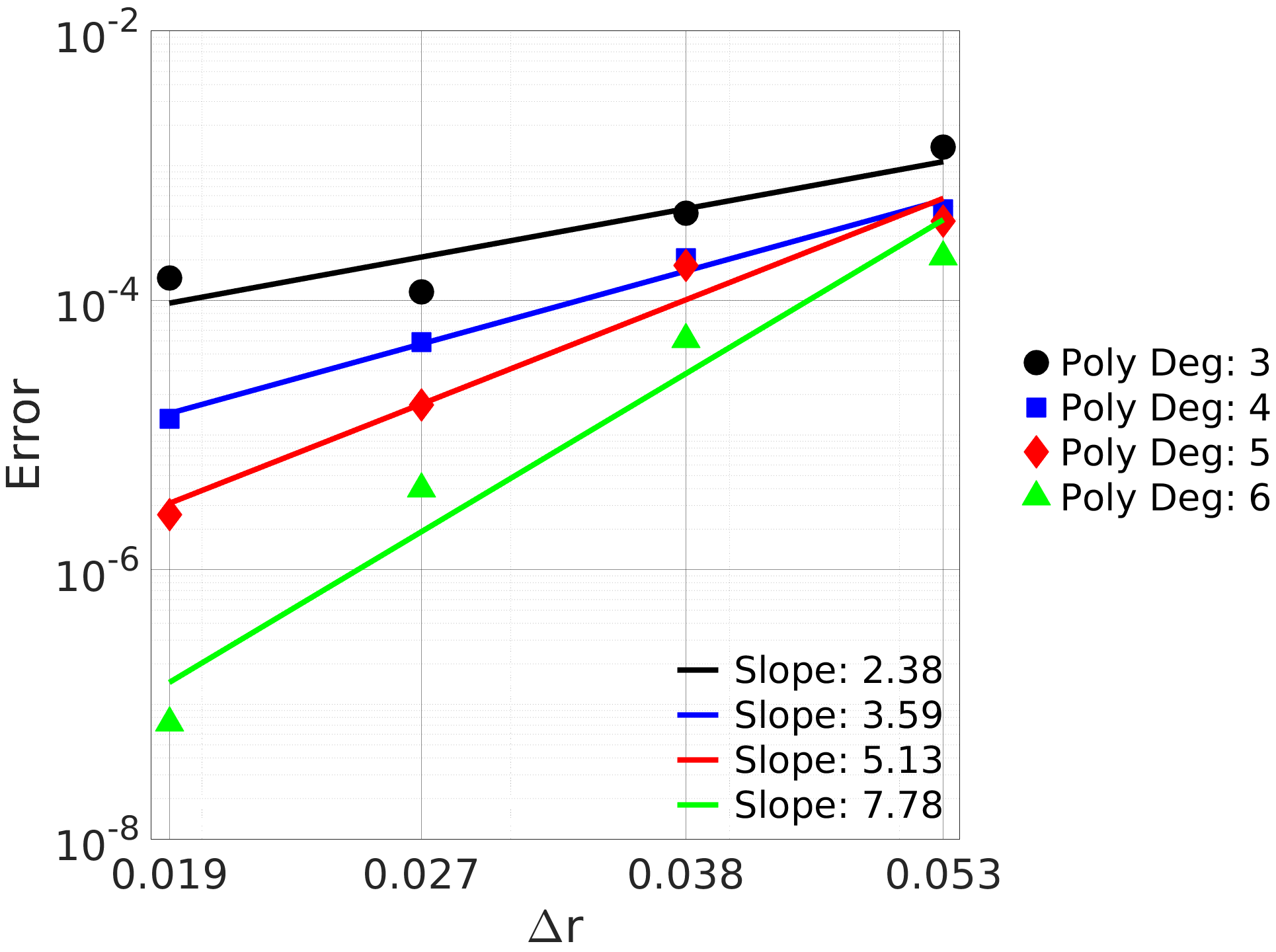}
	\caption{Astroidal domain inside a square: Variation of error at interface points for conductivity ratio of 10.}
	\label{Fig:astroid_in_rect_error_interface_k_10_1}
\end{figure}

\begin{figure}[H]
	\centering
	\subfigure[]{\includegraphics[width=0.47\textwidth]{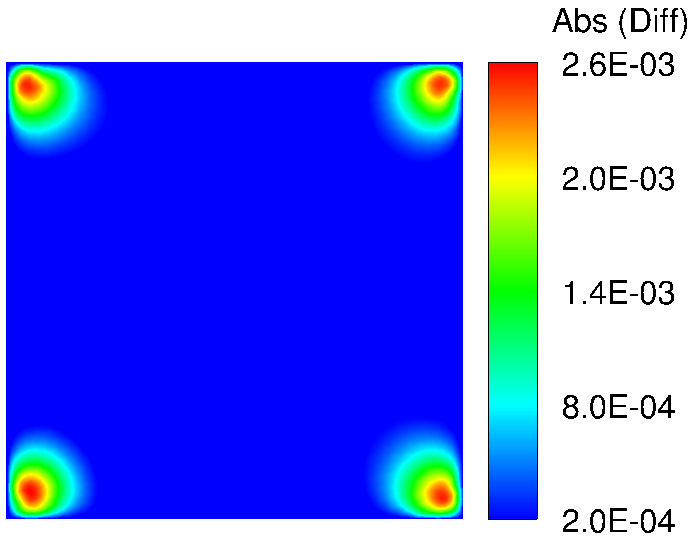}}\hspace{0.3cm}
	\subfigure[]{\includegraphics[width=0.47\textwidth]{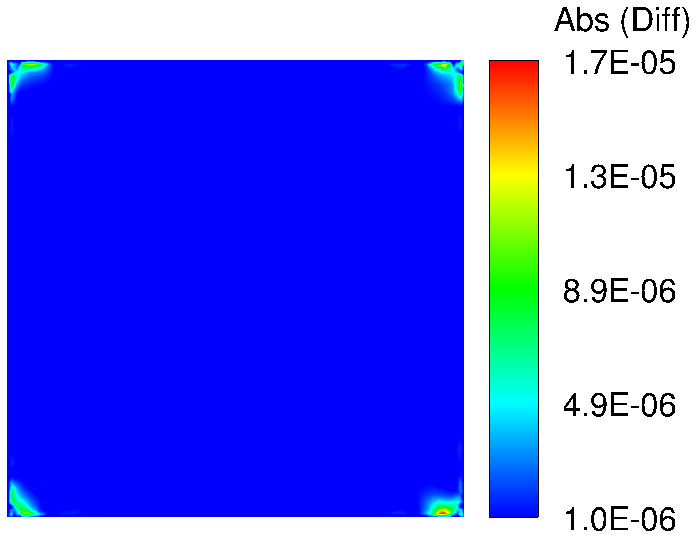}}
	\caption{Astroidal domain inside a square: Contours of absolute differences between computed and exact solutions for thermal conductivity ratio of 10 and with 12060 scattered points: (a) polynomial degree = 3 and (b) polynomial degree = 6.}
	\label{Fig:astroid_in_rect_three_domains_error_contours}
\end{figure}

In \cref{Fig:astroid_in_rect_error_interface_k_10_1}, we plot the averaged error only for the interface points. As stated earlier, the governing equation is not satisfied at the interface points. Only the flux is ensured to be continuous. Hence the interface is prone to larger errors than the domains where the governing equation is satisfied. However, we observe that even for the interface points, taken in isolation, the error is small and decreases as expected with the degree of the appended polynomial. The same conclusion can be drawn from the error contours (Abs (Diff) = $|T_{exact} - T_{calc}|$) shown in \cref{Fig:astroid_in_rect_three_domains_error_contours} with 12060 scattered points for polynomial degrees $3$ and $6$. Interestingly, the errors are largest near the four corners of the external boundary and not at the interface. 

\subsection{Sphere inside a cube}
\label{Sec:sphere_cuboid}

\begin{figure}[H]
	\centering
	\includegraphics[width = 0.55\textwidth]{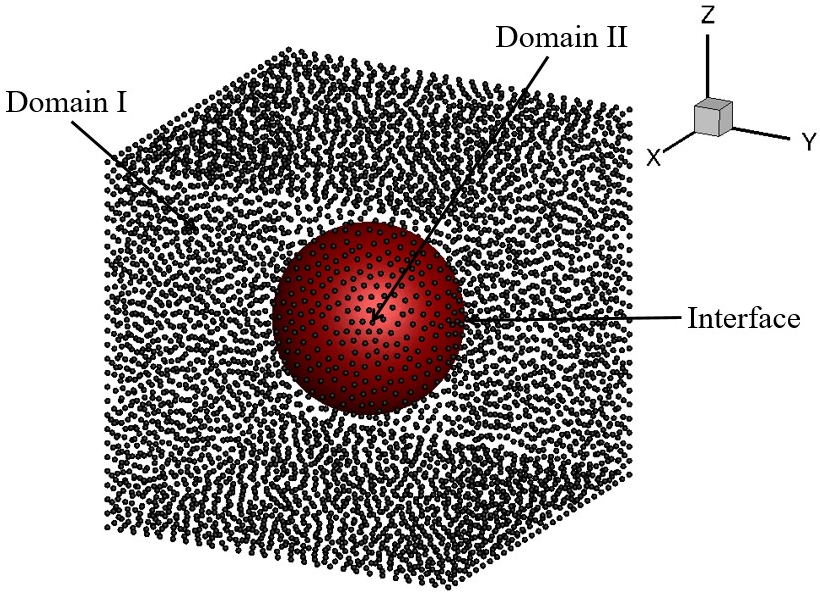}
	\caption{Sphere inside a cube: Illustrative point layout for sphere inside a cube.}
	\label{Fig:sphere_cuboid_domains}
\end{figure}

We next consider a three-dimensional problem of a sphere inside a cube. The sphere is of radius 0.5 units and the cube is 2 units on its side. The solution domain consists of three sets of points. Domain I is the region between the surface of the sphere and inside the cube. Domain II is the set of points inside the sphere. The interface is separately discretized by its own scattered points.  \Cref{Fig:sphere_cuboid_domains} shows the layout of the scattered points.

\begin{figure}[H]
	\centering
	\includegraphics[width = 0.6\textwidth]{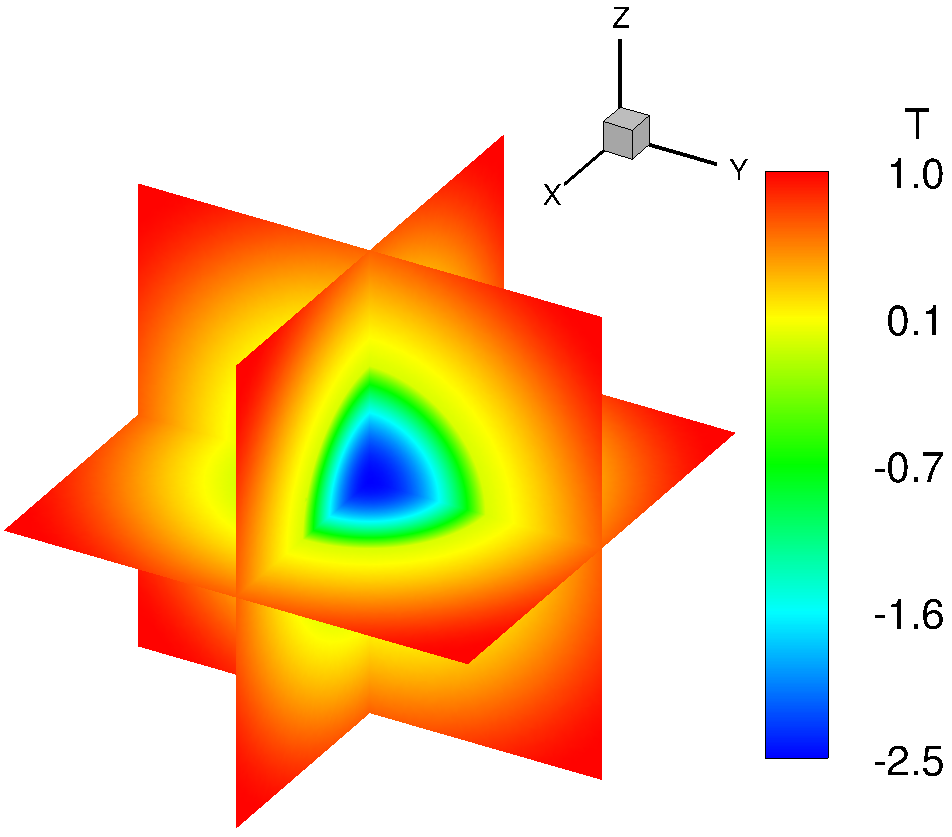}
	\caption{Sphere inside a cube: 3D view of the manufactured temperature distribution.}
	\label{Fig:sphere_cuboid_manufactured}
\end{figure}

 A manufactured temperature distribution given by \cref{Eq:manuf_solution_sphere_cuboid} and illustrated in \cref{Fig:sphere_cuboid_manufactured} for conductivity ratio of 10 is used for accuracy estimation. 

\begin{equation}
\begin{aligned}
T_{I} &= k_{II} \; sin(x^{2} + y^{2} + z^{2} - r^{2})\\ 
T_{II} &= k_{I} \; sin(x^{2} + y^{2} + z^{2} - r^{2})
\end{aligned}
\label{Eq:manuf_solution_sphere_cuboid}
\end{equation}

Four sets of scattered points consisting of 28961, 87680, 299798 and 711041 total points with average spacing of 0.057, 0.04, 0.026 and 0.018 units are considered for polynomial degree from 3 to 6. The number of points in different regions is given in \cref{Table:sphere_cuboid_nodes_list}.

\begin{table}[H]
	\centering
	\begin{tabular}{|c|c|c|c|c|}
		\hline
		$\Delta \text{r}$ & $N$ & $N_{I}$ & $N_{II}$ & $N_{Int}$ \\ \hline
		0.057 & 28961 & 14466 & 10326 & 4169 \\ 
		0.040 & 87680 & 50610 & 28826 & 8244 \\ 
		0.026 & 299798 & 170789 & 107603 & 21406 \\ 
		0.018 & 711041  & 379939  & 293838  & 37264 \\ \hline
	\end{tabular}%
	\caption{Sphere inside a cube: Average grid spacing ($\Delta \text{r}$) with the corresponding total number of points distributed within the entire domain ($N$) and subdomains i.e. Domain I ($N_{I}$), Domain II ($N_{II}$) and Interface ($N_{Int}$).}
	\label{Table:sphere_cuboid_nodes_list}
\end{table}
\newpage

\begin{figure}[H]
	\centering
	\subfigure[]{\includegraphics[width=0.49\textwidth]{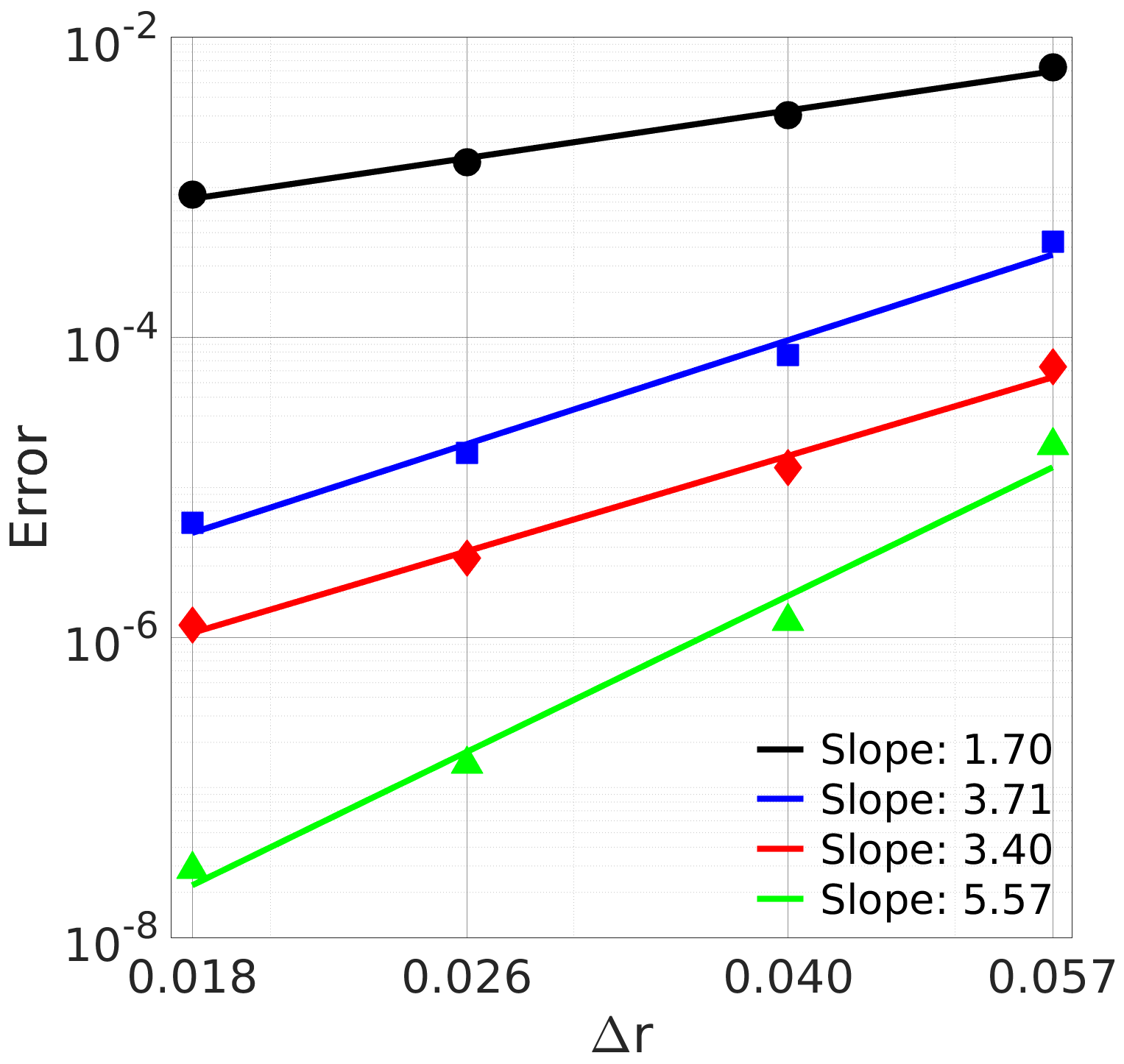}\label{fig:sphere_cuboid_error_full_domain_k_5_1}}
	\subfigure[]{\includegraphics[width=0.49\textwidth]{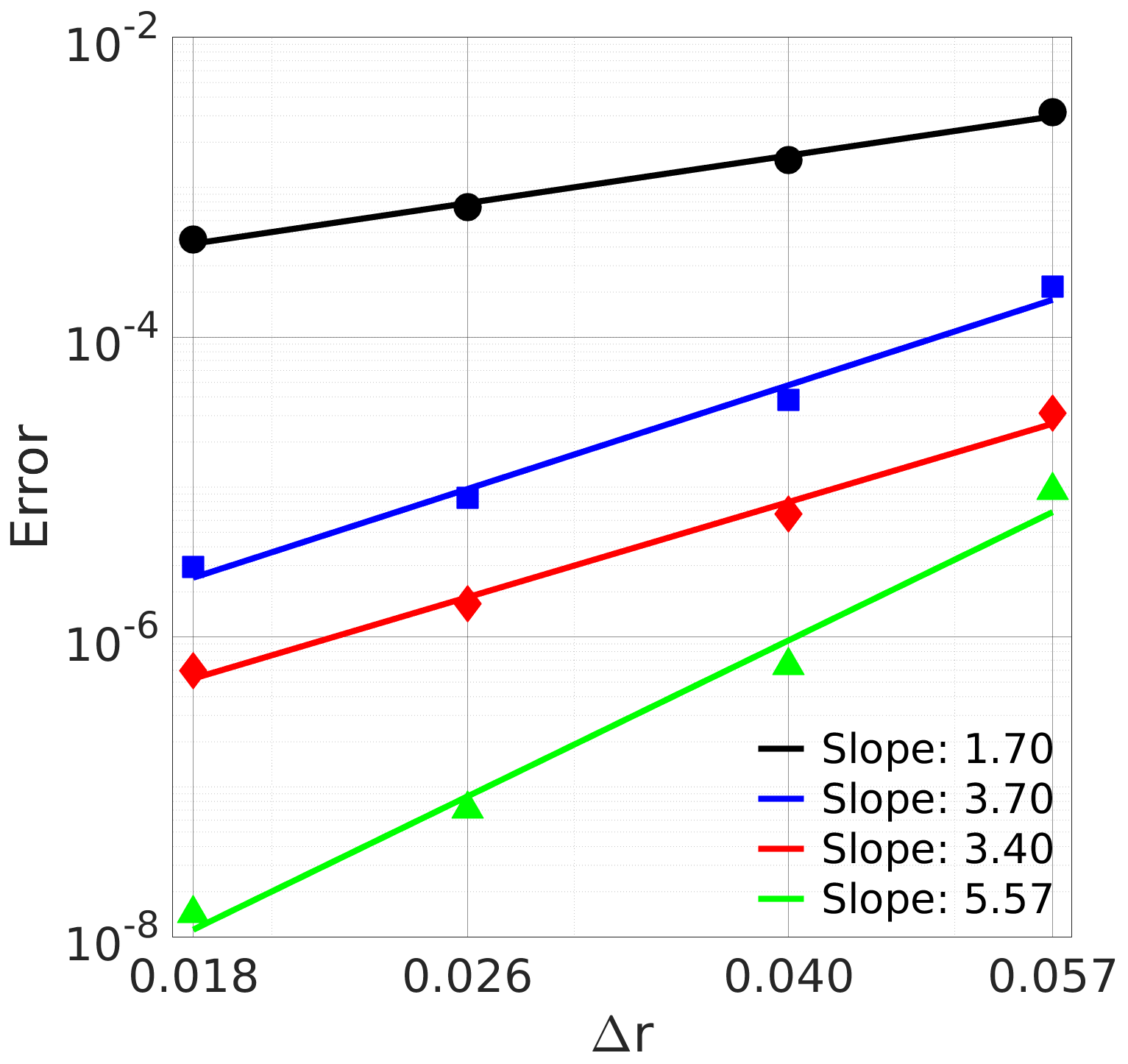}\label{fig:sphere_cuboid_error_full_domain_k_10_1}}
	\subfigure[]{\includegraphics[width=0.65\textwidth]{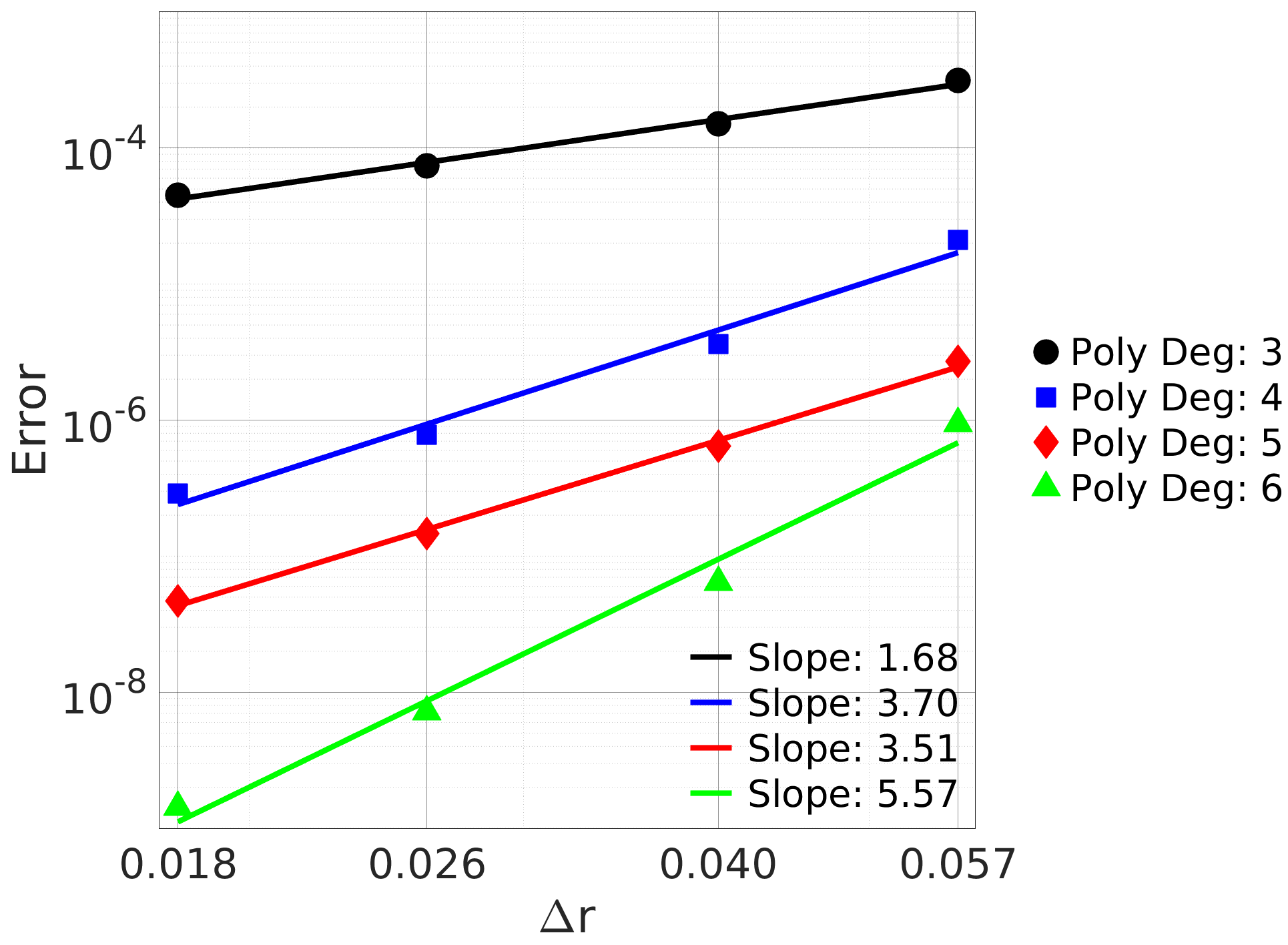}\label{fig:sphere_cuboid_error_full_domain_k_100_1}}
    \caption{Sphere inside a cube: Average error vs. $\Delta \text{r}$ in the entire domain for thermal conductivity ratios: (a) 5, (b) 10 and (c) 100.}
	\label{Fig:sphere_cuboid_error_full_domain}
\end{figure}

\begin{figure}[H]
	\centering
	\subfigure[]{\includegraphics[width=0.49\textwidth]{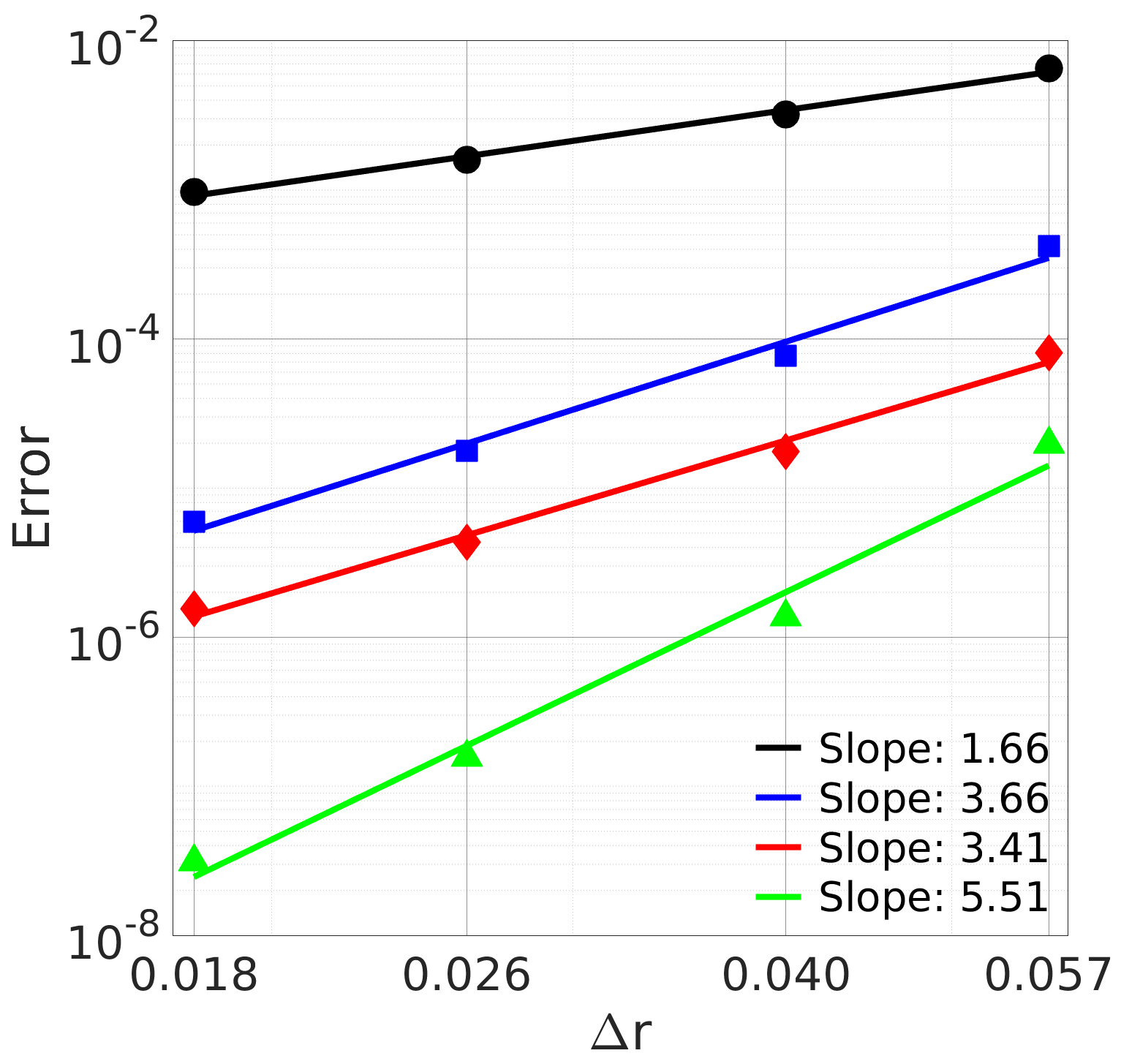}}
	\subfigure[]{\includegraphics[width=0.49\textwidth]{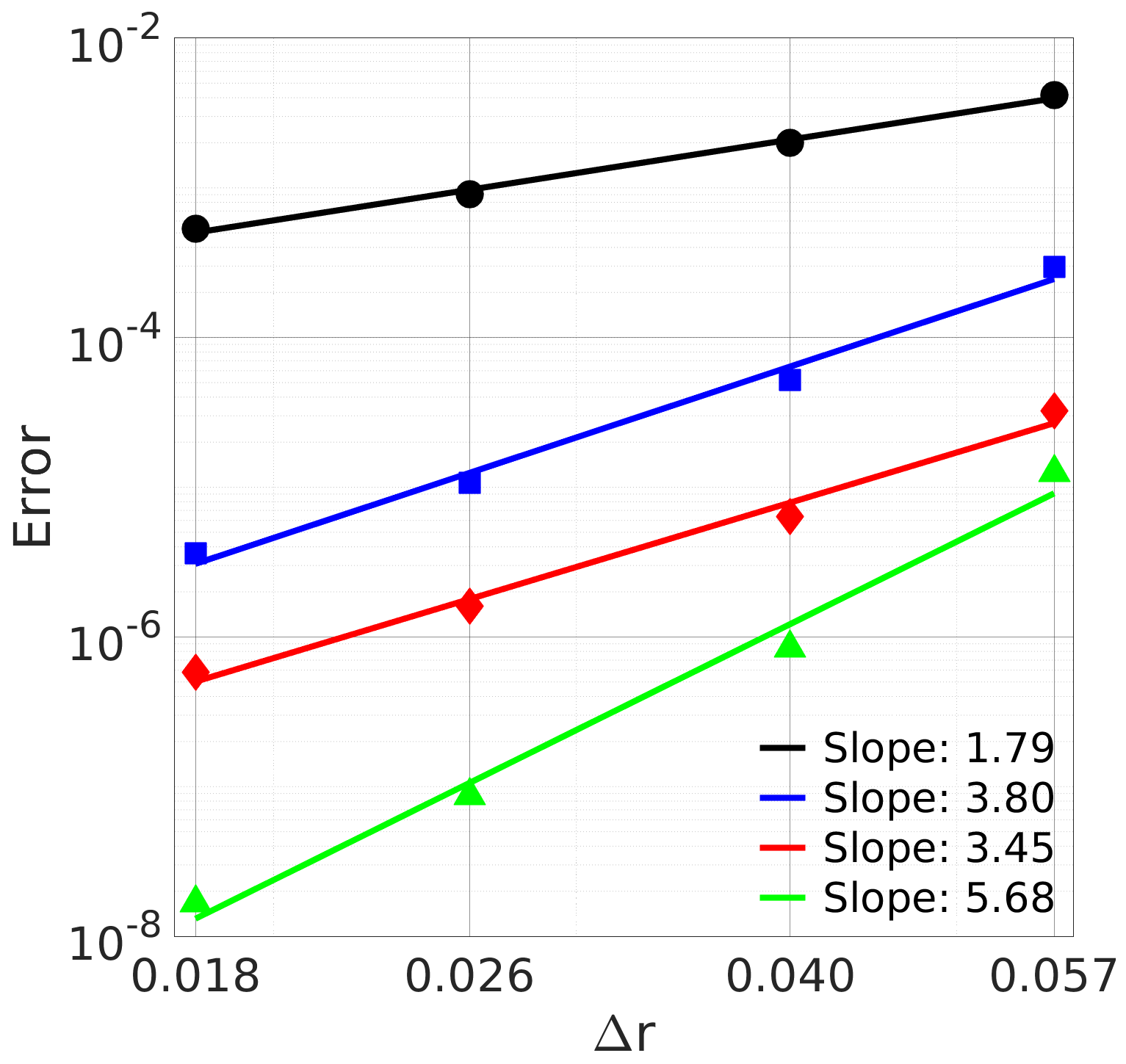}}\\
	\subfigure[]{\includegraphics[width=0.65\textwidth]{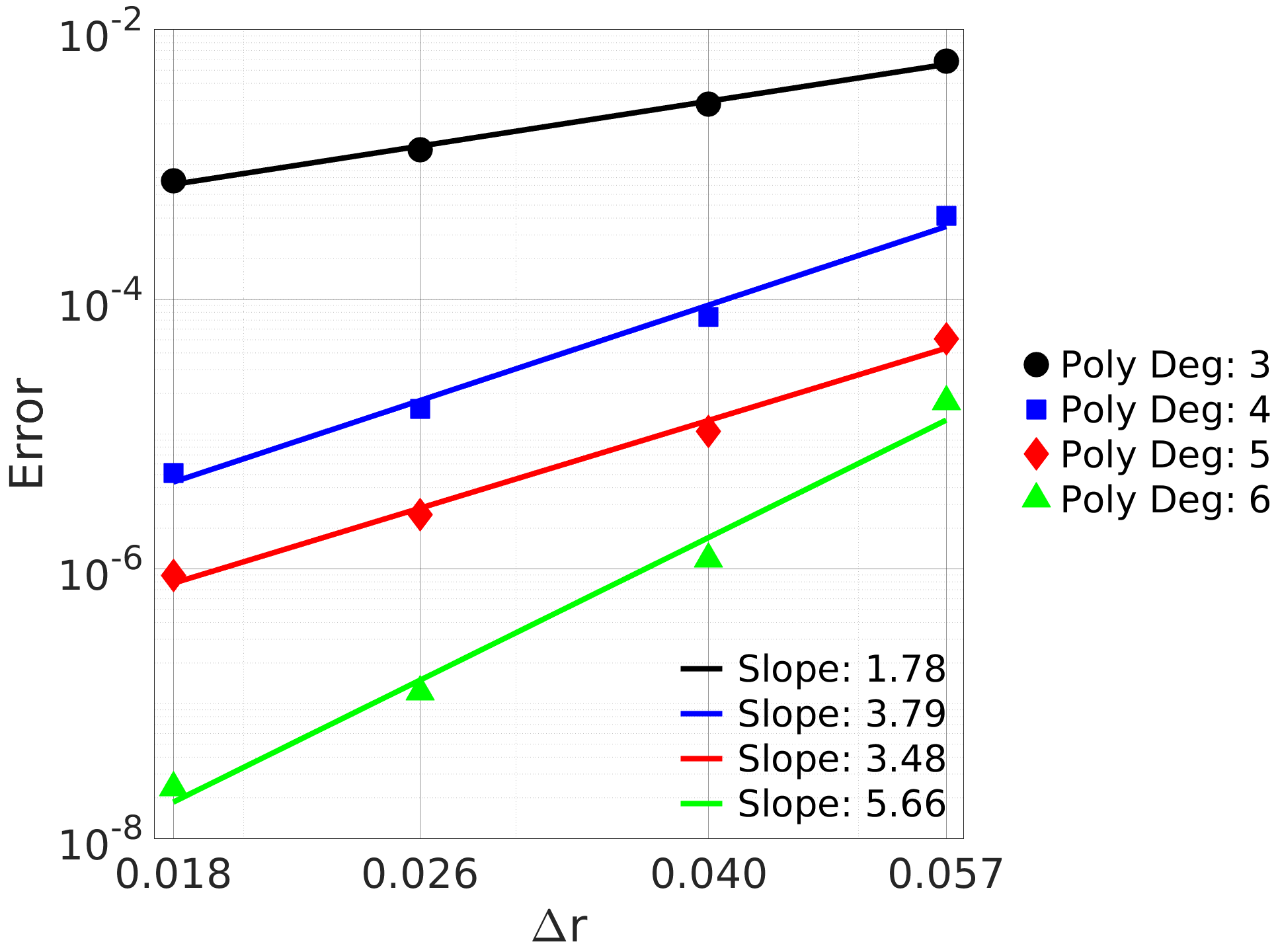}\label{fig:sphere_cuboid_interface_error}}
    \caption{Sphere inside a cube: Average error vs. $\Delta \text{r}$ in subdomains for thermal conductivity ratio of 10: (a) domain I, (b) domain II and (c) interface.}
	\label{Fig:sphere_cuboid_three_domains}
\end{figure}
\begin{figure}[H]
	\centering
	\subfigure[]{\includegraphics[width=0.49\textwidth]{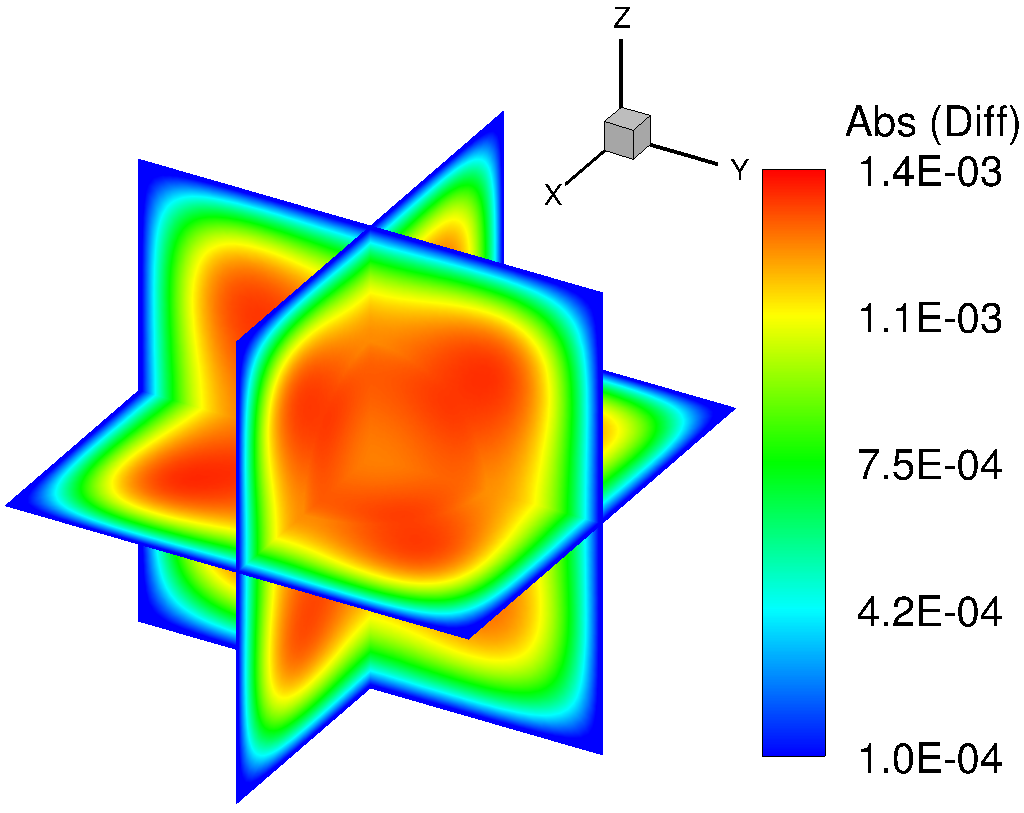}}
	\subfigure[]{\includegraphics[width=0.49\textwidth]{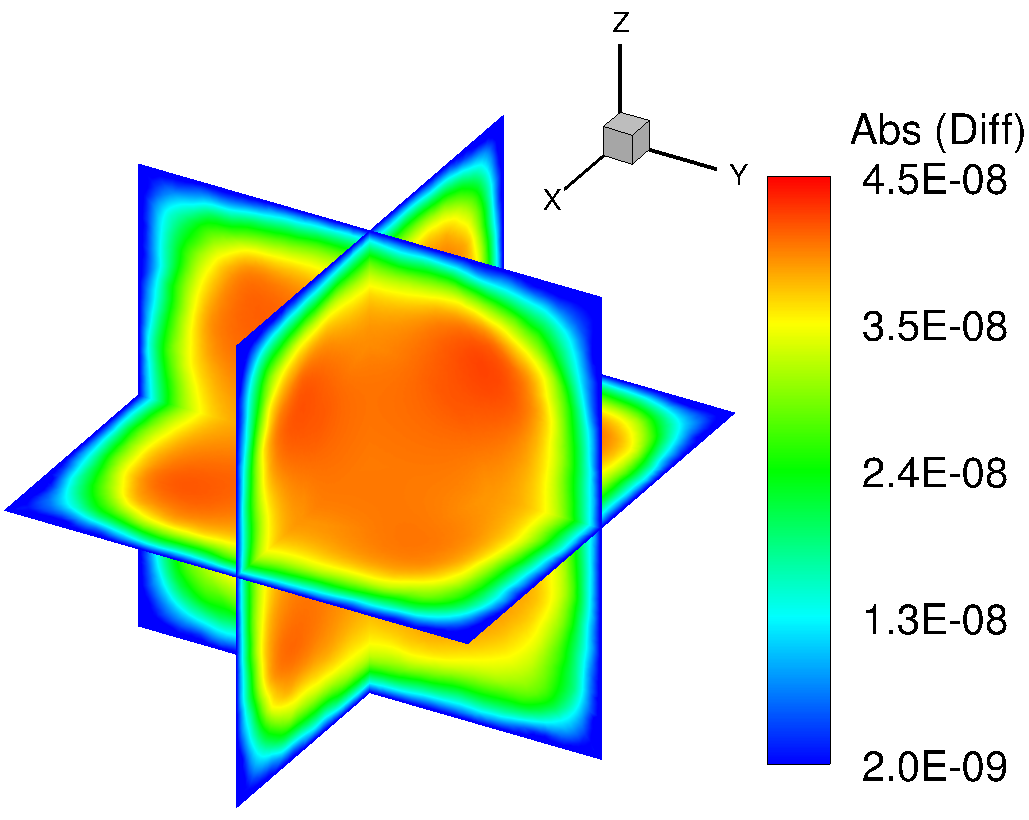}}
	\caption{Sphere inside a cube: Contours of absolute differences between computed and exact solutions for thermal conductivity ratio of 10 and with 711041 scattered points: (a) polynomial degree = 3 and (b) polynomial degree = 6.}
	\label{Fig:sphere_cuboid_three_domains_contours}
\end{figure}

\newpage
\Cref{Fig:sphere_cuboid_error_full_domain} shows the variation of the discretization error with the average grid spacing for thermal conductivity ratios of 5, 10 and 100. In \Cref{Fig:sphere_cuboid_three_domains}, we show the variation for the three subdomains separately for thermal conductivity ratio of 10. We demonstrate that the expected error reduction is also achieved in the three dimensional problems. It can be noted from \cref{fig:sphere_cuboid_interface_error} that at the interface, the order of convergence increases from $1.78$ at polynomial degree of $3$ to $5.66$ at polynomial degree of $6$. \Cref{Fig:sphere_cuboid_three_domains_contours} presents the  difference between numerically computed and manufactured temperature on three orthogonal planes at $X = 0$, $Y = 0$ and $Z = 0$ for two polynomial degrees $3$ and $6$ with 711041 points. The conductivity ratio is 10.

\newpage
\subsection{Convergence plots}
\begin{figure}[H]
	\centering
	\subfigure[]{\includegraphics[width=0.46\textwidth]{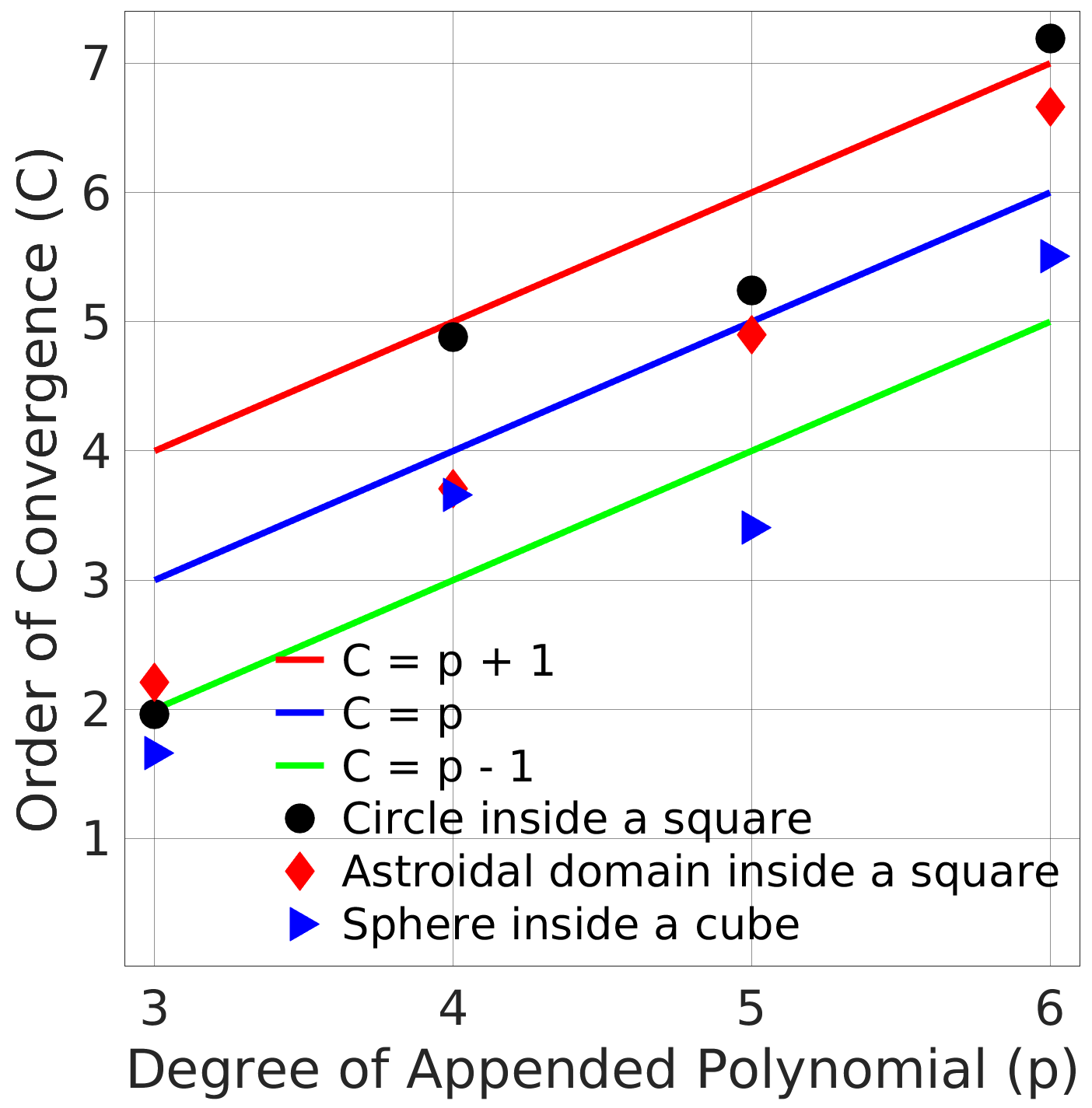}}\hspace{0.3cm}
	\subfigure[]{\includegraphics[width=0.46\textwidth]{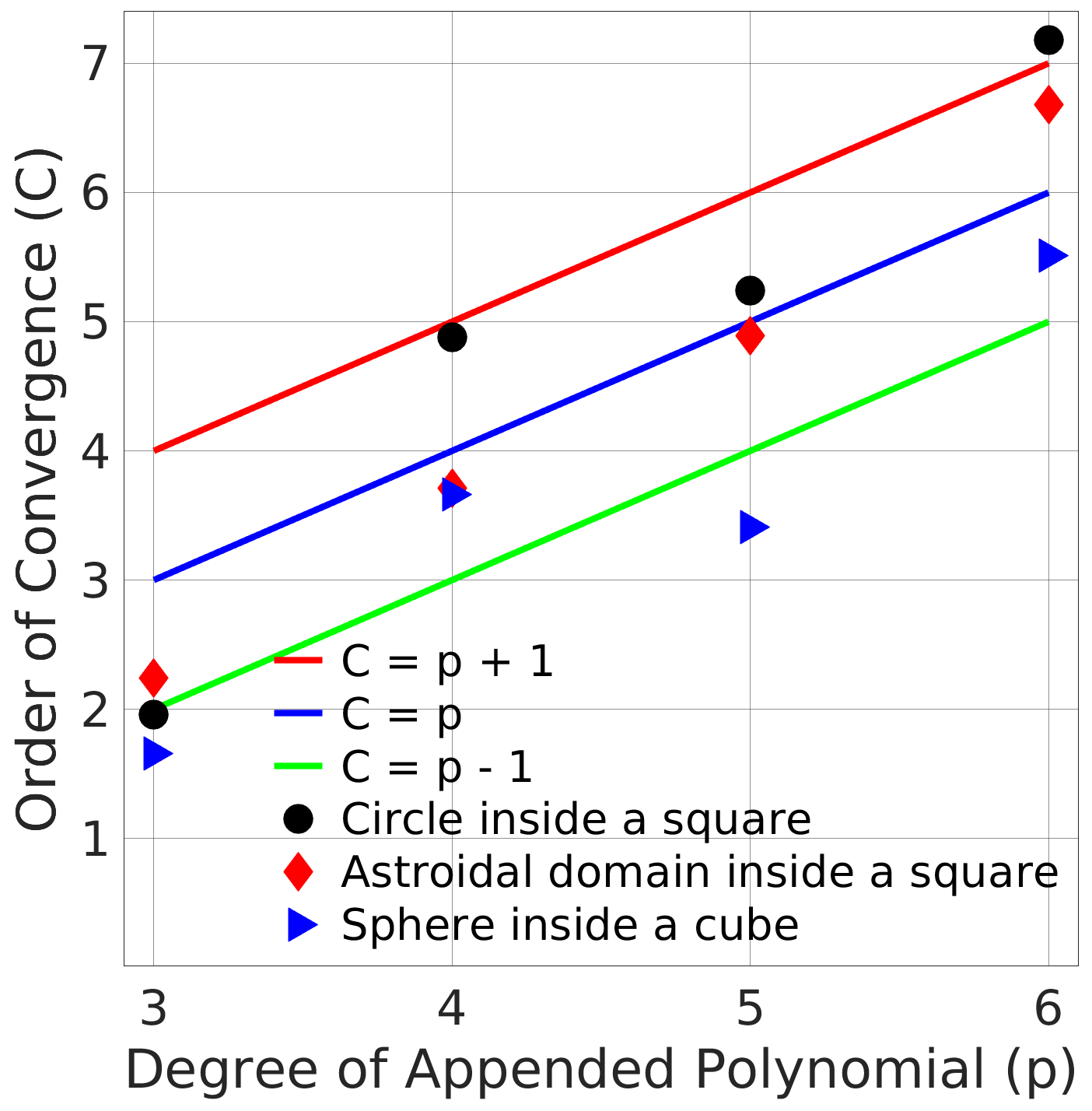}}
	\caption{Convergence plots: Order of convergence (C) vs. Degree of appended polynomial (p) for domain I with thermal conductivity ratios of (a) 10 and (b) 100.}
	\label{Fig:composite_plot_domain_I}
\end{figure}

\begin{figure}[H]
	\centering
	\subfigure[]{\includegraphics[width=0.46\textwidth]{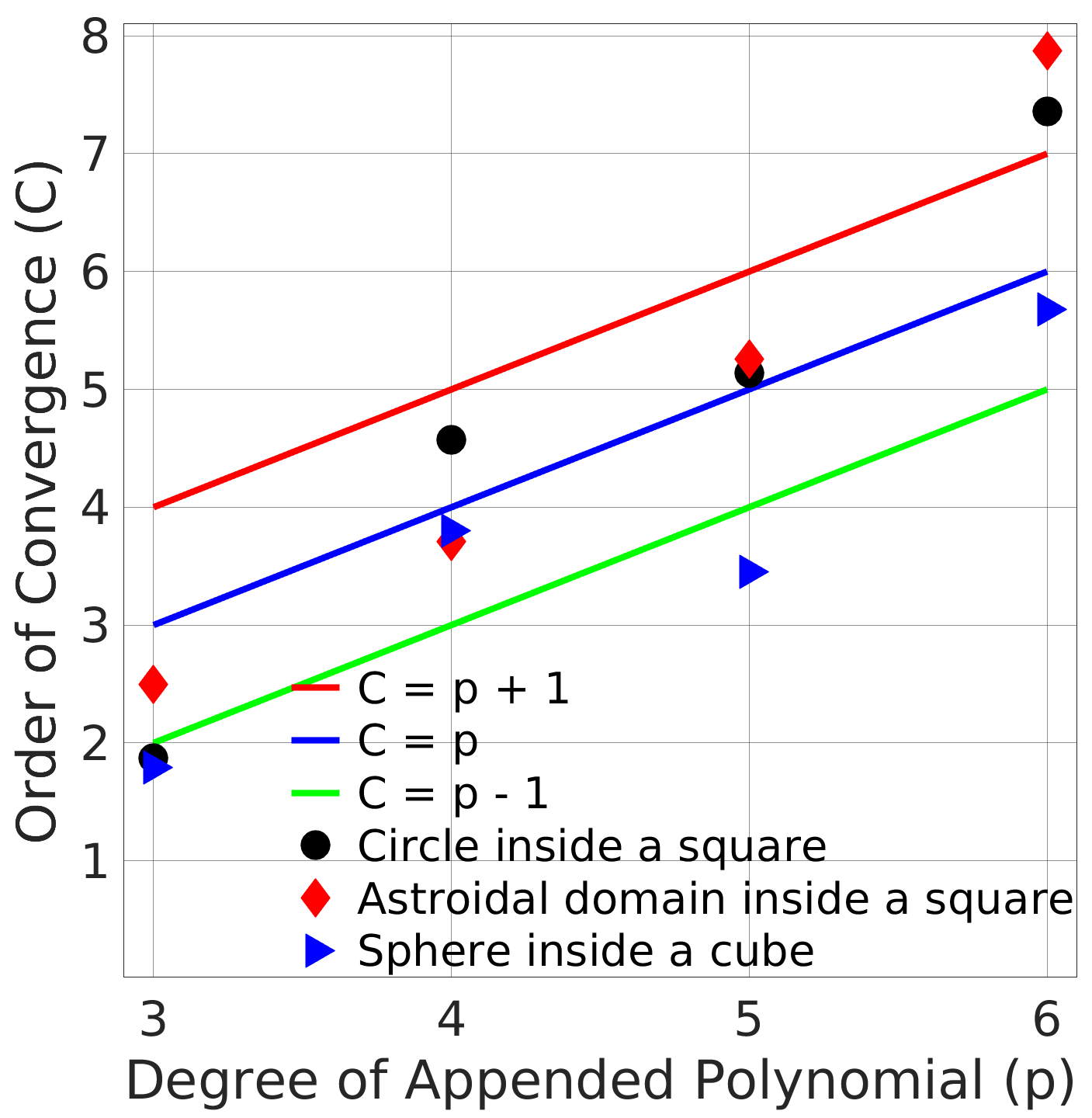}}\hspace{0.3cm}
	\subfigure[]{\includegraphics[width=0.46\textwidth]{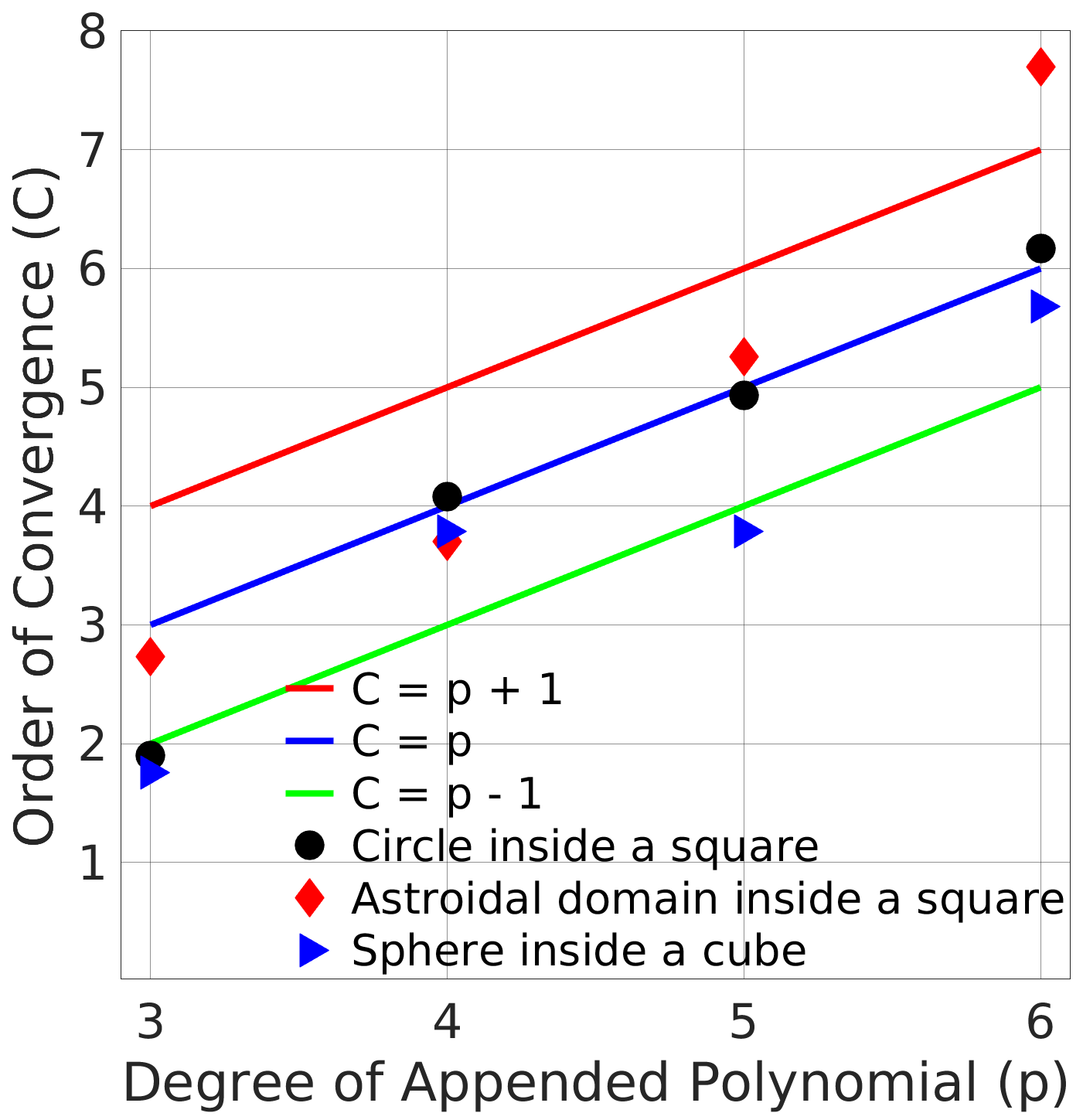}}
	\caption{Convergence plots: Order of convergence (C) vs. Degree of appended polynomial (p) for domain II with thermal conductivity ratios of (a) 10 and (b) 100.}
	\label{Fig:composite_plot_domain_II}
\end{figure}

\begin{figure}[H]
	\centering
	\subfigure[]{\includegraphics[width=0.46\textwidth]{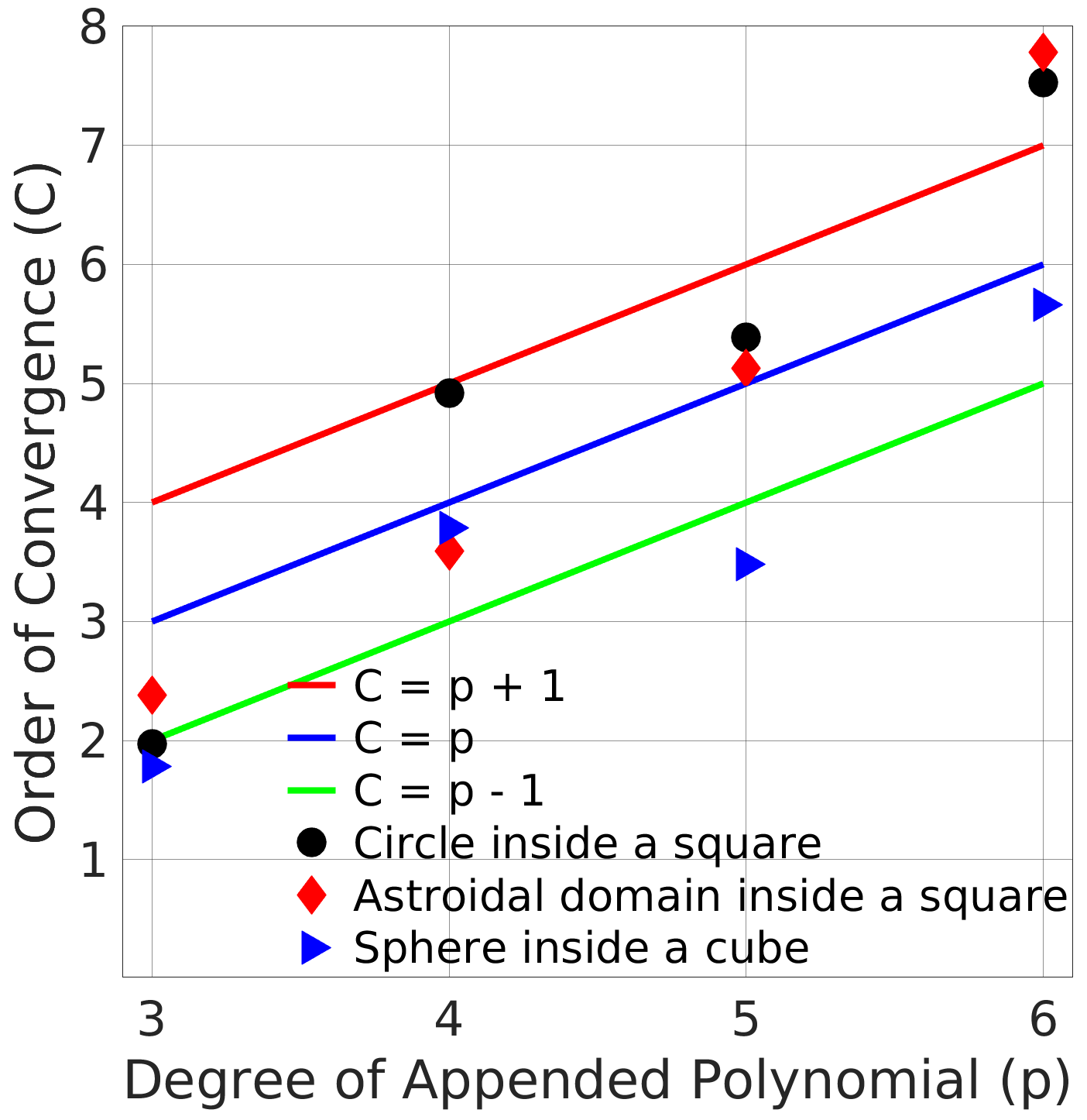}}\hspace{0.3cm}
	\subfigure[]{\includegraphics[width=0.46\textwidth]{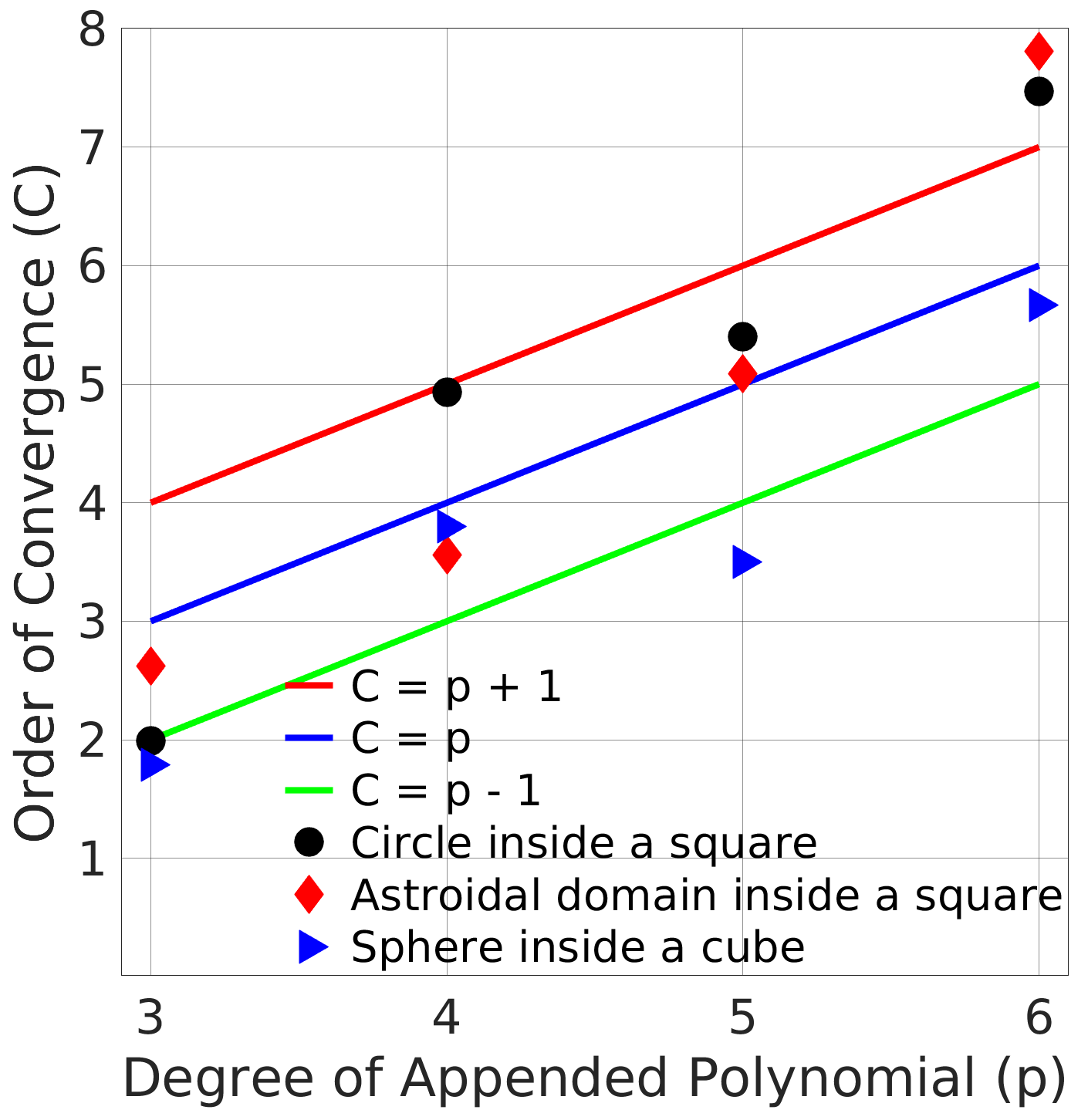}}
	\caption{Convergence plots: Order of convergence (C) vs. Degree of appended polynomial (p) for interface with thermal conductivity ratios of (a) 10 and (b) 100.}
	\label{Fig:composite_plot_interface}
\end{figure}

\begin{figure}[H]
	\centering
	\subfigure[]{\includegraphics[width=0.46\textwidth]{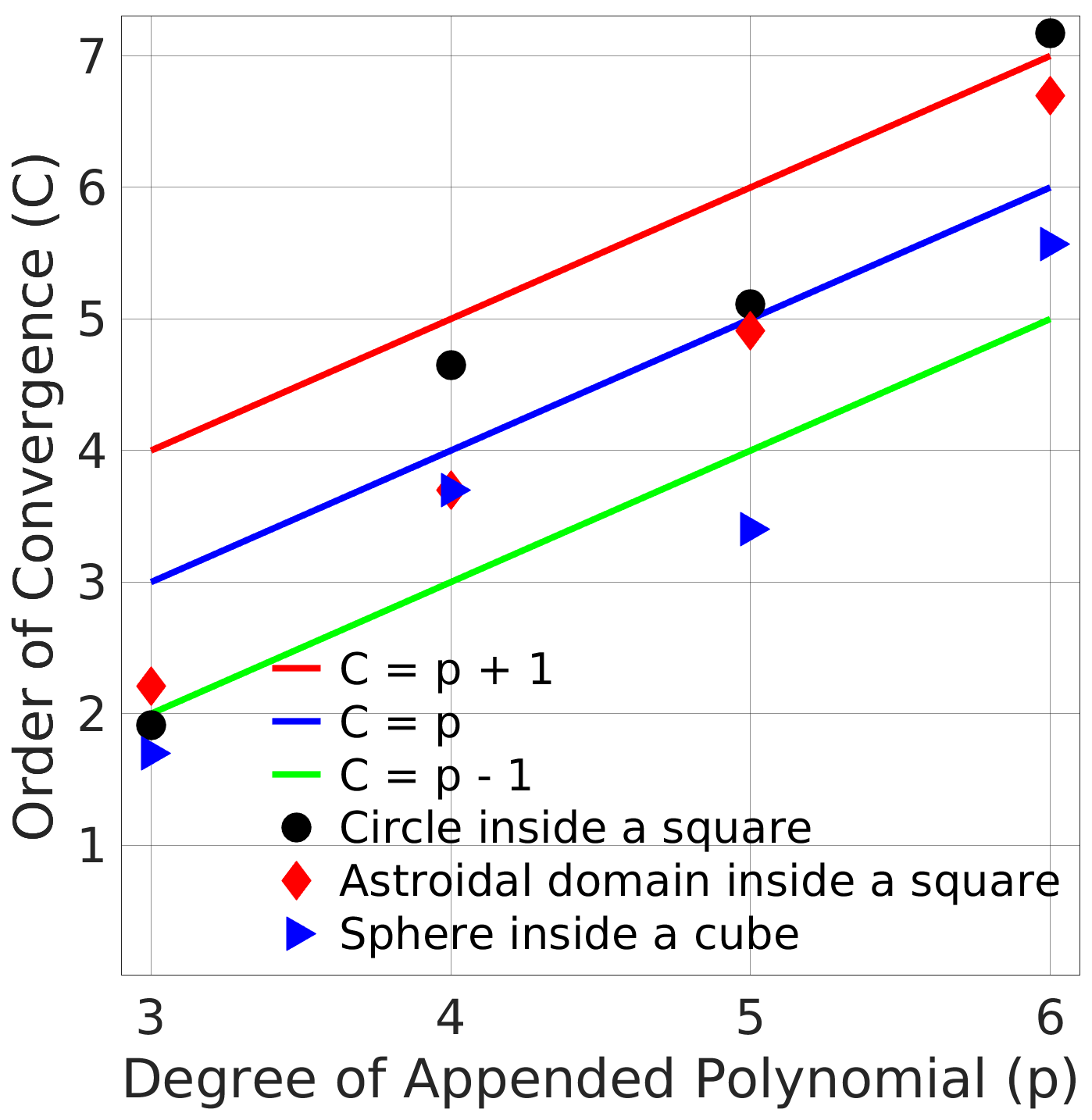}}\hspace{0.3cm}
	\subfigure[]{\includegraphics[width=0.46\textwidth]{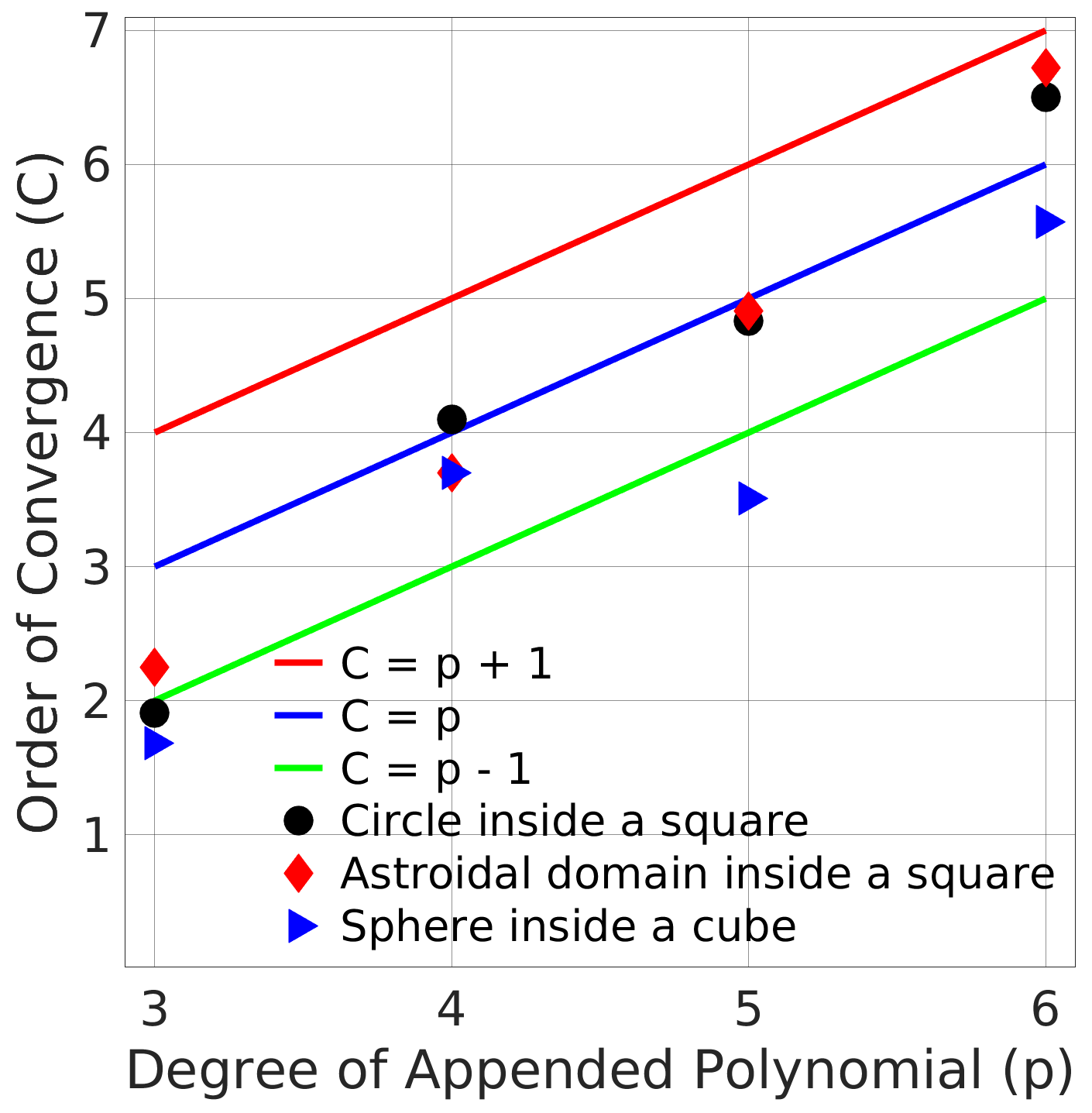}}
	\caption{Convergence plots: Order of convergence (C) vs. Degree of appended polynomial (p) for the entire domain with thermal conductivity ratios of (a) 10 and (b) 100.}
	\label{Fig:composite_plot_full_domain}
\end{figure}

A comprehensive analysis of the order of convergence ($C$) with the degree of appended polynomial ($p$) for the three problems discussed in \cref{Sec:circ_rectangle,Sec:sphere_cuboid,Sec:astroid_rect} is presented in \cref{Fig:composite_plot_full_domain,Fig:composite_plot_domain_I,Fig:composite_plot_domain_II,Fig:composite_plot_interface} for thermal conductivity ratio of 10 and 100. The composite plot of the error is presented individually for each subdomain and the interface, followed by the error for the entire domain. Results are presented for all the polynomials in the same plot. Two lines $C = p - 1$ (green line) and $C = p + 1$ (red line) are also drawn where $p$ is the degree of the appended polynomial. It is observed that the order of convergence lies almost within the regions bounded by $C = p - 1$ and $C = p + 1$ corresponding to a specified degree of appended polynomial. In the current study, the most critical subdomain is the interface. It can be noted from \cref{Fig:composite_plot_interface} that even at the interface, the error decreases with the increase in the degree of the polynomial.

\section{Application to complex problems}
\label{Sec:results}
In this section, we demonstrate the applicability of the multidomain meshless algorithm to complex geometries in two and three dimensions with realistic boundary conditions instead of manufactured solutions. First, we study a couple of two dimensional domains: (i) astroidal domain inside a square and (ii) multiple circular and elliptical subdomains within a square section. Grid independence tests are also performed for the given geometries. Later, a couple of three dimensional geometries: (i) wavy circular tube inside a straight circular tube and (ii) multiple cylindrical fillets inside a cube are simulated to demonstrate the robustness of PHS-RBF algorithm.

\subsection{Astroidal domain inside a square}
\begin{figure}[H]
	\centering
	\subfigure[]{\includegraphics[width=0.47\textwidth]{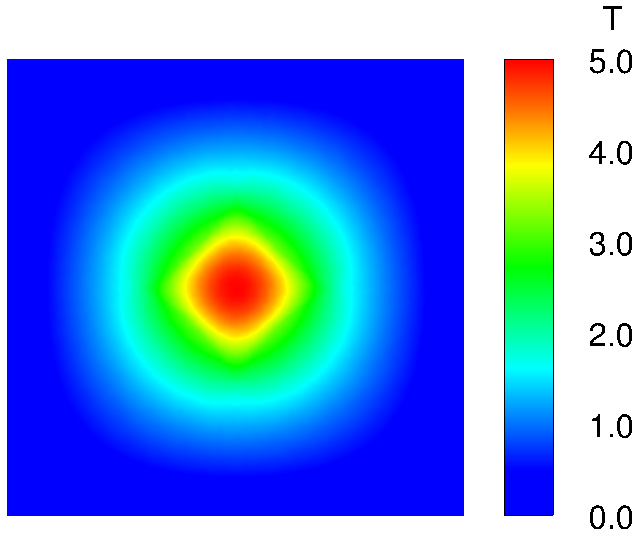}\label{fig:astroid_in_rect_contour_heatgen_100_k_2_1}}\hspace{0.3cm}
	\subfigure[]{\includegraphics[width=0.47\textwidth]{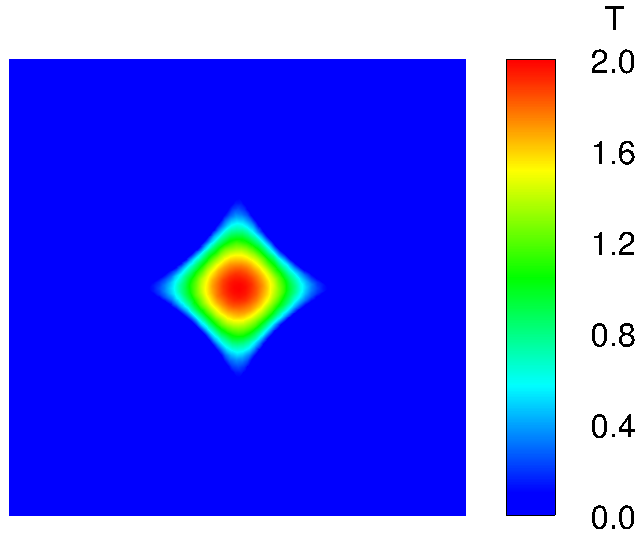}\label{fig:astroid_in_rect_contour_heatgen_100_k_100_1}}
	\caption{Astroidal domain inside a square: Temperature contours: (a) $\frac{k_I}{k_{II}} = 2$ and (b) $\frac{k_I}{k_{II}} = 100$.}
	\label{Fig:astroid_in_rect_contours_dirichlet}
\end{figure}

We have modified the two dimensional case presented earlier (\cref{Sec:astroid_rect}) by imposing $T = 0$ at the outer boundary and prescribing a heat generation rate of 100 units inside the astroidal domain. The conductivity ratio is set to either 2 or 100. The temperature distribution for both the cases are shown in \cref{Fig:astroid_in_rect_contours_dirichlet}. It can be seen that for the first case (\cref{fig:astroid_in_rect_contour_heatgen_100_k_2_1}) that there is a relatively large diffusion of heat as compared to the second case (\cref{fig:astroid_in_rect_contour_heatgen_100_k_100_1}) where the thermal conductivity in domain I is higher i.e. $k_{I} = 100$, hence the temperature gradient is lower. Therefore, the second case shows a sharper jump in temperature across the interface.

\begin{figure}[H]
	\centering
	\subfigure[$k_{I} = 2$, $k_{II} = 1$]{\includegraphics[width=0.48\textwidth]{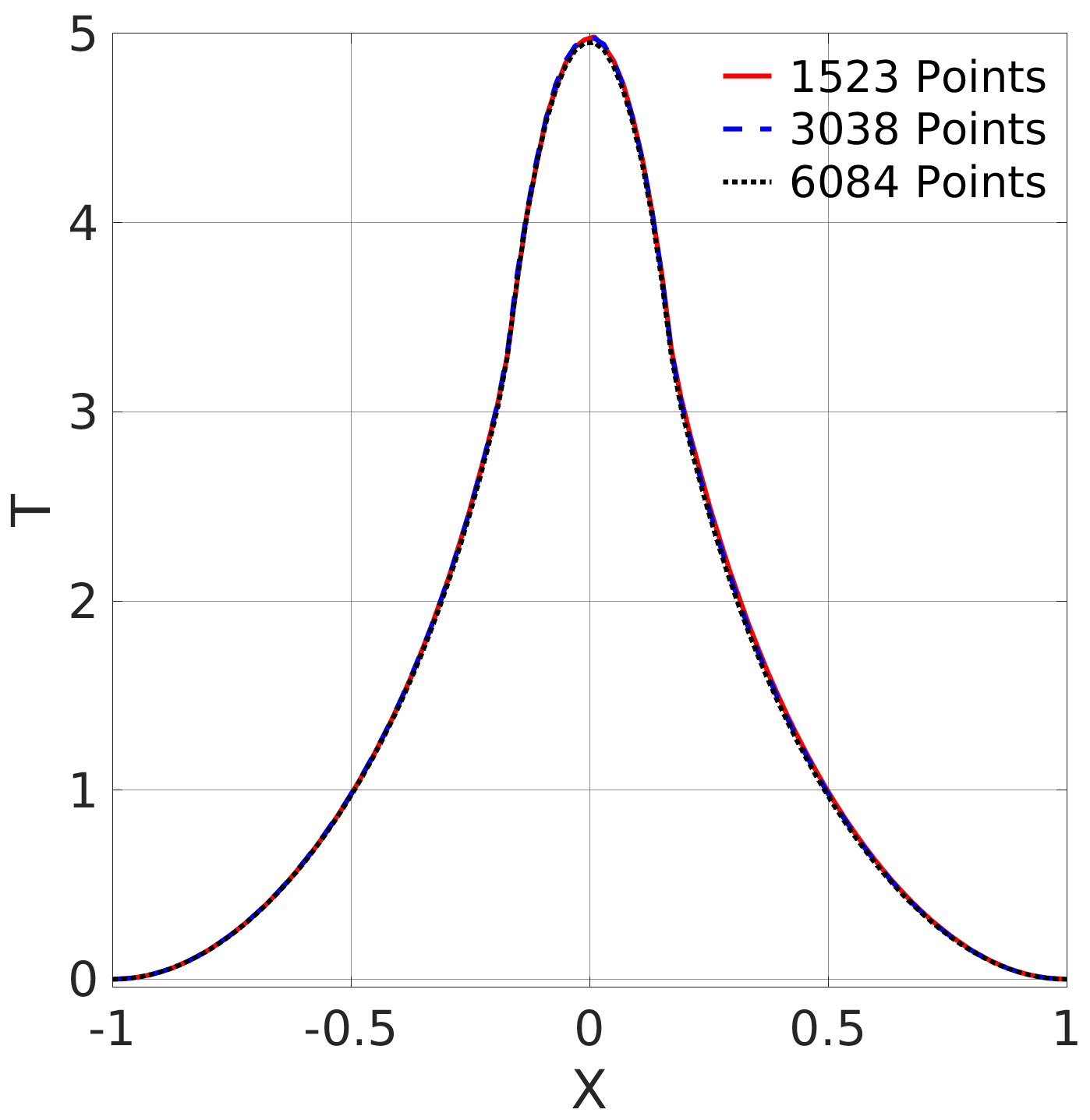}\label{fig:astroid_in_rect_grid_independence_heatgen_100_k_2_1}}
	\subfigure[$k_{I} = 100$, $k_{II} = 1$]{\includegraphics[width=0.495\textwidth]{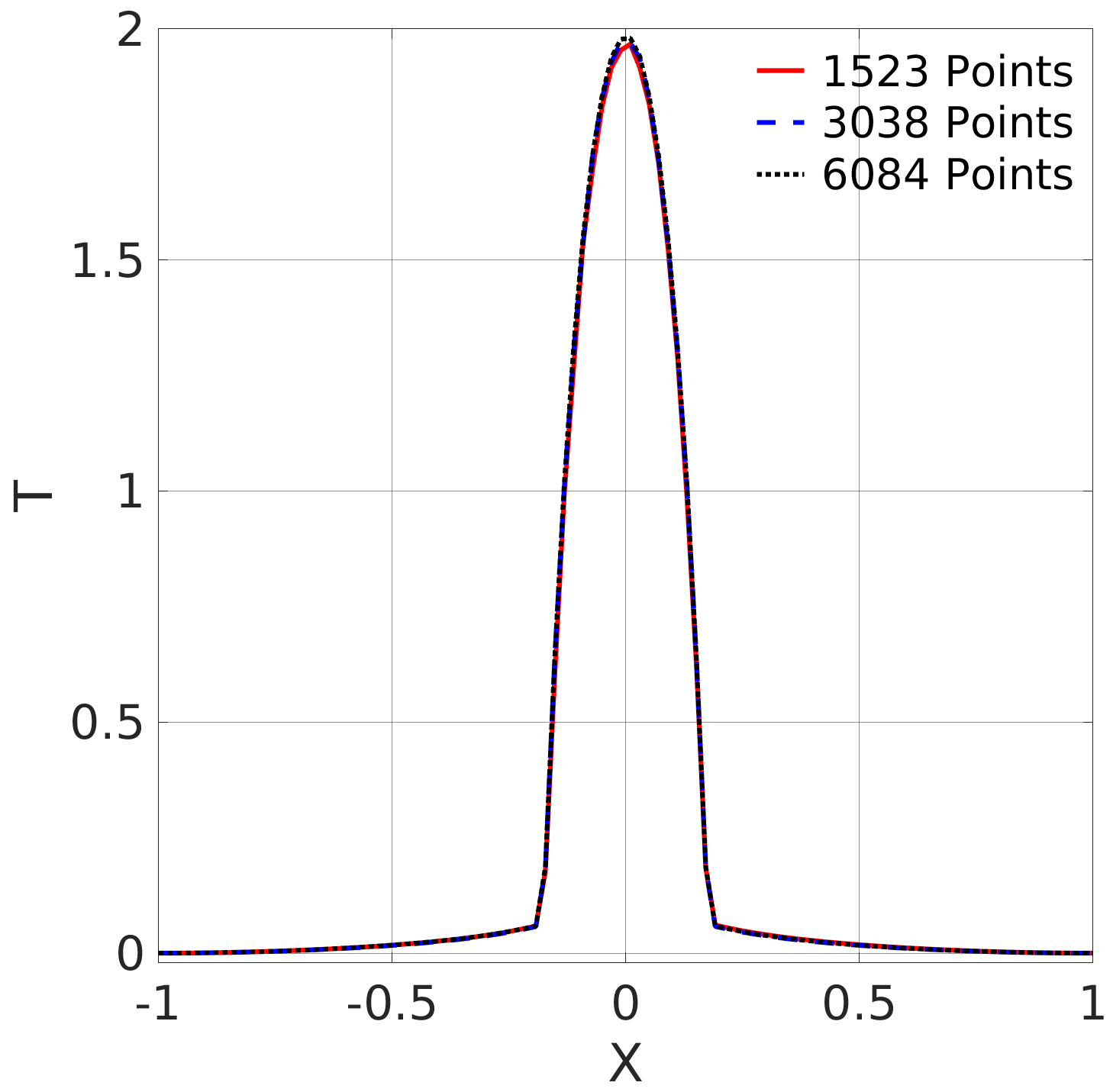}\label{fig:astroid_in_rect_grid_independence_heatgen_100_k_100_1}}
	\caption{Astroidal domain inside a square: Grid independence test for an astroidal domain inside a square.}
	\label{Fig:astroid_in_rect_grid_independence}
\end{figure}

Grid independence tests are also perfomed for both the cases as shown in \cref{fig:astroid_in_rect_grid_independence_heatgen_100_k_2_1,fig:astroid_in_rect_grid_independence_heatgen_100_k_100_1} for the degree of appended polynomial of $6$. The temperature variations along a diagonal of the square are plotted for the two conductivity ratios using up to a total of  6084 scattered points. As the horizontal ($X$) and vertical coordinates ($Y$) of the points along the diagonal are the same, the variation of temperature ($T$) is shown with respect to the horizontal distance ($X$) from the bottom left corner of the square.
\Cref{Table:astroid_in_rect_average_temperature} presents the average temperatures in the solution domain for conductivity ratios of 2 and 100 at polynomial degree of 6. Similar results
are seen for polynomial degrees from 3 to 5.

\begin{table}[H]
	\centering
	\begin{tabular}{|c|c|c|}
		\hline
		Points & \begin{tabular}[c]{@{}c@{}}Average temperature\\ for case (a)\end{tabular} & \begin{tabular}[c]{@{}c@{}}Average temperature\\ for case (b)\end{tabular} \\ \hline
		1523  & 0.966                                                                      & 0.071                                                                      \\ 
		3038  & 1.015                                                                      & 0.078                                                                      \\ 
		6084  & 1.011                                                                      & 0.079                                                                      \\ \hline
	\end{tabular}
	\caption{Astroidal domain inside a square: Average temperatures over the domain for case (a) and case (b).}
		\label{Table:astroid_in_rect_average_temperature}
	\end{table}
	
\subsection{A square with  embedded circular and elliptical domains}
\begin{figure}[H]
	\centering
	\includegraphics[width = 0.7\textwidth]{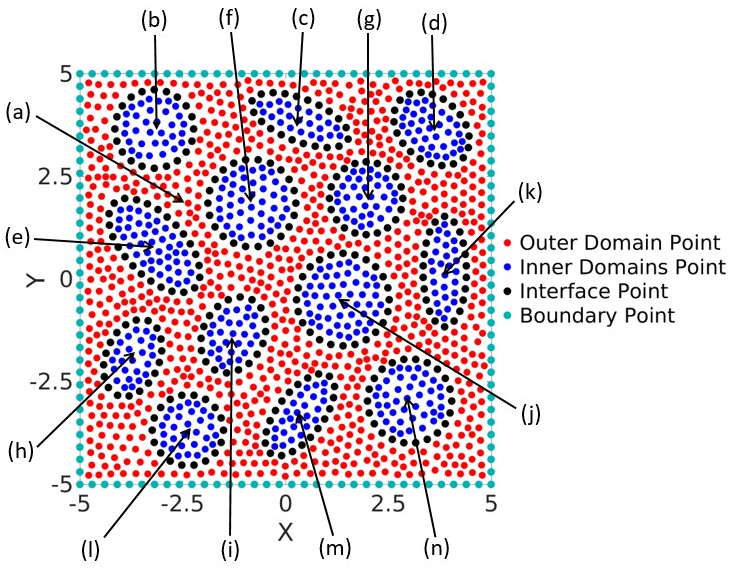}
	\caption{A square with  embedded circular and elliptical domains: Layout of 1574 scattered points.}
	\label{Fig:multiple_subdomain_nodes}
\end{figure}
We next consider a complex domain consisting of circular and elliptical subdomains embedded within a square. A square with a side of 10 units encloses 13 circular and elliptical subdomains of different sizes and aspect ratios. The geometry emulates composite materials having varying thermo-physical characteristics. \Cref{Fig:multiple_subdomain_nodes} shows the layout of 1574 points within the domain divided into multiple subdomains. The outer region is marked by red dots and the different subdomains resolving the ellipses and circles are  discretized by the blue coloured points. Black dots represent the interfaces of different subdomains with the outer domain. 

Within each subdomain, there is either a heat source or sink as tabulated in \cref{Table:multiple_subdomains_interior_conditions}. The external boundary is maintained at a temperature T = 0. \Cref{Table:multiple_subdomains_interior_conditions}  also lists the total number of scattered points distributed within each subdomain and on the corresponding interface. The total number of points in the entire domain is 15032.  

\begin{table}[H]
	\resizebox{\textwidth}{!}{%
		\begin{tabular}{|c|c|c|c|c|}
			\hline
			\multirow{2}{*}{Subdomain} & \multirow{2}{*}{$k$} & \multirow{2}{*}{$\dot{q}$} & \multirow{2}{*}{No. of points} & \multirow{2}{*}{\begin{tabular}[c]{@{}c@{}}No. of points on the\\ interface or boundary\end{tabular}} \\
			&  &  &  &  \\ \hline
			a & 100 & 0 & 8954 & 444 \\ 
			b & 1 & 100 & 408 & 68 \\ 
			c & 3 & -300 & 270 & 64 \\ 
			d & 5 & 500 & 349 & 64 \\ 
			e & 1 & -100 & 453 & 78 \\ 
			f & 3 & 300 & 482 & 74 \\ 
			g & 5 & -500 & 334 & 62 \\ 
			h & 1 & 100 & 279 & 58 \\ 
			i & 3 & -300 & 312 & 60 \\ 
			j & 5 & 500 & 541 & 78 \\ 
			k & 1 & -100 & 292 & 68 \\ 
			l & 3 & 300 & 316 & 60 \\ 
			m & 5 & -500 & 273 & 62 \\ 
			n & 1 & 100 & 457 & 72 \\ \hline
		\end{tabular}%
	}
	\caption{A square with  embedded circular and elliptical domains: Thermal conductivity ($k$), rate of heat generation ($\dot{q}$), number of interior points and interface points for the various domains.}
	\label{Table:multiple_subdomains_interior_conditions}
\end{table}

\begin{figure}[H]
	\centering
	\includegraphics[width = 0.6\textwidth]{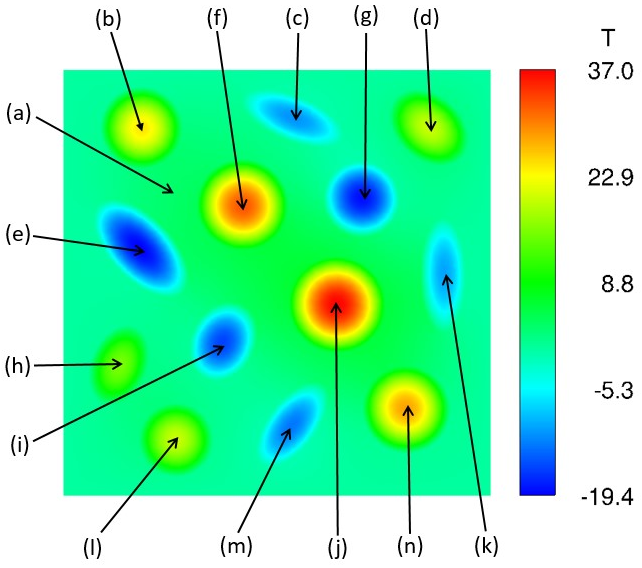}
	\caption{A square with  embedded circular and elliptical domains: Temperature contours in a square with heat generating elliptical and circular subdomains.}
	\label{Fig:multiple_subdomains_contour}
\end{figure}

\Cref{Fig:multiple_subdomains_contour} shows the temperature contours for the given boundary and interior conditions for a polynomial degree of $6$. Due to the high thermal conductivity prescribed in the outer domain, the temperature gradients in this region are relatively small. The temperature distributions in the interior subdomains are according to the heat generation prescribed. The contours of temperature display expected qualitative trends.

\begin{figure}[H]
	\centering
	\subfigure[]{\includegraphics[width=0.48\textwidth]{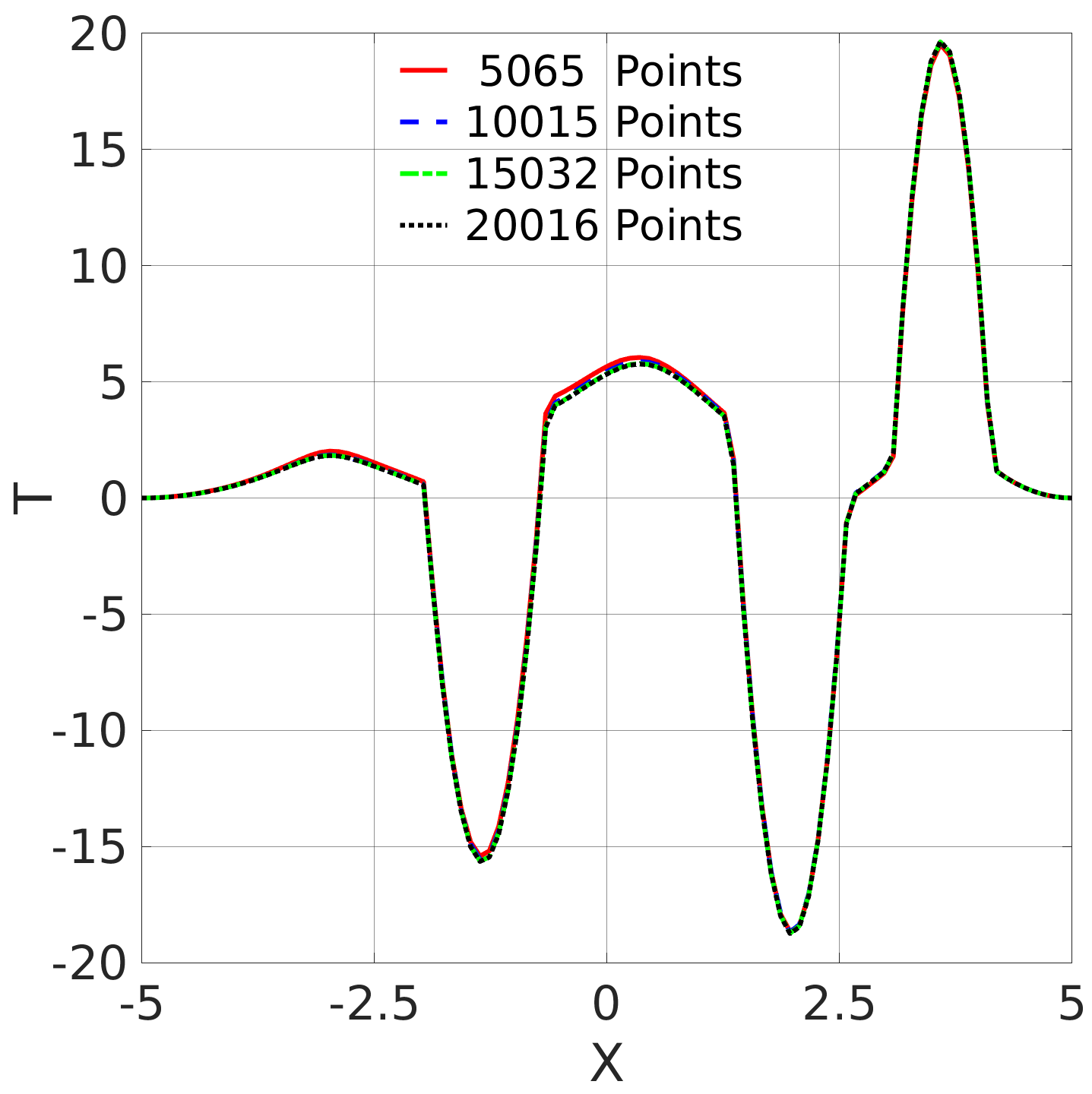}}\hspace{0.3cm}
	\subfigure[]{\includegraphics[width=0.48\textwidth]{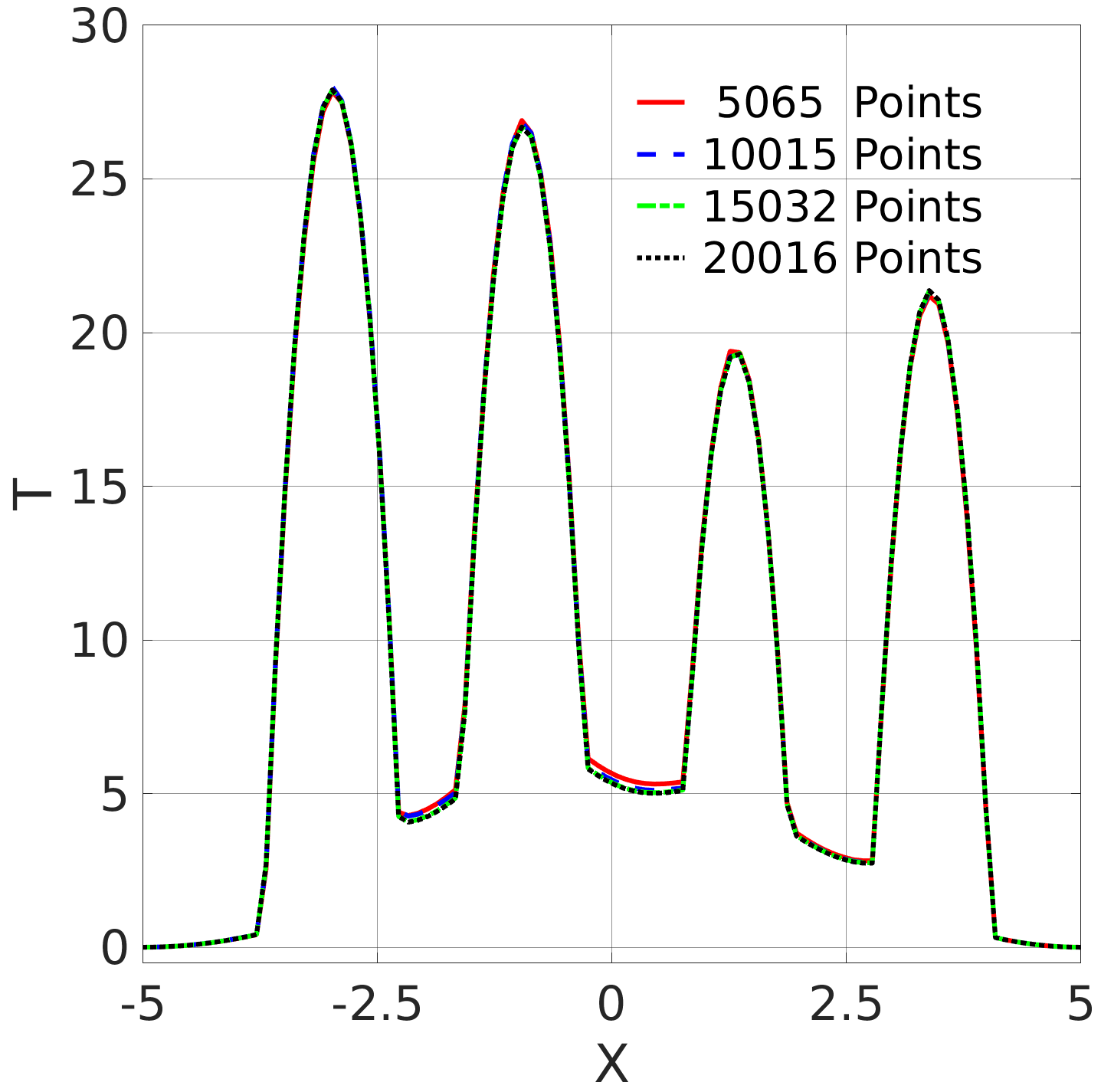}}
	\caption{A square with  embedded circular and elliptical domains: Temperature profiles along the two cross-diagonals for different point distributions.}
	\label{Fig:multiple_subdomains_grid_independence}
\end{figure}

Grid independency tests are also performed with four different point layouts consisting of 5065, 10015, 15032 and 20016 points for the degree of appended polynomial of $6$. \Cref{Fig:multiple_subdomains_grid_independence} shows the temperature variation on the diagonals from $(-5,-5)$ to $(5,5)$ and $(-5,5)$ to $(5,-5)$. The variation is shown with respect to the $X$ coordinates of the points along the diagonals. It can be seen that the temperature distributions overlap for the different point sets, implying fast convergence of the solutions. \Cref{Table:multiple_subdomains_average_temperature} gives the convergence of the average temperatures over the entire domain.  

\begin{table}[H]
	\centering
	\begin{tabular}{|c|c|}
		\hline
		Number of points & Average temperature \\ \hline
		5065  & 2.587               \\ 
		10015 & 2.744               \\ 
		15032 & 2.706               \\ 
		20016 & 2.706               \\ \hline
	\end{tabular}
	\caption{A square with  embedded circular and elliptical domains: Average temperatures for different point sets.}
	\label{Table:multiple_subdomains_average_temperature}
\end{table}

\newpage
\subsection{Wavy circular tube inside a straight circular tube}
\begin{figure}[H]
	\centering
	\includegraphics[width = 0.6\textwidth]{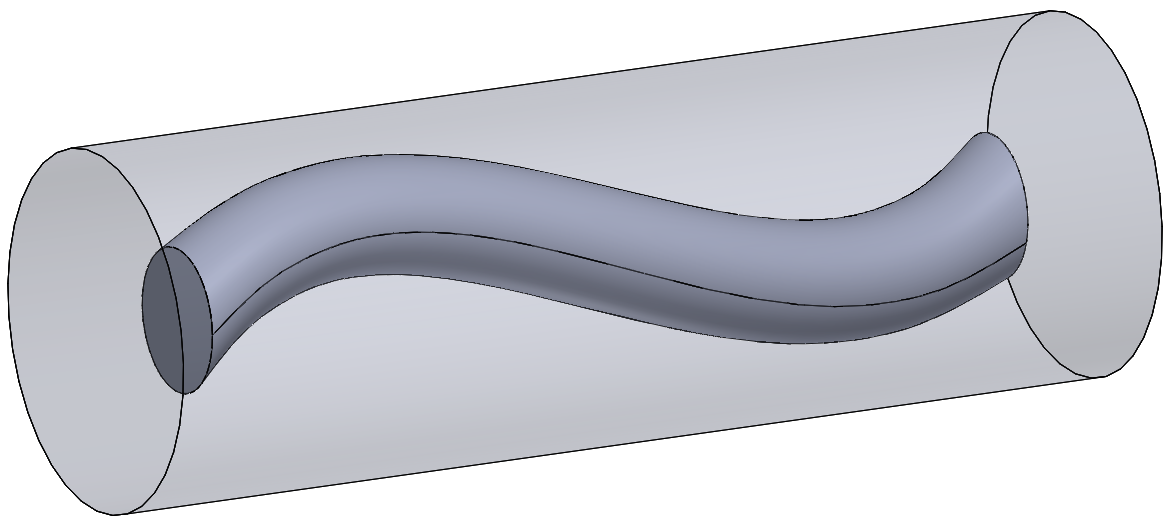}
	\caption{Wavy circular tube inside a straight circular tube.}
	\label{Fig:wavy_circular_tube}
\end{figure}

\begin{figure}[H]
	\centering
	\subfigure[]{\includegraphics[width=0.48\textwidth]{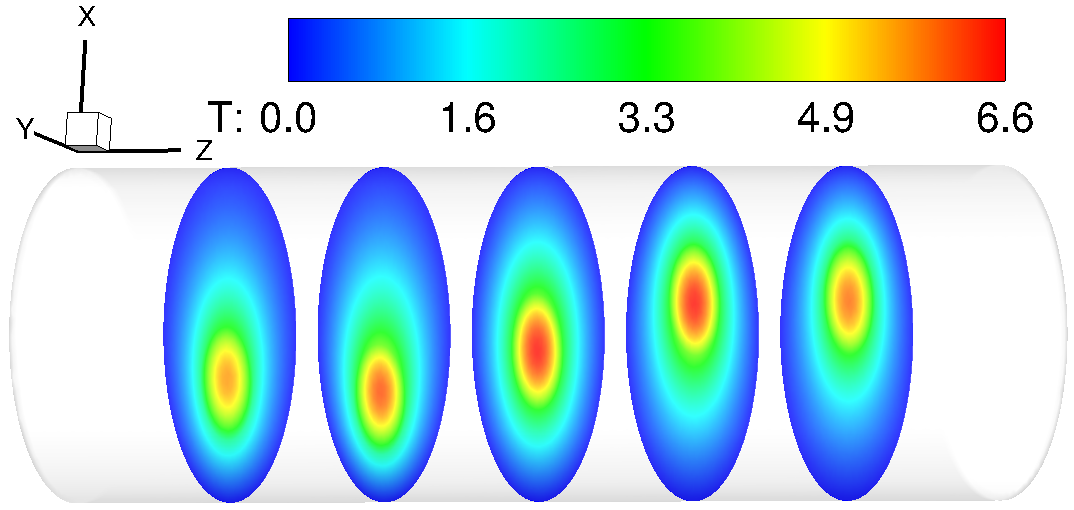}\label{fig:wavy_tube_temperature_contours_part_I}}\hspace{0.3cm}
	\subfigure[]{\includegraphics[width=0.48\textwidth]{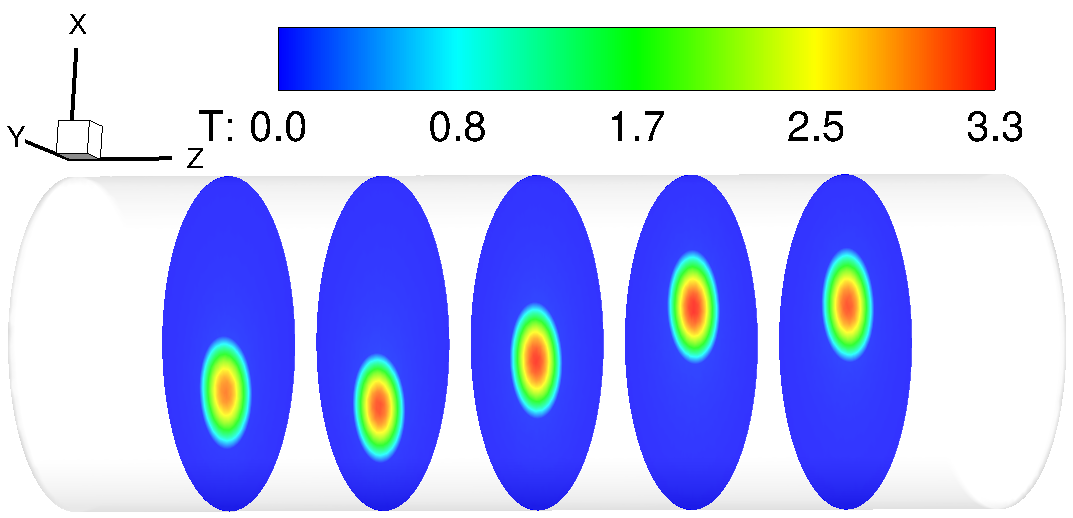}\label{fig:wavy_tube_temperature_contours_part_II}}
	\caption{Wavy circular tube inside a straight circular tube: Temperature contours: (a) $\frac{k_{I}}{k_{II}} = 2$ and (b) $\frac{k_{I}}{k_{II}} = 100$.}
	\label{Fig:wavy_tube_temperature_contours}
\end{figure}

\begin{figure}[H]
	\centering
	\subfigure[]{\includegraphics[width=0.48\textwidth]{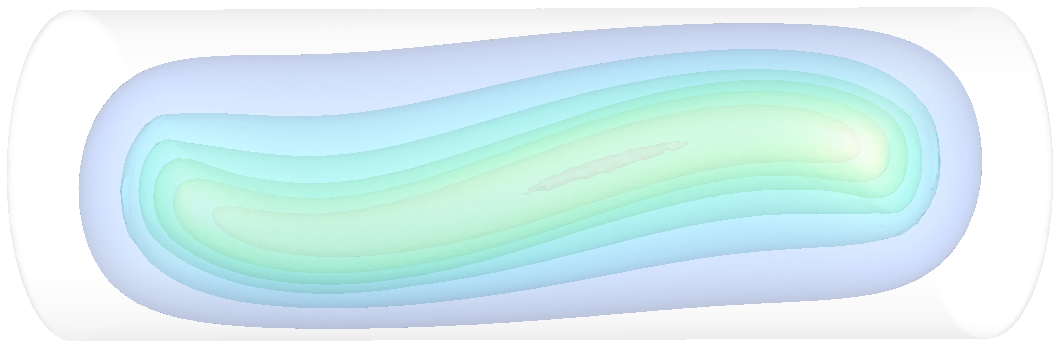}\label{fig:wavy_tube_isosurface_contours_part_I}}\hspace{0.3cm}
	\subfigure[]{\includegraphics[width=0.48\textwidth]{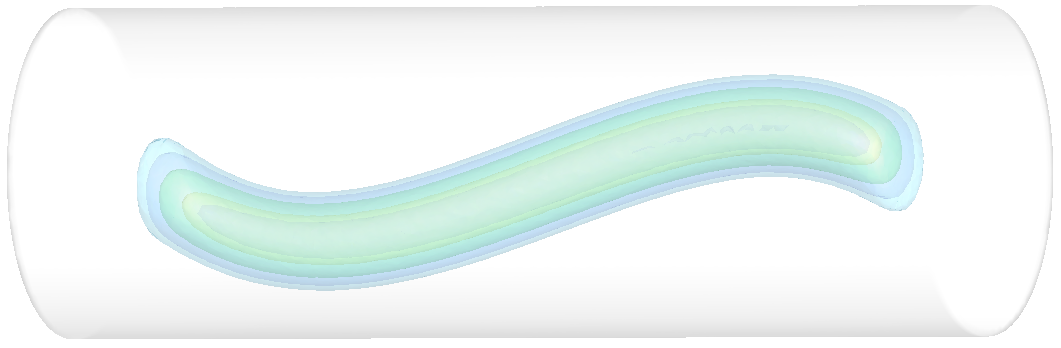}\label{fig:wavy_tube_isosurface_contours_part_II}}
	\caption{Wavy circular tube inside a straight circular tube: Isosurfaces: (a) $\frac{k_{I}}{k_{II}} = 2$ and (b) $\frac{k_{I}}{k_{II}} = 100$.}
	\label{Fig:wavy_tube_isosurface_contours}
\end{figure}

We next consider two three-dimensional geometries to demonstrate the applicability of the method to complex three-dimensional problems. \Cref{Fig:wavy_circular_tube} shows a wavy circular tube of radius 0.4 unit and length of 5 units confined in a straight circular tube of radius 1 unit and length of 6 units. The region outside the wavy circular tube is marked as domain I, the region within the wavy tube is marked as domain II, and the surface of the tube is marked as the interface. The wavy circular tube has an internal heat generation rate of 100 units. 
Two test cases with diffusion coefficient ratios of $\frac{k_{I}}{k_{II}}$ = 2 and 100 are considered. \Cref{Fig:wavy_tube_temperature_contours} shows the temperature contours for both the cases at five planes orthogonal to $Z$ axis i.e. at $Z$ = 1, 2, 3, 4 and 5 units with 100220 scattered points distributed throughout the domain. The temperature contours are more diffused in \cref{fig:wavy_tube_temperature_contours_part_I} as compared to \cref{fig:wavy_tube_temperature_contours_part_II} due to lower thermal conductivity ratio \Big($\frac{k_{I}}{k_{II}}$\Big), hence relatively higher values of temperature gradient in domain I. \Cref{Fig:wavy_tube_isosurface_contours} shows the isosurfaces of  temperature for both the cases.

\subsection{A complex dispersion of cylindrical fillets inside a cube}
\begin{figure}[H]
	\centering
	\includegraphics[width = 0.5\textwidth]{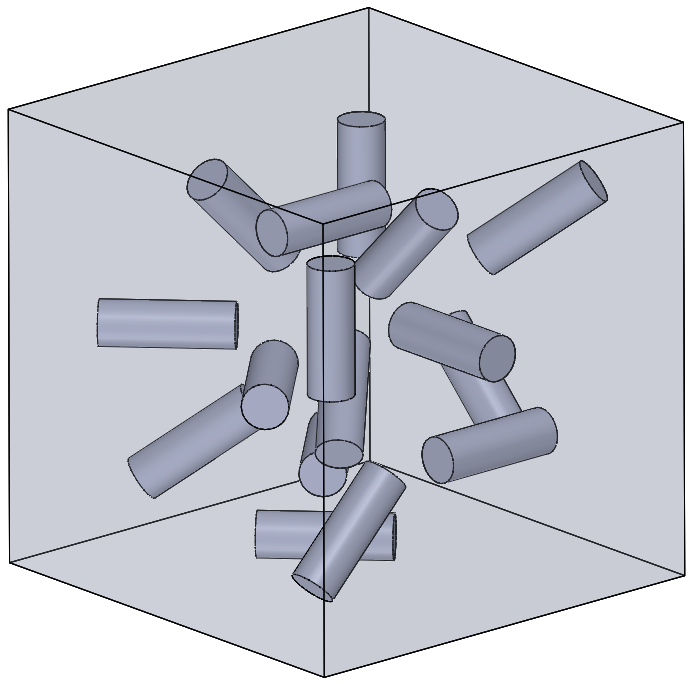}
	\caption{Geometry for multiple cylindrical fillets inside a cube.}
	\label{Fig:multiple_cylindrical_fillets}
\end{figure}

The last three dimensional problem consists of 16 cylindrical fillets of radius = 0.3 unit, length = 1.73 units within a cube of dimensions 6 units $\times$ 6 units $\times$ 6 units as shown in \cref{Fig:multiple_cylindrical_fillets}. The cylindrical fillets are arranged in such a way that adequate space exists between the fillets to place sufficient number of points. 

\begin{figure}[H]
	\centering
	\includegraphics[width = 0.7\textwidth]{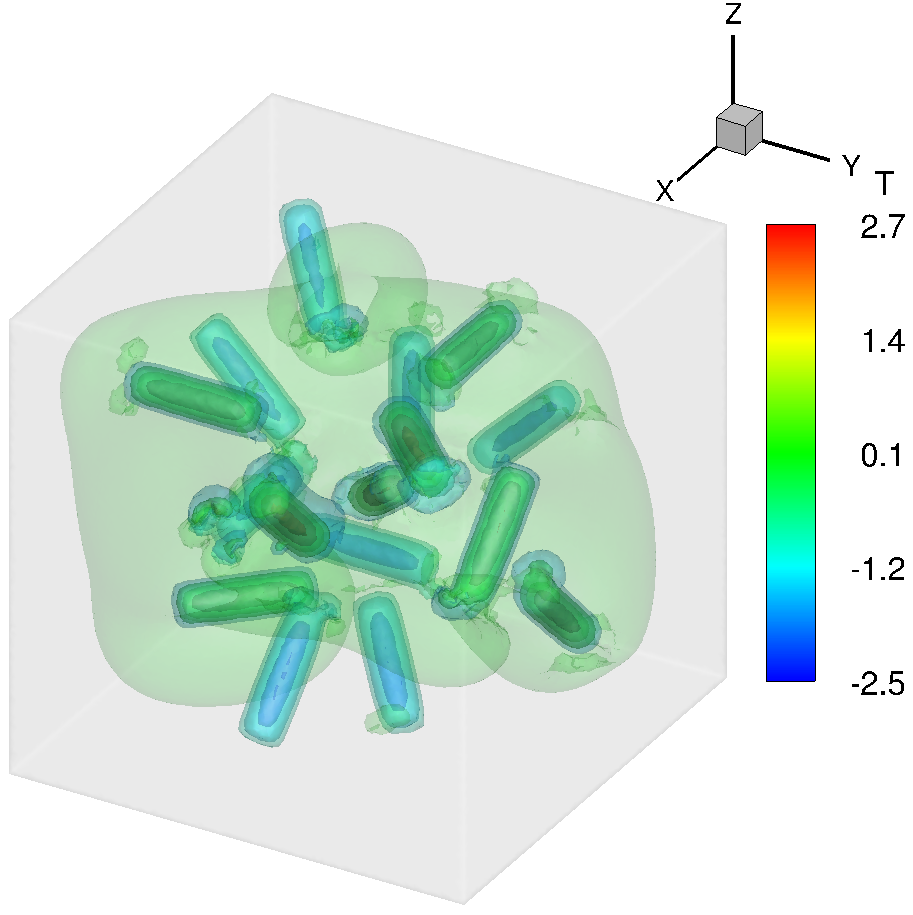}
	\caption{A complex dispersion of cylindrical fillets inside a cube: Three-dimensional distributions of temperature for the fillet geometry.}
	\label{Fig:multiple_cylindrical_fillets_contour}
\end{figure}

The cylindrical fillets are assigned to be either a heat source or a heat sink with different thermal conductivities. The values of the diffusion coefficients have been assigned to be higher in the non-fillet regions than inside the fillets. A total of 198545 points are used to discretize the three-dimensional region using vertices of a finite element grid. A polynomial degree of 6 is used for the PHS-RBF. The outer boundary is assigned a temperature of zero on all six faces. A steady state temperature distribution is computed and the isosurfaces of the same  is  shown in \cref{Fig:multiple_cylindrical_fillets_contour}. 

\section{Summary and Future work}
\label{Sec:summary}
In this paper, we have described a multidomain cloud based numerical algorithm for computing heat transfer in complex domains with significant variations of thermo-physical properties. The algorithm interpolates variables at scattered points using polyharmonic splines and appended polynomial. At the interfaces between subdomains, the fluxes are balanced, and the clouds of interpolation points are restricted to be within individual subdomains. We show that the imposition of the flux balance condition at the interface points (and not satisfying the governing equation) does not degrade the expected high accuracy of the overall algorithm. Several two and three-dimensional simulations with manufactured solutions have been performed for different conductivity ratios, degree of the appended polynomial and inter-point spacings. The expected high accuracy is observed for all the problems. Further, a couple of two and three dimensional problems in complex domains are simulated to demonstrate the capability of the algorithm. The proposed multidomain algorithm has the potential to study fluid flow and convective heat transfer in complex domains, including natural convection, phase change and multiphysics phenomena such as coupling with electromagnetics, radiative heat transfer, etc. These will be considered in the future. 

\section*{Acknowledgement}
Naman Bartwal would like to thank his labmate Siddharth D. Sharma for fruitful discussions on modeling of geometries. The authors acknowledge SPARC (Scheme for Promotion of Academic and Research Collaboration $\rightarrow$ Project code: P1211) for its financial support.

\section*{Declaration of Interests}
The authors declare that they have no known competing financial interests or personal relationships that could have appeared to influence the work reported in this paper.
\newpage
\printnomenclature

\section*{References}
\bibliography{mybibfile}

\end{document}